\documentclass[]{amsart}

\usepackage[Symbol]{upgreek}

\usepackage{graphicx}

\usepackage{amssymb}
\usepackage{epic}
\usepackage{mathrsfs}
\usepackage{accents}
\usepackage{array}
\usepackage{stmaryrd}
\usepackage{enumitem}
\usepackage{longtable}

\usepackage[hidelinks]{hyperref}
\usepackage{etoolbox} 

\usepackage{tikz}

\input{xy}
\xyoption{matrix}
\xyoption{arrow}

\numberwithin{equation}{section}

\makeatletter
\renewcommand{\tocsection}[3]
 { \indentlabel{\@ifnotempty{#2}{\parbox{2.5em}{\ignorespaces#1 #2.}\quad}}#3}

\makeindex

\newcommand{\bbold}{\mathbb}
\newcommand{\cal}{\mathcal}

\def\R { {\bbold R} }
\def\Q { {\bbold Q} }
\def\Z { {\bbold Z} }
\def\C { {\bbold C} }
\def\N { {\bbold N} }

\def\c {\mathcal{C}}

\def \I{\operatorname{I}}

\def \Ex{\operatorname{E}}
\def \Dx{\operatorname{D}}
\def \Exp{\operatorname{E}}

\def \ex{\operatorname{e}}
\def \wr {\operatorname{wr}}

\def \Frac {\operatorname{Frac}}

\renewcommand\epsilon{\varepsilon}

\def \d{\operatorname{d}}

\def \ev{\operatorname{e}}
\def \bar {\overline}
\def \<{\langle}
\def \>{\rangle}

\def \tilde {\widetilde}

\def \((  {(\!(}
\def \)) {)\!)}

\def \Li{\operatorname{Li}}
\def \res{\operatorname{res}}

\def \k {{{\boldsymbol{k}}}}
\def \flatter{\mathrel{\prec\!\!\!\prec}}
\DeclareMathSymbol{\precequ}{\mathrel}{symbols}{"16}
\DeclareMathSymbol{\succequ}{\mathrel}{symbols}{"17}
\def \flattereq{\mathrel{\precequ\!\!\!\precequ}}
\def \steeper{\mathrel{\succ\!\!\!\succ}}

\def \comp{\mathrel{-{\hskip0.06em\!\!\!\!\!\asymp}}}

\def \nasymp{\not\asymp}

\renewcommand{\Re}{\operatorname{Re}}
\renewcommand{\Im}{\operatorname{Im}}

\newtheorem{theorem}{Theorem}[section]

\newtheorem{lemma}[theorem]{Lemma}
\newtheorem{prop}[theorem]{Proposition}
\newtheorem{cor}[theorem]{Corollary}

\newtheorem{corintro}{Corollary}

\newtheorem*{theoremUnnumbered}{Theorem}

\theoremstyle{definition}

\theoremstyle{remark}
\newtheorem*{example}{Example}

\newtheorem{exampleNumbered}[theorem]{Example}

\newtheorem*{notation}{Notation}

\newtheorem*{remarks}{Remarks}
\newtheorem*{remark}{Remark}

\newtheorem*{conjecture}{Conjecture}
\newtheorem{remarkNumbered}[theorem]{Remark}

\newcommand{\abs}[1]{\lvert#1\rvert}
\newcommand{\dabs}[1]{\lVert#1\rVert}

\def \sgn {\operatorname{sign}}

\let\oldi\i
\let\oldj\j
\renewcommand\i{\relax\ifmmode{\boldsymbol{i}}\else\oldi\fi}
\renewcommand\j{\relax\ifmmode{\boldsymbol{j}}\else\oldj\fi}

\renewcommand\leq{\leqslant}
\renewcommand\geq{\geqslant}
\renewcommand\preceq{\preccurlyeq}
\renewcommand\succeq{\succcurlyeq}
\renewcommand\le{\leq}
\renewcommand\ge{\geq}

\DeclareMathAlphabet{\mathbf}{OML}{cmm}{b}{it}

\DeclareFontFamily{U}{fsy}{}
\DeclareFontShape{U}{fsy}{m}{n}{<->s*[.9]psyr}{}
\DeclareSymbolFont{der@m}{U}{fsy}{m}{n}
\DeclareMathSymbol{\der}{\mathord}{der@m}{182}

\DeclareSymbolFont{der@m}{U}{fsy}{m}{n}
\DeclareMathSymbol{\derdelta}{\mathord}{der@m}{100}

\newcommand\nwt{\operatorname{nwt}}
\newcommand\ndeg{\operatorname{ndeg}}

\DeclareSymbolFont{imag@m}{OT1}{cmr}{m}{ui}
\DeclareMathSymbol{\imag}{\mathord}{imag@m}{105}

\DeclareFontFamily{OMS}{smallo}{}
\DeclareFontShape{OMS}{smallo}{m}{n}{<->s*[.65]cmsy10}{}
\DeclareSymbolFont{smallo@m}{OMS}{smallo}{m}{n}
\DeclareMathSymbol{\smallo}{\mathord}{smallo@m}{79}

\DeclareFontFamily{OMS}{largerdot}{}
\DeclareFontShape{OMS}{largerdot}{m}{n}{<->s*[.8]cmsy10}{}
\DeclareSymbolFont{largerdot@m}{OMS}{largerdot}{m}{n}
\DeclareMathSymbol{\largerdot}{\mathord}{largerdot@m}{15}

\DeclareMathSymbol{\llambda}{\mathord}{der@m}{108}
\DeclareMathSymbol{\rrho}{\mathord}{der@m}{114}

\def \upg{\upgamma}
\def \Upg{\Upgamma}
\def \upl{\uplambda}
\def \Upl{\Uplambda}
\def \upo{\upomega}
\def \Upo{\Upomega}

\def \Upd{\Updelta}

\def\HLO{\Upl\Upo}

\newcommand{\equationqed}[1]{\[\pushQED{\qed}#1 \qedhere\popQED\]\let\qed\relax}
\newcommand{\alignqed}[1]{\begin{align*}\pushQED{\qed} #1 \qedhere\popQED\end{align*}\let\qed\relax}

\makeatletter
\newcommand{\dminus}{\mathbin{\text{\@dminus}}}

\newcommand{\@dminus}{%
  \ooalign{\hidewidth\raise1ex\hbox{\bf.}\hidewidth\cr$\m@th-$\cr}%
}
\makeatother

\def \Car{\mathcal{C}^r_a}
\def \Caz{\mathcal{C}^0_a}
\def \Cao{\mathcal{C}^1_a}
\def \Cat{\mathcal{C}^2_a}
\def \Cainf{\mathcal{C}^{\infty}_a}
\def \Caom{\mathcal{C}^{\omega}_a}

\def \Gr{\mathcal{C}^r}
\def \Gn{\mathcal{C}^n}
\def \Gz{\mathcal{C}^0}
\def \Go{\mathcal{C}^1}
\def \Gt{\mathcal{C}^2}
\def \Gi{\mathcal{C}^{<\infty}}
\def \Ginf{\mathcal{C}^{\infty}}
\def \Gom{\mathcal{C}^{\omega}}
\def \inv{\operatorname{inv}}
\def \Sol{\operatorname{Sol}}

\def \Car{\mathcal{C}^r_a}
\def \Carl{\mathcal{C}^{r-1}_a}

\def \Caz{\mathcal{C}^0_a}

\def \Cao{\mathcal{C}^1_a}
\def \Cat{\mathcal{C}^2_a}
\def \Cainf{\mathcal{C}^{\infty}_a}

\def \Can{\mathcal{C}^n_a}
\def \Caom{\mathcal{C}^{\omega}_a}
\def \Calinf{\mathcal{C}^{<\infty}}

\def \Caln{\mathcal{C}^{n}}

\makeatletter
\renewcommand\part{\@startsection{part}{0}%
  \z@{\linespacing\@plus\linespacing}{.5\linespacing}%
  {\normalfont\bfseries\centering}}
\makeatother

\makeatletter
\renewcommand\theindex{\@restonecoltrue\if@twocolumn\@restonecolfalse\fi
  \columnseprule\z@ \columnsep 35\p@
  \twocolumn[\@xp\part\@xp*\@xp{\bf Index}\bigskip]%
  \let\item\@idxitem
  \parindent\z@  \parskip\z@\@plus.3\p@\relax
  \small}

\makeatletter
\newcommand{\smallbullet}{} 
\DeclareRobustCommand\smallbullet{%
  \mathord{\mathpalette\smallbullet@{0.6}}%
}
\newcommand{\smallbullet@}[2]{%
  \vcenter{\hbox{\scalebox{#2}{$\m@th#1\bullet$}}}%
}
\makeatother

\begin{document}

\title{Constructing $\upo$-free Hardy fields}
\author[Aschenbrenner]{Matthias Aschenbrenner}
\address{Kurt G\"odel Research Center for Mathematical Logic\\
Universit\"at Wien\\
1090 Wien\\ Austria}
\email{matthias.aschenbrenner@univie.ac.at}

\author[van den Dries]{Lou van den Dries}
\address{Department of Mathematics\\
University of Illinois at Urbana-Cham\-paign\\
Urbana, IL 61801\\
U.S.A.}
\email{vddries@illinois.edu}

\author[van der Hoeven]{Joris van der Hoeven}
\address{CNRS, LIX (UMR 7161)\\ 
Campus de l'\'Ecole Polytechnique\\  91120 Palaiseau \\ France}
\email{vdhoeven@lix.polytechnique.fr}

\date{March, 2026}

\begin{abstract} 
We show that every Hardy field extends to an $\upo$-free Hardy field. This result relates to classical oscillation criteria for second-order homogeneous linear differential equations. It is essential in \cite{ADH5}, and here we apply it to answer questions of Boshernitzan, and to generalize a theorem of his. 
\end{abstract}

\dedicatory{To the memory of Michael Boshernitzan \textup{(}1950--2019\textup{)}}

\pagestyle{plain}
 
\maketitle



\section*{Introduction}

\noindent
We let $a$, $b$, $c$ range over $\R$ in this introduction. Let  $f\colon  [a,+\infty)\to\R$ be continuous, and consider the second-order linear differential equation
\begin{equation}\tag{$\ast$}\label{eq:ast}
Y''+fY\ =\ 0.
\end{equation}
A (real) {\it solution}\/ to \eqref{eq:ast} is a $\c^2$-function $y\colon [a,+\infty)\to\R$   such that~${y''+fy=0}$, and such a solution is either (identically) zero, or its zero set as a subspace of~$[a,+\infty)$ is discrete. 
Some equations \eqref{eq:ast} have oscillating solutions. Here,
a  continuous function~$g\colon [a,+\infty)\to\R$      {\it oscillates}\/ if~$g(t)=0$ for arbitrarily large $t\geq a$, and~$g(t)\ne 0$ for arbitrarily large~$t\geq a$. Every   oscillating solution to \eqref{eq:ast} has arbitrarily large isolated zeros, whereas each
nonzero non-oscillating solution to~\eqref{eq:ast} has  only finitely many zeros.
  We say that {\it $f$ generates oscillation}\/ if \eqref{eq:ast} has an oscillating solution. In this case,  by 
 Sturm~\cite{Sturm},    {\it every}\/ nonzero  solution   to \eqref{eq:ast} oscillates.
This is really a property of the germ
of $f$ at $+\infty$: for $b\geqslant a$, 
$f$ generates oscillation iff~$f|_{[b,+\infty)}$ does.  Below ``germ'' means ``germ at $+\infty$" and ``oscillates'' and ``generates oscillation''  will hold by convention for the germ of $f$ iff it holds for~$f$. 

\medskip
\noindent
There is an extensive literature giving criteria on $f$ to generate oscillation (see for example \cite{Hinton,Swanson,Willett}),
some of which have their root in another fundamental result of Sturm, his   Comparison Theorem: for continuous $g: [a,+\infty)\to \R$,
$$\textit{if $f$ generates oscillation and $f\leq g$ on $[a,+\infty)$, then $g$   generates oscillation.}$$
To see this theorem in action, let $\ell_0:=x$ be the germ of the identity function on~$\R$ and recursively define the germs 
$\ell_1,\ell_2,\dots$ by $\ell_{n+1}:=\log \ell_n$ for each $n$. Also set
$$\upg_n\, :=\,  \frac{1}{\ell_0\cdots\ell_n},\qquad \upo_n \, :=\, \upg_0^2+\upg_1^2+\cdots+\upg_n^2 \,=\, \frac{1}{\ell_0^2}+\frac{1}{(\ell_0\ell_1)^2}+\cdots+\frac{1}{(\ell_0\cdots\ell_{n})^2}.$$
Computation shows that the germ $y=\upg_n^{-1/2}$  satisfies the equation~$y''+(\upo_n/4)y=0$, so $\upo_n/4$ doesn't
generate oscillation. 
More precisely, a germ
$$\frac{\upo_n+c\upg_{n}^2}{4}\ =\ \frac{1}{4}\left(\frac{1}{\ell_0^2}+\frac{1}{(\ell_0\ell_1)^2}+\cdots+\frac{1}{(\ell_0\cdots\ell_{n-1})^2} + \frac{c+1}{(\ell_0\cdots\ell_{n})^2}\right)$$ 
generates oscillation if and only if $c>0$.
(See \cite[Chapter~6, Theorem~10]{Bellman} or Corollary~\ref{cor:17.7 generalized} below.)
In the rest of this introduction $\c$ is the ring of germs (at~$+\infty$) of continuous real-valued functions whose domain is a subset of $\R$ containing some interval $(a,+\infty)$,  and $f$, $g$ range over $\c$.
We partially order $\c$ by: 
$$f\leqslant g\ :\Longleftrightarrow\ f(t)\leqslant g(t), \text{ eventually,}$$
where tacitly $f$ and $g$ also denote representatives of their germs, and ``eventually'' means ``for all sufficiently large real $t$". 
So if  $f\leq\upo_n/4$ for some~$n$, then~\eqref{eq:ast} has no oscillating solutions,
whereas if~${f\geq (\upo_n+c\upg_{n}^2)/4}$ for some~$n$ and~${c>0}$, then~$f$ generates oscillation.
For~$n=0$ this was first noted by A.~Kneser~\cite{AKneser}: if~$f\leq 1/{4x^2}$, then~$f$ does not generate oscillation,
but if $f\geq (1+c)/4x^2$ for some~$c>0$, then~$f$ generates oscillation.
The general case is implicit in
Riemann-Weber~\cite[p.~63]{Weber}, and was rediscovered by various authors~\cite{Hartman48, Hille48, Miller}. 

\medskip
\noindent
Now $\upo_n\leqslant \upo_{n+1}\leqslant \upo_{n+1}+\upg_{n+1}^2 \leqslant \upo_{n}+c\upg_n^2$ for all $n$ and all $c>0$, and
it is not difficult to obtain a germ $\upo$  
such that $\upo_n \leq \upo \leq \upo_n+\upg_n^2$ for all $n$, hence
the Riemann-Weber criterion is inconclusive for $f=\upo/4$.
(One can even take such $\upo$ to be the germ of an analytic function $(a, +\infty)\to \R$, by our
Example~\ref{ex:non-upo-free}.)
 However, Hartman~\cite{Hartman48} observed that if
 $f$ is the germ of a {\it logarithmico-exponential function}\/ ({\it LE-function}\/) in the sense of Hardy~\cite{Har12a,Ha}, then
 this criterion applies:
\begin{equation}\label{eq:Hartman}\tag{H}
\textit{$f$ does not generate oscillation} \qquad\Longleftrightarrow\qquad\textit{$f\leq\upo_n/4$   for some~$n$.}
\end{equation}
Roughly speaking, LE-functions are the real-valued functions
  obtained in finitely many steps from real constants and the variable~$x$ using addition, multiplication, division, exponentiation, and taking logarithms. Examples: every rational function in $\R(x)$, each power function  $x^r$ ($r\in\R$), each  iterated logarithm $\ell_n$,
  $\ex^{\sqrt{\ell_1}/\ell_2}$, etc.
Hardy showed that each LE-function, defined on some interval $(a,+\infty)\to\R$, has eventually constant sign, and so the
germs at $+\infty$ of such functions form what Bourbaki~\cite{Bou} later called a {\it Hardy field}\/: a subfield $H$ of the ring $\c$  consisting of germs of continuously differentiable functions 
$(a,+\infty)\to\R$  such that the germ of the derivative of the function is also in $H$ (so $H$ is a differential field). 

\medskip
\noindent
The Hardy field~$H_{\operatorname{LE}}$ of LE-functions 
has  good uses (see~\cite{CS,HL14,Mac81}), but overall  is
too small for many analytic purposes: for example, every $h\in H_{\operatorname{LE}}$ is {\em differentially algebraic over $\R$} (that is, satisfies a differential equation $P\big(h, h',\dots, h^{(n)}\big)=0$ with~$P$ a nonzero polynomial over $\R$ in $1+n$ indeterminates), and yet  $H_{\operatorname{LE}}$ contains no antiderivative of $\ex^{x^2}$ (Liouville~\cite{Liouville1, Liouville2}).
Boshernitzan~\cite[Theorem~17.7]{Boshernitzan82} (see also Corollary~\ref{cor:17.7})   generalized Hartman's result and showed that the equivalence~\eqref{eq:Hartman} continues to hold provided~$f$ is merely assumed to be  {\it hardian}\/ (i.e.,
contained in some  Hardy field) and to be differentially algebraic over $\R$. In~\cite[Conjecture~17.11]{Boshernitzan82} he conjectured  a version of \eqref{eq:Hartman}  with the increasing sequence
$(\upo_n)$ replaced by the decreasing sequence~$(\upo_n+\upg_n^2)$, and ``for some" replaced accordingly by ``for all'':
if $f$ is hardian and differentially algebraic over $\R$, then
\begin{equation}\label{eq:Boshernitzan}\tag{B}
\textit{$f$ does not generate oscillation} \qquad\Longleftrightarrow\qquad\textit{$f\leq (\upo_n+\upg_n^2)/4$   for all~$n$.}
\end{equation}
Corollary~\ref{cor:17.11} below establishes this conjecture along the way to our main result. 
Important here is (cf.~Theorem~\ref{thm:Ros83} and Corollary~\ref{cor:Bosh 14.11}) that if $f$ is hardian and differentially algebraic over $\R$, then there is an $n$ such that
$\ell_n\leq g$
for all positive infinite $g$ in $\R\<f\>:=\R(f, f',f'',\dots)$ (=~the Hardy field generated by $f$ over $\R$); here
{\it $g$ is positive infinite}\/ ($g\in \c$) means that
  $g(t)\to+\infty$ as~$t\to+\infty$.

\medskip
\noindent
The germ $f$ being hardian and not generating oscillation has nice consequences:
for example, each Hardy field  containing such an~$f$  extends to
one which contains a fundamental system of solutions of \eqref{eq:ast}, and a   Hardy-type inequality
with weight~$f$ holds. (See Proposition~\ref{prop:2nd order Hardy field} and Remarks~\ref{rem:Hardy inequ, 1},~\ref{rem:Hardy inequ, 2}, respectively.) Hence it is desirable to have versions of \eqref{eq:Hartman} and \eqref{eq:Boshernitzan}
 for arbitrary hardian~$f$.
Many natural functions, for example the restrictions of
Euler's $\Gamma$-function and Riemann's $\zeta$-function to $(1,+\infty)$, have non-oscillating germs that are hardian but are not differentially algebraic over $\R$~\cite{Rosenlicht83}.  
More trouble are  hardian germs
  $\upo$ as above, that is, $\upo_n \leq \upo\leq \upo_n+\upg_n^2$ for all~$n$, as in our Example~\ref{ex:non-upo-free}.
  In this paper we show that nevertheless, {\it   versions of~\eqref{eq:Hartman} and~\eqref{eq:Boshernitzan} can be restored if we extend  the relevant Hardy field 
  and prolong the sequence~$(\ell_n)$ of iterated logarithms  accordingly.}\/
  
\medskip
\noindent   
To make this precise, let~$H$ be a Hardy field.  The  partial ordering~$\leq$ on~$\c$ restricts to a total ordering on~$H$, making~$H$
an ordered field.  We also assume: $H\supseteq\R(x)$
and~$H$ is  {\it log-closed,}\/ that is, $\log h\in H$ for all~$h>0$ in $H$.
(This holds for $H=H_{\operatorname{LE}}$, and every Hardy field extends to one with these properties.)  Then~$H\supseteq \R(\ell_0,\ell_1,\ell_2,\dots)$ and
 we extend the sequence~$\ell_0,\ell_1,\ell_2,\dots$   by transfinite recursion to a sequence~$(\ell_\rho)$ of positive infinite elements of $H$,
 indexed by all ordinals $\rho$ less than some infinite limit ordinal~$\kappa$, as follows: $\ell_{\rho+1}:=\log\ell_\rho$, and
 if $\lambda$ is an infinite limit ordinal such that all~$\ell_\rho$ with~$\rho<\lambda$ have already been chosen, then we pick
 $\ell_\lambda$ to be any positive infinite element of~$H$   such that $\ell_\lambda \leq \ell_\rho$ for all  $\rho<\lambda$,
 if there is such an element; otherwise we put~$\kappa:=\lambda$. 
 Given~$z\in H$  we   set~$\omega(z):=-(2z'+z^2)$; then for   $y\in H\setminus\{0\}$ we have~$y''+fy=0$ iff $z:=2y'/y$ satisfies
 the (Riccati) equation~$\omega(z)=4f$. 
We now define 
 $$\upg_\rho\, :=\, \ell_\rho'/\ell_\rho,\qquad  \upo_\rho\,:=\,\omega(-\upg_\rho'/\upg_\rho).$$
We have $\upg_n:=\ell_n'/\ell_n=1/(\ell_0\cdots\ell_n)$, and
 taking $y:=1/\sqrt{\upg_n}$ we obtain 
 $$z\ :=\ 2y'/y\ =\ -\upg_n'/\upg_n\ =\ \upg_0+\cdots+\upg_n, \qquad
 \omega(z)\ =\ \upg_0^2+\upg_1^2+\cdots+\upg_n^2.$$ (Note that these $\upg_n$, $\upo_n$ agree with the $\upg_n, \upo_n$ given earlier.) In the beginning of Section~\ref{sec:upo-free Hardy fields} we show that the sequences $(\upo_{\rho})$ and $(\upo_{\rho}+\upg_{\rho}^2)$ in $H$ are
 strictly increasing and strictly decreasing, respectively, and that $\upo_{\lambda} < \upo_{\mu}+\upg_{\mu}^2$ for all indices $\lambda, \mu$. 

\medskip
\noindent
Using results from \cite{Boshernitzan82,Ros} it is easy (see Corollary~\ref{hstar} and subsequent comment) to extend any Hardy field to an $H$ as above (that is, $H\supseteq \R(x)$ and $H$ is log-closed) such that the following variant of~\eqref{eq:Hartman} holds for all $f\in H$:  
\begin{equation}\label{eq:Hartman H}\tag{H$^*$}
 \textit{$f$ does not generate oscillation}\qquad\Longleftrightarrow\qquad\textit{$f \leq \upo_\rho/4$ for some $\rho$.}
\end{equation}
Restoring \eqref{eq:Boshernitzan} requires the concept of a Hardy field being
{\em $\upo$-free}. This (first-order) concept was introduced in a more general setting in our book~[ADH, 11.7], where it was shown to be very robust. (For example, by [ADH, 13.6.1] it is preserved under passage to   differentially algebraic Hardy field extensions.) We repeat the formal definition in Section~\ref{sec:prelim}, and Corollary~\ref{bstar} says that $H$ is  $\upo$-free if and only if for all~$f\in H$ the equivalence \eqref{eq:Hartman H} holds as well as the following equivalence:
\begin{equation}\label{eq:Boshernitzan H}\tag{B$^*$}
\textit{$f$ does not generate oscillation}\qquad\Longleftrightarrow\qquad\textit{$f \leq (\upo_\rho+\upg_\rho^2)/4$ for all $\rho$.}
\end{equation}
Keeping in mind that $H$ ranges over log-closed Hardy fields containing $\R(x)$, here is the main result of this paper, already announced  in~\cite{HIVP,ADH2}:

\begin{theoremUnnumbered}
Every Hardy field is contained in some $\upo$-free $H$. 
\end{theoremUnnumbered}

\noindent
A more precise version is given by Theorem~\ref{upo}. 
There are $H$ that are not $\upo$-free, but those with a natural origin usually are. (For example,  
by Hartman's and Boshernitzan's oscillation criteria,
 $H_{\operatorname{LE}}$ is $\upo$-free.) The proof of Theorem~\ref{upo} takes nevertheless considerable effort.
 
After the preliminary Section~\ref{sec:prelim} we give in
Section~\ref{sec:germs} basic definitions and facts about germs
of one-variable (real- or complex-valued) functions, and in Section~\ref{sec:differentiable germs} we collect
the main facts we need about second-order linear differential equations.
In Section~\ref{sec:Hardy fields} we introduce Hardy fields in more detail and review some extension results  due to
Boshernitzan~\cite{Boshernitzan81, Boshernitzan82, Boshernitzan86} and Rosenlicht~\cite{Ros}. In
Section~\ref{sec:upper lower bds} we discuss upper and lower bounds on the growth of hardian germs  from \cite{Boshernitzan82,Boshernitzan86,Rosenlicht83}, and   Section~\ref{sec:order 2 Hardy fields} focusses on second-order linear differential
equations over Hardy fields.
In Section~\ref{sec:upo-free Hardy fields} we review $\upo$-freeness, prove the theorem above, and some refinements.

\medskip
\noindent
By Zorn, each Hardy field is contained in  one which is {\it maximal}\/, that is,
which has no proper Hardy field extension.
By the   theorem above,    maximal Hardy fields are $\upo$-free.
This is a first important step towards showing that
  they  are {\it $H$-closed fields}\/, in the terminology of \cite{ADH2}. This requires serious further work, which is in~\cite{ADH5}.

\medskip
\noindent
In    Section~\ref{sec:applications} we show that
our main theorem  by itself,
combined with results about $\upo$-freeness from [ADH], 
already     has some  applications.
First, it yields a positive answer to
a question about maximal Hardy fields posed by Boshernitzan~\cite[\S{}7]{Boshernitzan86}:
 
\begin{corintro}
Every maximal Hardy field contains a positive infinite  germ~$\ell_\omega$ such that $\ell_\omega \leq \ell_n$ for all~$n$.
\end{corintro}

\noindent
This corollary is actually much weaker than~\cite[Corollary~4.8]{ADHfgh}, which however ultimately relies on deeper results from \cite{ADH5} that in turn depend on Theorem~\ref{upo}. 

\medskip
\noindent
In the remainder of this introduction we let $H$ range over arbitrary Hardy fields.
The intersection $\Ex(H)$ of all maximal Hardy fields that contain $H$ is a Hardy field that is log-closed and properly contains $H_{\operatorname{LE}}$.  
These Hardy fields~$\Ex(H)$ were studied in detail by Boshernitzan~\cite{Boshernitzan81,Boshernitzan82}, who proved, among other things, that 
the sequence~$(\ell_n)$ is coinitial in the set of positive infinite elements of~$\Ex(\Q)=\Ex(H_{\operatorname{LE}})$.
Theorem~\ref{thm:coinitial in E(H)} generalizes this fact as follows:

\begin{corintro}
If $H\supseteq \R(x)$ is log-closed and $\upo$-free, then any log-sequence $(\ell_\rho)$ in~$H$ as above
is coinitial in the set $\Ex(H)^{>\R}$ of  positive infinite elements of $\Ex(H)$.
\end{corintro}

\subsection*{Notations and conventions}
We generally follow the conventions from [ADH]. In particular,  $m$,~$n$   range over the set~$\N=\{0,1,2,\dots\}$ of natural numbers.
Given an additively written abelian group $A$ we set $A^{\ne}:=A\setminus\{0\}$. By convention, the ordering of an ordered abelian group or ordered field is {\em total}. For an ordered abelian group 
$A$ and $b\in A$ we put $A^{>b}:=\{a\in A:a>b\}$ and $A^>:=A^{>0}$, and likewise
with $\geq$, $<$, or $\leq$ in place of $>$.
 Rings are associative  with  identity~$1$ (and almost always commutative). For a commutative 
 ring~$R$ we let~$R^\times$ be the multiplicative group of units of $R$.

\subsection*{Acknowledgements}
We thank the anonymous referee for a careful reading of the paper and useful comments. 
 
\section{Preliminaries on Asymptotic Fields}\label{sec:prelim}

\noindent
In this section we first collect some basic definitions from [ADH] needed throughout this paper.
We then review some general facts on iterated logarithmic derivatives, iterated exponentials, and the asymptotic behavior of ``large'' solutions of algebraic differential equations
in $H$-asymptotic fields.
We do not need these facts to achieve our main objective,  but they will be used at a few points for
  applications and corollaries; see Sections~\ref{sec:upper lower bds} and~\ref{sec:applications}.

\subsection*{Differential rings and fields}
Let $R$ be a {\em differential ring\/}, that is, a commutative ring $R$ containing
(an isomorphic copy of) $\Q$ as a subring equipped with a derivation~${\der\colon R \to R}$. Then $C_R:=\ker\der$ is a subring of~$R$, called the ring of constants of~$R$, and~${\Q\subseteq C_R}$.  If~$R$ is a field, then so is $C_R$. A {\it differential field}\/ is a differential ring that happens to be a field.
When the derivation $\der$ of $R$ is clear from the context and $a\in R$, then we denote~$\der(a),\der^2(a),\dots,\der^n(a),\dots$ by~$a', a'',\dots, a^{(n)},\dots$, and for~$a\in R^\times$ we set $a^\dagger:=a'/a$ (the {\it logarithmic derivative}\/  of $a$), so $(ab)^\dagger=a^\dagger + b^\dagger$ for~$a,b\in R^\times$.

\medskip
\noindent
We have the differential ring $R\{Y\}=R[Y, Y', Y'',\dots]$ of differential polynomials in a differential indeterminate $Y$ over $R$. We say that  
$P=P(Y)\in R\{Y\}$ has {\it order}\/ at most~$r\in \N$ if
$P\in R[Y,Y',\dots, Y^{(r)}]$; in this case $P=\sum_{\i}P_{\i}Y^{\i}$, as in [ADH, 4.2], with~$\i$ ranging over tuples
$(i_0,\dots,i_r)\in \N^{1+r}$, $Y^{\i}:= Y^{i_0}(Y')^{i_1}\cdots (Y^{(r)})^{i_r}$,     coefficients~$P_{\i}$   in $R$, and $P_{\i}\ne 0$ for only finitely many $\i$. 

\medskip
\noindent
For $P\in R\{Y\}$ and~$a\in R$ we let $P_{\times a}:=P(aY)$.
For $\phi\in R^\times$ we let $R^{\phi}$ be the {\it compositional conjugate}\/ of~$R$ by~$\phi$: the differential ring
with the same underlying ring as~$R$ but with derivation~$\phi^{-1}\der$ (usually denoted by $\derdelta$) instead of $\der$. 
We have an $R$-algebra isomorphism~$P\mapsto P^\phi\colon R\{Y\}\to R^\phi\{Y\}$ such that $P^\phi(y)=P(y)$ for all $y\in R$;
see [ADH, 5.7].  

\medskip
\noindent
Let $K$ be a differential field. Then $K\{Y\}$ is an integral domain, and the differential fraction field of  $K\{Y\}$ is denoted by~$K\<Y\>$. 
Let $y$ be an element of a differential field extension $L$ of $K$. We let~$K\{y\}$ be the  
differential subring of $L$ generated by~$y$ over $K$, and let $K\<y\>$ be the differential fraction field of $K\{y\}$
in $L$. We say that~$y$ is 
{\it differentially algebraic\/} over $K$ if  $P(y)=0$ for some
 $P\in K\{Y\}^{\ne}$; otherwise $y$ is  called {\it differentially transcendental\/} over $K$.
  As usual in [ADH], the prefix~``$\d$'' abbreviates ``differentially'', so ``$\d$-algebraic'' means ``differentially algebraic''.
We say that $L$ is 
 {\it $\d$-algebraic\/} over $K$ if each $y\in L$ is $\d$-algebraic over $K$.
 See [ADH, 4.1] for more on this. We set $K^\dagger:= (K^\times)^\dagger$, a subgroup of the additive group of $K$.

\subsection*{Valued fields}
For a field $K$ we have $K^\times=K^{\ne}$, and
a (Krull) valuation on $K$ is a surjective map 
$v\colon K^\times \to \Gamma$ onto an ordered abelian 
group $\Gamma$ (additively written) satisfying the usual laws, and extended to
$v\colon K \to \Gamma_{\infty}:=\Gamma\cup\{\infty\}$ by $v(0)=\infty$,
where the ordering on $\Gamma$ is extended to a total ordering
on $\Gamma_{\infty}$ by $\gamma<\infty$ for all~$\gamma\in \Gamma$. 
A {\em valued field\/} $K$ is a field (also denoted by $K$) together with a valuation ring~$\mathcal O$ of that field,
and the corresponding valuation $v\colon K^\times \to \Gamma$  on the underlying field is
such that $\mathcal O=\{a\in K:va\geq 0\}$ as explained in [ADH, 3.1].

\medskip
\noindent
Let $K$ be a valued field
with valuation ring $\mathcal O$ and valuation $v\colon K^\times \to \Gamma$. Then~$\mathcal O$ is a local ring  with maximal ideal $\smallo:=\{a\in K:va>0\}$. 
In this paper $K$ always has {\it equicharacteristic zero}\/, that is, the
residue field $\res(K):=\mathcal{O}/\smallo$ of $K$   has characteristic zero.
In asymptotic differential algebra, sometimes the following  notations are more natural:
with $a$, $b$ ranging over $K$, set 
\begin{align*} a\asymp b &\ :\Leftrightarrow\ va =vb, & a\preceq b&\ :\Leftrightarrow\ va\ge vb, & a\prec b &\ :\Leftrightarrow\  va>vb,\\
a\succeq b &\ :\Leftrightarrow\ b \preceq a, &
a\succ b &\ :\Leftrightarrow\ b\prec a, & a\sim b &\ :\Leftrightarrow\ a-b\prec a.
\end{align*}
It is easy to check that if $a\sim b$, then $a, b\ne 0$ and $a\asymp b$, and that
$\sim$ is an equivalence relation on $K^\times$.\label{p:asymprels} 
Let  $L$ be a valued field extension of $K$; then
the relations~$\asymp$,~$\preceq$,  etc.~on $L$ restrict to the corresponding relations on $K$, and
we identify in the usual way the value group of $K$ with an ordered subgroup of the value group of~$L$
and $\res(K)$ with a subfield of $\res(L)$. Such a valued field extension is called {\it immediate} if 
for every $a\in L^\times$ there is a $b\in K^\times$ with $a\sim b$. We use {\em pc-sequence\/} to abbreviate
{\em pseudocauchy sequence}, and~${a_\rho\leadsto a}$ indicates that~$(a_\rho)$ is a  pc-sequence pseudoconverging to~$a$;
here the $a_{\rho}$ and $a$ lie in a valued field understood from the context, see [ADH, 2.2,~3.2].\label{p:pc} 

\medskip
\noindent
A binary relation $\preceq$ on a field $K$ for which there is a valuation $v$ on $K$ such that~$a\preceq b\Leftrightarrow va\geq vb$
for each $a,b\in K$ is called a {\it dominance relation}\/ on $K$. See~[ADH, 3.1] for  an axiomatization of dominance relations.

\subsection*{Valued differential fields}
As in~[ADH],  a {\em valued differential field\/} is a valued field of equicharacteristic zero together with a derivation, generally denoted by~$\der$, on the underlying field.
The derivation $\der$ of a valued differential field $K$ is said to be {\it small}\/ if~$\der\smallo\subseteq\smallo$; then~${\der\mathcal O\subseteq\mathcal O}$~[ADH, 4.4.2], so $\der$ induces a derivation on~$\res(K)$ making the residue map ${\mathcal O\to\res(K)}$ into a morphism of differential rings.
A valued differential field $K$ in this paper is usually an 
\textit{asymptotic field}\/, that is, for all nonzero  $f,g\prec 1$ in $K$ we have: $f\preceq g\Longleftrightarrow f'\preceq g'$. 
Every compositional conjugate of an asymptotic field is asymptotic.

\medskip
\noindent
Let $K$ be an asymptotic field,  with constant field $C=C_K$ and valuation ring $\mathcal{O}$. Then $C\subseteq \mathcal{O}$, and we say that $K$ is
  {\it $\d$-valued}\/ if for all $f\in K$ with $f\asymp 1$ there is a~$c\in C$ with $f\sim c$.
 Let $\I(K)$ be the $\mathcal O$-submodule of $K$ generated by $\der\mathcal O$.
Then~$K$ is called {\it pre-$\d$-valued}\/ if $\I(K) \cap (K^\times)^\dagger = (\mathcal O^\times)^\dagger$.
(This is not exactly the definition from [ADH, 10.1], but equivalent to it.) 
Pre-$\d$-valued fields are exactly the valued differential subfields of $\d$-valued fields, by  \cite[4.4]{AvdD2}.

\medskip
\noindent
We associate to~$K$ its {\em asymptotic couple\/}   $(\Gamma,\psi)$,
where $\psi\colon\Gamma^{\neq}\to\Gamma$ is given by 
$$\psi(vg)\ =\ v(g^\dagger)\ \text{ for $g\in K^\times$ with $vg\ne 0$}.$$ 
We put $\Psi:=\psi(\Gamma^{\neq})$. 
If we want to stress the dependence on $K$, we write~$(\Gamma_K,\psi_K)$ and $\Psi_K$ instead of $(\Gamma,\psi)$ and $\Psi$, respectively.
An {\it asymptotic couple}\/ (without mentioning any asymptotic field) is a pair $(\Gamma,\psi)$ consisting 
of an ordered abelian group $\Gamma$ and a map $\psi\colon\Gamma^{\neq}\to\Gamma$  subject to natural axioms
obeyed by the asymptotic couples of asymptotic fields, see [ADH, 6.5]. 
We  
extend~$\psi\colon\Gamma^{\neq}\to\Gamma$ to a map~${\Gamma_\infty\to\Gamma_\infty}$
by $\psi(0):=\psi(\infty):=\infty$.   If $(\Gamma,\psi)$ is understood from the context and $\gamma\in\Gamma$ we write $\gamma^\dagger$ and $\gamma'$ instead of $\psi(\gamma)$ and $\gamma+\psi(\gamma)$, respectively.
An {\it $H$-asymptotic couple}\/ is an asymptotic couple $(\Gamma,\psi)$ such that for all $\gamma,\delta\in\Gamma$ we have:
$0<\gamma\leq\delta\Rightarrow \psi(\gamma)\geq\psi(\delta)$.
An asymptotic field whose asymptotic couple is $H$-asymptotic is called an
{\it $H$-asymptotic field}\/ (or an asymptotic field of {\it $H$-type}\/).

Let $(\Gamma,\psi)$ be an asymptotic couple and $\Psi:=\psi(\Gamma^{\neq})$. Then $\gamma\in\Gamma$ is said to be a {\it gap}\/ in $(\Gamma,\psi)$
if~$\Psi<\gamma<(\Gamma^>)'$. (There is at most one such $\gamma$.)
We also say that~$(\Gamma,\psi)$ is {\it grounded}\/ if $\Psi$ has a largest element, 
and $(\Gamma,\psi)$ has {\it asymptotic integration}\/ if~$(\Gamma^{\neq})'=\Gamma$.
An asymptotic field is said to have a gap if  its asymptotic couple does,  and likewise with ``grounded'' or ``asymptotic integration'' in place of ``has a gap''.
See [ADH, 9.1, 9.2] for more on this,
in particular for  the following important trichotomy: every $H$-asymptotic couple  either
 has a gap, or  is grounded, or has asymptotic integration [ADH, 9.2.16].
  An element $\phi$ of an asymptotic field~$K$ is said to be {\em active in $K$\/} if $\phi\succeq f^\dagger$ for some $f\nasymp 1$ in $K^\times$; in that case the derivation~$\phi^{-1}\der$ of the compositional conjugate $K^\phi$ is small, cf.~[ADH, 11.1]. 

Next two concepts from [ADH, 11.6, 11.7] that may seem technical but that are key to understanding subtler aspects of Hardy fields. Let $K$ be an asymptotic field. We say that $K$ is {\it $\upl$-free\/}  if $K$ is $H$-asymptotic and ungrounded, and for all~$f\in K$ there exists $g\succ 1$ in $K$ such that~$f-g^{\dagger\dagger}\succeq g^\dagger$. We say that $K$ is {\it $\upo$-free\/} if $K$ is $H$-asymptotic, ungrounded, and for all $f\in K$ there exists $g\succ 1$ in~$K$ such that~$f-\omega(g^{\dagger\dagger})\succeq (g^\dagger)^2$, where
$\omega(z):=-2z'-z^2$ for $z\in K$.  This notion of $\upo$-freeness is clearly first-order in the logical sense. 
If $K$ is $\upo$-free, then $K$ is $\upl$-free~[ADH, 11.7.3], and if $K$ is $\upl$-free, then $K$ has
asymptotic integration~[ADH, 11.6.8].
(We do not use $\upl$-freeness or $\upo$-freeness before Section~\ref{sec:order 2 Hardy fields}.)
  
\subsection*{Ordered differential fields}
  An  {\em ordered differential field\/} is a differential field $K$ with an ordering on $K$ making $K$ an ordered field. Likewise, an  {\em ordered valued differential field\/} is a  valued differential field $K$ equipped an ordering on $K$ making~$K$ an ordered field (no relation between derivation, valuation, or ordering
being assumed).  Let $K$ be an ordered differential field. Then we have the convex subring 
$$\mathcal{O}\ :=\ \big\{g\in K:\text{$\abs{g} \le c$ for some $c\in C$}\big\},$$ 
which is a valuation ring of $K$ and has maximal ideal
$$\smallo\ =\ \big\{g\in K:\text{$\abs{g} < c$ for all positive $c\in C$}\big\}.$$ We call $K$ an {\it $H$-field}\index{ordered differential field!$H$-field}\index{H-field@$H$-field} 
 if for all $f\in K$ with $f>C$ we have $f'>0$, and
  $\mathcal{O}=C+\smallo$. We view such an $H$-field $K$ as an ordered valued differential field with its valuation given by $\mathcal O$. 
{\it Pre-$H$-fields}\/  are the ordered valued differential subfields of $H$-fields. Every pre-$H$-field is   $H$-asymptotic, and each $H$-field is $\d$-valued of $H$-type.
 See [ADH, 10.5] for basic facts about (pre-)$H$-fields. 
 An $H$-field $K$ is said to be {\it Liouville closed}\/
 if $K$ is real closed and for all $f,g\in K$ there exists~$y\in K^\times$ with~$y'+fy=g$.
 Every $H$-field extends to  a Liouville closed one; see [ADH, 10.6].\index{H-field@$H$-field!Liouville closed}\index{closed!Liouville}

\medskip
\noindent
{\it In the rest of this section $K$ is an $H$-asymptotic field, and $f$, $g$ range over $K$.}\/

\subsection*{Iterated logarithmic derivatives}
Let $(\Gamma,\psi)$ be an $H$-asymp\-to\-tic couple. We let~$\gamma$ range over~$\Gamma$, and we denote by
$$[\gamma]=\big\{\delta\in\Gamma: \text{$\abs{\gamma} \leq n\abs{\delta}$ and $\delta\leq n\abs{\gamma}$ for some $n\geq 1$} \big\}$$ the archimedean class of $\gamma$; cf.~[ADH, 2.4].
We define $\gamma^{\langle n\rangle}\in\Gamma_\infty$ recursively by~$\gamma^{\langle 0\rangle}:=\gamma$
and $\gamma^{\langle n+1\rangle}:=
\psi(\gamma^{\langle n\rangle})$.
The following is~\cite[Lem\-ma~5.2]{A}; for the convenience of the reader we include a
proof: 

\begin{lemma}\label{lem:itpsi}
Suppose that $0\in (\Gamma^<)'$, $\gamma\neq 0$,  and $n\geq 1$. If $\gamma^{\langle n\rangle}<0$, then~$\gamma^{\langle i\rangle}<0$ for $i=1,\dots,n$ and
$[\gamma] > [\gamma^\dagger] >\cdots > [\gamma^{\langle n-1\rangle}] >  [\gamma^{\langle n\rangle}]$.
\end{lemma}
\begin{proof}
By [ADH, 9.2.9] we have $(\Gamma^>)'\subseteq\Gamma^>$,  so the case $n=1$
follows from~[ADH, 9.2.10(iv)].
Assume inductively that the lemma holds for a certain value of~$n\geq 1$, and suppose $\gamma^{\langle n+1\rangle}<0$. Then $\gamma^{\langle n\rangle}\neq 0$, so we can apply the case~$n=1$ to
$\gamma^{\langle n\rangle}$ instead of $\gamma$ and get $[\gamma^{\langle n\rangle} ]>
 [\gamma^{\langle n+1 \rangle}]$. By the inductive assumption the remaining
inequalities will follow from 
$\gamma^{\langle n\rangle}<0$. 
From $0\in (\Gamma^<)'$ we obtain an element~$1$ of~$\Gamma^>$ with $0=(-1)'=-1+1^\dagger$. Suppose $\gamma^{\langle n\rangle}\ge 0$. 
Then $\gamma^{\langle n\rangle}\in \Psi$, thus~$0<\gamma^{\langle n\rangle}<1+1^\dagger=1+1$ and so
$[\gamma^{\langle n\rangle}]\leq [1]$.
Hence $0> \gamma^{\langle n+1\rangle} \geq 1^\dagger=1$, 
a contradiction.
\end{proof}
 
\noindent
Suppose now that $(\Gamma,\psi)$ is the asymptotic couple of $K$ and $y\in K^\times$.
In [ADH, p.~213],  we defined the $n$th iterated logarithmic derivative of    $y$: 
$y^{\langle 0\rangle}:=y$, and recursively, if~$y^{\langle n\rangle}\in K$ is defined and nonzero, then
 $y^{\langle n+1\rangle}:=(y^{\langle n\rangle})^\dagger$, while otherwise~$y^{\langle n+1\rangle}$ is not defined. (Thus if~$y^{\langle n\rangle}$ is defined, then so are $y^{\langle 0\rangle},\dots,y^{\langle n-1\rangle}$.)
If  
$(vy)^{\langle n \rangle}\ne \infty$, then $y^{\langle n\rangle}$  is defined and
 $v(y^{\langle n\rangle})= (vy)^{\langle n\rangle}\in\Gamma$. 
 Recall from [ADH, p.~383] that for $f,g\neq 0$,
$$f\flatter g\  :\Leftrightarrow\ f^\dagger\prec g^\dagger,\quad f\flattereq g\ :\Leftrightarrow\ f^\dagger\preceq g^\dagger, \quad f\comp g\ :\Leftrightarrow\ f^\dagger\asymp g^\dagger,$$ hence, assuming also $f, g\nasymp 1$,
$$f\flatter g\  \Rightarrow\ [vf] < [vg], \qquad[vf] \leq [vg] \ \Rightarrow \  f\flattereq g.$$ 
{\it In the rest of this  section we are given $x\succ 1$ in~$K$ with $x'\asymp 1$.}\/
Then $0\in (\Gamma^<)'$, so from  the previous lemma we obtain:

\begin{cor}\label{cor:itpsi}
If $y\in K^\times$,  $y\nasymp 1$, $n\ge 1$, and   $(vy)^{\langle n\rangle}<0$, then~$y^{\langle i\rangle}\succ 1$ for~$i=1,\dots,n$ and $[vy] > \big[v(y^\dagger)\big] > \dots > \big[v( y^{\langle n-1\rangle})\big] > \big[v(y^{\langle n \rangle})\big]$.
\end{cor}

\noindent
Let $\i=(i_0,\dots,i_n)\in\Z^{1+n}$ and $y\in K^\times$ be such that $y^{\langle n\rangle}$ is defined;   we   put 
$$y^{\langle\i\rangle}\ :=\ (y^{\langle 0\rangle})^{i_0}\cdots (y^{\langle n\rangle})^{i_n}\in K.$$ 
If $y^{\langle n\rangle}\ne 0$, then $\i\mapsto y^{\langle\i\rangle}\colon\Z^{1+n}\to K^\times$ is a group morphism.
Suppose now that~$y\in K^\times$, $(vy)^{\langle n\rangle}<0$, and $\i=(i_0,\dots,i_n)\in\Z^{1+n}$, 
$\i\neq 0$, and $m\in\{0,\dots,n\}$ is minimal with $i_m\neq 0$.
Then by Corollary~\ref{cor:itpsi}, $\big[v(y^{\langle\i\rangle})\big]=\big[v(y^{\langle m\rangle})\big]$. Thus if~$y\succ 1$, we have the equivalence
$y^{\langle\i\rangle} \succ 1\  \Leftrightarrow\ i_m\geq 1$.  
If $K$ is equipped with an ordering making it a pre-$H$-field and $y\succ 1$, then $y^\dagger >0$, so
$y^{\<i\>}>0$ for   $i=1,\dots,n$, and thus~$\sgn y^{\<\i\>} = \sgn y^{i_0}$.

\subsection*{Iterated exponentials}  
{\it In this subsection we assume   that $\Psi$ is downward closed.}\/
For $f\succ 1$ we have $f'\succ f^\dagger$, so we can and do choose $\Exp(f)\in K^\times$ such that~${\Exp(f)\succ 1}$ and~$\operatorname{E}(f)^\dagger \asymp f'$, hence $f\prec \Exp(f)$ and $f \flatter \Exp(f)$. 
Moreover,  if  $f,g \succ 1$, then $${f \prec g} \quad\Longleftrightarrow\quad \Exp(f) \flatter \Exp(g).$$
For~$f\succ 1$
  define~$\Exp_n(f)\in K^{\succ 1}$ recursively
by 
$$ \Exp_0(f)\ :=\ f,\qquad \Exp_{n+1}(f)\ :=\ \Exp\!\big(\!\Exp_n(f)\big),$$ 
and thus by induction
$$\Exp_n(f)\ \prec\ \Exp_{n+1}(f) \quad\text{ and }\quad\Exp_{n}(f)\ \flatter\ \Exp_{n+1}(f)\qquad \text{ for all $n$.}$$
{\it In the rest of this subsection $f\succeq x$, and $y$ ranges over elements of $H$-asymptotic extensions of $K$.}\/
The proof of the next lemma is like that of~\cite[Lem\-ma~1.3(2)]{AvdD3}.


\begin{lemma}\label{lem:expn, 1}
If $y\succeq\Exp_{n+1}(f)$, $n\geq 1$, then $y\neq 0$ and $y^\dagger\succeq\Exp_n(f)$.
\end{lemma}
\begin{proof}
If $y\succeq\Exp_2(f)$, then
$y\neq 0$, and using $\Exp_2(f)\succ 1$ we obtain
$$y^\dagger\ \succeq\ \Exp_2(f)^\dagger\  \asymp\ \Exp(f)'\  =\  \Exp(f)\Exp(f)^\dagger \ \asymp\  \Exp(f)f' \ \succeq\  \Exp(f),
$$
Thus the lemma holds for $n=1$. In general,~$\Exp_{n-1}(f) \succeq f \succeq x$, hence  the lemma follows from the case $n=1$ applied
to~$\Exp_{n-1}(f)$ in place of $f$.
\end{proof}

\noindent
An obvious induction on $n$ using Lemma~\ref{lem:expn, 1} shows: if $y\succeq \Exp_n(f)$, then $(vy)^{\langle n \rangle}\le vf <0$. We shall use this fact without further reference.

\begin{lemma}\label{lem:expn, 2}
If   $y \succeq \Exp_{n+1}(f)$, then $y^{\langle n\rangle}$ is defined and $y^{\langle n\rangle}  \succeq\Exp(f)$.
\end{lemma}
\begin{proof}
First note that if $y\neq 0$, $n\ge 1$, and $(y^\dagger)^{\langle n-1\rangle}$ is defined,
then $y^{\langle n\rangle}$ is defined and $y^{\langle n\rangle}=(y^\dagger)^{\langle n-1\rangle}$.
Now use induction on $n$ and Lemma~\ref{lem:expn, 1}.
\end{proof}

\begin{lemma}\label{lem:expn, 3}
If
$y\succeq\Exp_n(f^2)$, then  $y^{\langle n\rangle}$ is defined and $y^{\langle n\rangle}  \succeq f$, with
$y^{\langle n\rangle}  \succ f$ if~$f\succ x$.
\end{lemma} 
\begin{proof}
This is clear if $n=0$, so suppose    $y\succeq\Exp_{n+1}(f^2)$. Then by Lemma~\ref{lem:expn, 2} (applied with $f^2$ in place of $f$)  we have $y^{\langle n\rangle}  \succeq\Exp(f^2)\succ 1$, so
$$y^{\<n+1\>}\ =\ (y^{\<n\>})^\dagger\ \succeq\ \Exp(f^2)^\dagger\ \asymp\ (f^2)'\ =\ 2ff'\ \succeq\ f,$$
with~$y^{\<n+1\>} \succ f$ if $f\succ x$,
 as required. 
\end{proof}

\begin{cor}\label{cor:expn}
Suppose  $y\succeq \operatorname{E}_{n}(f^2)$, and let $\i\in\Z^{1+n}$ be such that $\i>0$ lexicographically. Then 
$y^{\<n\>}$ is defined and $y^{\langle \i\rangle}\succeq f$, with
$y^{\langle \i\rangle}\succ f$ if~$f\succ x$.
\end{cor}
\begin{proof}
By Lemma~\ref{lem:expn, 3}, $y^{\langle n\rangle}$ is defined with $y^{\langle n\rangle}  \succeq f$, and $y^{\langle n\rangle}  \succ f$ if $f\succ x$.
Let~$m\in\{0,\dots,n\}$ be minimal such that $i_m\neq 0$; so $i_m\geq 1$.
If $m=n$ then~$y^{\<\i\>}=(y^{\<n\>})^{i_n}\succeq y^{\<n\>}$, hence $y^{\langle \i\rangle}\succeq f$, with $y^{\langle \i\rangle}\succ  f$ if $f\succ x$.
Suppose $m<n$.
Then~$y\succeq\Exp_{m+1}(f^2)$ and hence $y^{\<m\>}\succeq\Exp(f^2)$ by 
Lemma~\ref{lem:expn, 2}. 
Also,   $f\comp f^2\flatter\Exp(f^2)$, thus~$y^{\langle m\rangle}\steeper f$. 
The remarks following Corollary~\ref{cor:itpsi} now yield $y^{\langle \i\rangle}\succ f$.
\end{proof}

\subsection*{Asymptotic behavior of $P(y)$ for large $y$}
In this subsection~$\i$,~$\j$,~$\k$ range over~$\N^{1+n}$.
Let~$P_{\<\i\>}\in K$ be such that $P_{\<\i\>}=0$ for all but finitely many~$\i$ and~${P_{\<\i\>}\neq 0}$ for some~$\i$, and set
$P:=\sum_{\i} P_{\<\i\>} Y^{\<\i\>}\in K\<Y\>$. 
So if
  $P\in K\{Y\}$, then~$P=\sum_{\i} P_{\<\i\>} Y^{\<\i\>}$ is the logarithmic decomposition of the differential polynomial $P$ as defined in
[ADH, 4.2]. 

\begin{example}
$Y=Y^{\<0\>},\quad Y'=Y^{\<0\>}Y^{\<1\>},\quad Y''=Y^{\<0\>}(Y^{\<1\>})^2+Y^{\<0\>}Y^{\<1\>}Y^{\<2\>},$
and for all $m$, 
 $Y^{(m)}\in \Z\big[Y^{\<0\>},Y^{\<1\>},\dots,Y^{\<m\>}\big] $. 
Thus  $P=2Y^3+Y'Y''$ has logarithmic decomposition
$$P\  =\  2(Y^{\<0\>})^3+(Y^{\<0\>})^2(Y^{\<1\>})^3+(Y^{\<0\>})^2(Y^{\<1\>})^2Y^{\<2\>}.$$
\end{example}

\noindent
If $y$ is an element in a differential field extension~$L$ of $K$ such that~$y^{\<n\>}$ is defined, then
we put~$P(y):=\sum_{\i} P_{\<\i\>} y^{\<\i\>}\in L$ (and for~$P\in K\{Y\}$ this has the usual value).
Let $\j$ be lexicographically maximal such that~$P_{\<\j\>}\neq 0$, and 
choose $\k$ so that
$P_{\<\k\>}$ has   minimal valuation.
If $P_{\<\k\>}/P_{\<\j\>} \succ x$,  set
$f:= P_{\<\k\>}/P_{\<\j\>}$; otherwise set $f:=x^2$. Then $f\succ x$ and $f\succeq P_{\<\i\>}/P_{\<\j\>}$ for all $\i$.
The following is a more precise version of~[ADH, 16.6.10] and
\cite[(8.8)]{JvdH}:

\begin{prop}\label{prop:val at infty}
Suppose $\Psi$ is    downward closed, and $y$ in an $H$-asymptotic extension of $K$  satisfies  $y\succeq\Exp_{n}(f^2)$. Then  $y^{\<n\>}$ is defined and~$P(y)\sim P_{\<\j\>}y^{\langle \j\rangle}$.
\end{prop}

\begin{proof}
Let $\i<\j$.
We have   $f\succ x$, so 
$y^{\langle \j-\i\rangle}\succ f \succeq P_{\<\i\>}/P_{\<\j\>}$ by Corollary~\ref{cor:expn}. Hence 
$P_{\<\j\>}y^{\langle \j\rangle} \succ P_{\<\i\>}y^{\langle \i\rangle}$.
\end{proof}

\noindent
From Corollary~\ref{cor:itpsi}, Lemma~\ref{lem:expn, 3}, and Proposition~\ref{prop:val at infty} we obtain:

\begin{cor}\label{cor:val gp at infty}
Suppose $\Psi$ is   downward closed and $y$ in an $H$-asymptotic extension of $K$ satisfies $y\succ K$. Then
  $y$ is $\d$-transcendental over~$K$, and  for all $n$, $y^{\langle n\rangle}$ is defined, $y^{\langle n\rangle}\succ K$, and~$y^{\langle n+1\rangle}\flatter y^{\langle n\rangle}$. The  $H$-asymptotic extension $K\langle y\rangle$ of $K$ has residue field~$\res K\langle y\rangle=\res K$ and value group~$\Gamma_{K\langle y\rangle} = \Gamma \oplus \bigoplus_n \Z v(y^{\langle n\rangle})$ \textup{(}internal direct sum\textup{)}, and $\Gamma_{K\langle y\rangle} $
contains $\Gamma$ as a convex subgroup.
\end{cor}

\noindent
Suppose now that $K$ is equipped with an ordering making it a pre-$H$-field.  From Proposition~\ref{prop:val at infty} we recover~\cite[Theorem~3.4]{AvdD3} in slightly stronger form: 

{\sloppy
\begin{cor}\label{cor:val at infty}
Suppose $y$  lies in
a Liouville closed $H$-field extension of~$K$. If~${y\succeq\Exp_{n}(f^2)}$, then $y^{\<n\>}$ is defined and $\sgn P(y)=\sgn P_{\<\j\>}y^{j_0}$. 
In particular, if~$y^{\<n\>}$ is defined and~$P(y)=0$, then~$y\prec \Exp_{n}(f^2)$.
\end{cor}}

\begin{example}
Suppose $P\in K\{Y\}$. Using [ADH, 4.2, subsection on logarithmic decomposition] we obtain $j_0=\deg P$, and the logarithmic decomposition
$$P(-Y)\ =\ \sum_{\i} P_{\langle\i\rangle} (-1)^{i_0}Y^{\langle\i\rangle}.$$
If $\deg P$ is odd, and  $y>0$ lies in 
a Liouville closed $H$-field extension of $K$ such that~$y\succeq\Exp_{n}(f^2)$, 
then $$\sgn P(y)\ =\ \sgn P_{\<\j\>}, \qquad
\sgn P(-y)\ =\ -\sgn P_{\langle\j\rangle}\ =\ -\sgn P(y).$$
\end{example}

\section{Germs of Continuous Functions}\label{sec:germs}

\noindent
Hardy fields consist of germs of one-variable differentiable real-valued functions. In this section we first consider
the ring $\c$ of germs of {\it continuous}\/ real-valued functions, and its complex counterpart $\c[\imag]$. With an eye towards
applications to Hardy fields, we pay particular attention to extending subfields of $\c$.
 
\subsection*{Germs} As in [ADH, 9.1]  we let $\mathcal{G}$ be the ring
of germs at $+\infty$ of real-valued functions whose domain is
a subset of $\R$
containing an interval~$(a, +\infty)$, $a\in \R$; the domain may vary 
and the ring operations are defined as usual.\index{germ} 
If $g\in \mathcal{G}$ is the germ of a real-valued function 
on a subset of $\R$ containing an interval $(a, +\infty)$, $a\in \R$, then we simplify notation
by letting $g$ also denote this function if the resulting ambiguity is harmless.
With this convention, given a property~$P$ of real numbers
and $g\in \mathcal{G}$ we say that {\em $P\big(g(t)\big)$ holds eventually\/} if~$P\big(g(t)\big)$ holds for all sufficiently large real $t$. 
Thus for $g\in \cal G$ we have $g=0$ iff $g(t)=0$ eventually (and so $g\neq 0$ iff $g(t)\neq 0$ for arbitrarily large $t$). 
Note that the multiplicative group $\mathcal{G}^\times$ of units of $\mathcal{G}$ consists of the $f\in \mathcal{G}$ such that $f(t)\ne 0$, eventually.
We identify each real number~$r$ with the germ at~$+\infty$ of the function~$\R\to \R$ that takes the constant value~$r$. This  makes the field~$\R$ into a subring of $\mathcal{G}$. 
Given $g,h\in\cal G$, we set
\begin{equation}\label{eq:germs partial ordering}
g \leq h \quad :\Longleftrightarrow\quad \text{$g(t)\leq h(t)$, eventually.}
\end{equation}
This defines a partial ordering $\leq$ on $\mathcal G$ which restricts to the usual ordering of $\R$.

Let $g,h\in\mathcal{G}$. Then $g,h\geq 0\Rightarrow g+h,\,g\cdot h,\,g^2\geq 0$, and $g\geq r\in \R^{>}\Rightarrow g\in\mathcal{G}^\times$.
We define $g<h:\Leftrightarrow g\leq h$ and~$g\neq h$. Thus
if $g(t)<h(t)$, eventually, then $g < h$; the converse is not generally valid.

\subsection*{Continuous germs}
We call a germ $g\in \mathcal{G}$ {\em continuous}\/ if it is the germ of
a continuous function $(a,+\infty)\to \R$ for some $a\in \R$, and we let $\mathcal{C}\supseteq \R$ be the subring of $\mathcal{G}$ consisting of the continuous germs $g\in \mathcal{G}$.\index{germ!continuous}\label{p:cont} 
We have $\c^\times=\mathcal{G}^\times\cap\c$;
thus for~$f\in \c^\times$, we have $f(t)\neq 0$, eventually,   hence either $f(t)>0$, eventually, or~$f(t)<0$, eventually, and so $f>0$ or $f<0$.
More generally, if $g,h\in\cal C$ and~$g(t)\neq h(t)$, eventually, then 
 $g(t)<h(t)$, eventually, or  $h(t)<g(t)$, eventually.
 We let $x$ denote the germ at $+\infty$ of the identity function on $\R$, so $x\in\c^\times$.
 
\subsection*{The ring $\c[\imag]$} In analogy with $\c$ we define its complexification $\c[\imag]$ as the ring of
germs at $+\infty$ of $\C$-valued continuous functions whose domain is a subset of $\R$ containing an interval $(a,+\infty)$, $a\in \R$. It has 
$\c$ as a subring. Identifying each complex number $c$ with the germ at $+\infty$ of the function $\R\to \C$ that takes the constant value $c$ makes $\C$ also a subring of $\c[\imag]$ with $\c[\imag]=\c+\c\imag$, justifying the notation $\c[\imag]$. 
The ``eventual'' terminology for germs $f\in \c$ (like ``$f(t)\ne 0$, eventually'') is extended in the obvious way to germs $f\in \c[\imag]$. Thus for $f\in \c[\imag]$ we have:
$f(t)\ne 0$, eventually, if and only if $f\in \c[\imag]^\times$. 
In particular $\c^\times=\c[\imag]^\times\cap \c$.

\medskip
\noindent
Let $\Phi\colon U\to\C$ be a continuous function where $U\subseteq\C$, and let $f\in\c[\imag]$ be such that~$f(t)\in U$, eventually; then $\Phi(f)$ denotes the germ in $\c[\imag]$ with
$\Phi(f)(t)=\Phi\big(f(t)\big)$, eventually.
For example, taking $U=\C$, $\Phi(z)=\ex^z$, we obtain for $f\in\c[\imag]$ the germ  $\exp f = \ex^f\in\c[\imag]$
with $(\ex^f)(t)=\ex^{f(t)}$, eventually. 
Likewise, for~$f\in \c$ with~$f(t)>0$, eventually, we have the germ $\log f\in\c$. 
For $f\in\c[\imag]$  we have~$\overline{f}\in\c[\imag]$ with  $\overline{f}(t)=\overline{f(t)}$, eventually; the map $f\mapsto\overline{f}$ is an automorphism of the ring $\c[\imag]$
with $\overline{\overline{f}}=f$ and~$f\in\c\Leftrightarrow \overline{f}=f$. 
For $f\in\c[\imag]$ we also have~$\Re f,\Im f,\abs{f}\in\c$ with~$f(t)=(\Re f)(t)+(\Im f)(t)\imag$ and $\abs{f}(t)=\abs{f(t)}$, eventually.

\subsection*{Asymptotic relations on $\c[\imag]$.} 
Although $\c[\imag]$ is not a valued field, it will be convenient to equip $\c[\imag]$ with the asymptotic relations $\preceq$,~$\prec$,~$\sim$ (which are defined on any valued field [ADH, 3.1]) as follows: for $f,g\in \c[\imag]$,
\begin{align*} f\preceq g\quad &:\Longleftrightarrow\quad \text{there exists $c\in \R^{>}$ such that $|f|\le c|g|$,}\\
f\prec g\quad &:\Longleftrightarrow\quad \text{$g\in \c[\imag]^\times$ and $\lim_{t\to \infty} f(t)/g(t)=0$} \\
 &\phantom{:} \Longleftrightarrow\quad \text{$g\in \c[\imag]^\times$ and $\abs{f}\leq c\abs{g}$ for all $c\in\R^>$},\\
f\sim g\quad &:\Longleftrightarrow\quad \text{$g\in \c[\imag]^\times$ and
$\lim_{t\to \infty} f(t)/g(t)=1$}\\ 
\quad&\phantom{:} \Longleftrightarrow\quad f-g\prec g.
\end{align*}
If $h\in\c[\imag]$ and $1\preceq h$, then~${h\in\c[\imag]^\times}$. Also, for~$f,g\in\c[\imag]$ and $h\in\c[\imag]^\times$ we have
$$f\preceq g\ \Leftrightarrow\ fh\preceq gh, \qquad f\prec g\ \Leftrightarrow\ fh\prec gh, \qquad f \sim g\ \Leftrightarrow\ fh \sim gh.$$
The binary relation $\preceq$ on
$\c[\imag]$ is reflexive and transitive, and $\sim$ is an equivalence relation on $\c[\imag]^\times$.
Moreover, for $f,g,h\in \c[\imag]$ we have
$$f\prec g\ \Rightarrow\ f\preceq g, \qquad f\preceq g \prec h\ \Rightarrow\ f\prec h, \qquad f\prec g \preceq h\ \Rightarrow\ f\prec h.$$ 
Note that $\prec$ is a transitive binary relation
on $\c[\imag]$.  
For $f,g\in \c[\imag]$ we also set
$$f\asymp g:\ \Leftrightarrow\ f\preceq g \ \&\ g\preceq f,\qquad f\succeq g:\ \Leftrightarrow\ g\preceq f,\qquad f\succ g:\ \Leftrightarrow\ g\prec f,
$$
so $\asymp$ is an equivalence relation on $\c[\imag]$, and $f\sim g\Rightarrow f\asymp g$.  
Thus for $f,g,h\in\c[\imag]$, 
$$f\preceq g\ \Rightarrow\  fh\preceq gh,\quad f\preceq h\ \&\  g\preceq h\ \Rightarrow\ f+g \preceq h,
\quad f\preceq 1\ \&\ g\prec 1 \ \Rightarrow\ fg\prec 1,$$
hence\label{p:contpreceq} 
$$\c[\imag]^{\preceq}\ :=\ \big\{f\in\c[\imag]: f\preceq 1\big\}\ = \ \big\{ f\in\c[\imag]: \text{$|f|\leq n$ for some $n$}\big\}$$
is a subalgebra of the $\C$-algebra $\c[\imag]$ and
$$\c[\imag]^{\prec}\	:=\ \big\{f\in\c[\imag]: f\prec 1\big\}\ 
					 =\ \left\{f\in\c[\imag]: \lim_{t\to\infty} f(t)=0 \right\}$$
is an ideal of $\c[\imag]^{\preceq}$. The group of units of $\c[\imag]^{\preceq}$ is
$$\c[\imag]^{\asymp} \ := \ \big\{ f\in\c[\imag]: f\asymp 1\big\} \ = \ \big\{ f\in\c[\imag]: \text{$1/n\leq |f|\leq n$ for some $n\geq 1$}\big\} $$
and has the subgroup 
$$\C^\times\big(1+\c[\imag]^{\prec}\big)\ =\ \left\{f\in\c[\imag]: \lim_{t\to\infty} f(t)\in\C^\times \right\}.$$ 
We set 
$\c^{\preceq}:=\c[\imag]^{\preceq}\cap\c$, and similarly with~$\prec$,~$\asymp$ in place of $\preceq$.

\begin{lemma}\label{lem:sim props}  
Let $f,g,f^*,g^*\in\c[\imag]^\times$ with $f\sim f^*$ and $g\sim g^*$. Then
$1/f\sim 1/f^*$ and~$fg\sim f^*g^*$. Moreover, $f\preceq g\Leftrightarrow f^*\preceq g^*$, and similarly with
$\prec$, $\asymp$, or $\sim$ in place of~$\preceq$.
\end{lemma}

\noindent
This follows easily from the observations above. For later reference we also note:

\begin{lemma}\label{lem:log preceq} 
Let $f,g\in\c^\times$ be such that $1\prec f\preceq g$; then $\log |f| \preceq \log|g|$.
\end{lemma}
\begin{proof}
Clearly $\log|g| \succ 1$.
Take $c\in\R^>$ such that $|f| \leq c|g|$. Then $\log|f|\leq \log c+\log|g|$ where $\log c+\log|g|\sim\log|g|$;
hence $\log |f| \preceq \log|g|$.
\end{proof}

\begin{lemma}\label{lem:fsimg criterion}  
Let $f,g,h\in\c^\times$ be such that  $f-g\prec h$ and
$(f-h)(g-h)=0$. Then~$f\sim g$.
\end{lemma}
\begin{proof}
Take $a\in\R$ and representatives $(a,+\infty)\to\R$ of $f$, $g$, $h$, denoted by the same symbols, such that
for each $t > a$ we have $f(t),g(t),h(t)\neq 0$, and $f(t)=h(t)$ or $g(t)=h(t)$. Let $\varepsilon\in\R$ with $0<\varepsilon \leq 1$ be given, and choose $b\geq a$ such that~$|f(t)-g(t)| \leq \frac{1}{2}\varepsilon |h(t)|$ for all $t > b$.
Set $q:=f/g$ and
let $t > b$; we claim that then $|q(t)-1|\leq \varepsilon$.
This is clear if $g(t)=h(t)$, so suppose otherwise; then~$f(t)=h(t)$, and  
$|1-1/q(t)| \leq \frac{1}{2}\varepsilon\leq \frac{1}{2}$. In particular, $0<q(t)\leq 2$ and so~$|1-q(t)| = |1-1/q(t)|\cdot q(t) \leq \varepsilon$ as claimed. 
\end{proof}

\subsection*{Subfields of $\mathcal{C}$}
Let $H$ be a {\em Hausdorff field}, that is, a subring
of $\mathcal{C}$ that happens to be a field; see \cite{ADH2}.\index{Hausdorff field} 
Then $H$ has the subfield~$H\cap \R$. 
If $f\in H^\times$, then~$f(t)\neq 0$ eventually, hence
either $f(t)<0$ eventually or $f(t)>0$ eventually. The partial ordering of $\cal G$ from \eqref{eq:germs partial ordering} thus restricts to a total
ordering on $H$ making~$H$ an ordered field in the usual sense of that term. 
By \cite[Propositions~3.4 and~3.6]{Boshernitzan81}: 

\begin{prop}\label{b1} 
Let $H^{\operatorname{rc}}$ consist of the
$y\in \mathcal{C}$ with $P(y)=0$ for some $P\in H[Y]^{\ne}$. Then $H^{\operatorname{rc}}$ is the
unique real closed  Hausdorff field that
extends $H$ and is algebraic over~$H$. In particular, $H^{\operatorname{rc}}$ is a real closure of the ordered field $H$.  
\end{prop} 

\noindent
Boshernitzan~\cite{Boshernitzan81} assumes $H\supseteq \R$ for this result, but this is not really needed in the proof, much of
which already occurs in Hausdorff~\cite{Hau09}.

\medskip
\noindent 
Note that $H[\imag]$ is a subfield of $\c[\imag]$, and by Proposition~\ref{b1} and [ADH, 3.5.4], the subfield~$H^{\operatorname{rc}}[\imag]$ of $\c[\imag]$ is an algebraic closure of the field $H$.
If $f\in\c[\imag]$ is integral over~$H$, then so is $\overline{f}$, and hence so are the elements $\Re f=\frac{1}{2}(f+\overline{f})$ and
$\Im f=\frac{1}{2\imag}(f-\overline{f})$ of $\c$ [ADH, 1.3.2]. Thus
$H^{\operatorname{rc}}[\imag]$  consists of the
$y\in \c[\imag]$ with~$P(y)=0$ for some $P\in H[Y]^{\ne}$. 

\medskip
\noindent
The ordered field $H$ has a convex subring 
$$ \mathcal{O}\ =\ 
\big\{f\in H :\ \text{$|f|\le n$ for some $n$}\big\}\ =\ \c^{\preceq}\cap H,$$
which is a valuation ring of $H$, and we
consider $H$ accordingly as a valued ordered field. We make $\res(H)$ into an ordered field such that $a\ge 0 \Rightarrow a+\smallo\ge 0$ for $a\in \mathcal O$, where $\smallo=\mathcal C^\prec \cap H$ is the maximal ideal of $\mathcal O$. The residue morphism $\mathcal O\to\res(H)$ restricts to an ordered field
embedding $H\cap\R\to \res(H)$, which is bijective if $\R\subseteq H$.
Restricting the binary relations
$\preceq$,~$\prec$,~$\sim$ from the previous subsection to $H$ gives exactly the asymptotic relations $\preceq$,~$\prec$,~$\sim$ on $H$
that it comes equipped with as a valued field. By [ADH, 3.5.15], 
$$\mathcal{O}+\mathcal{O}\imag\ =\ \big\{f\in H[\imag]: \text{$|f|\le n$ for some $n$}\big\}\ =\ \c[\imag]^{\preceq}\cap H[\imag]$$ is the unique valuation ring of $H[\imag]$ whose intersection with $H$ is $\mathcal{O}$. In this way we consider
$H[\imag]$ as a valued field extension of $H$. The maximal ideal of $\mathcal{O}+\mathcal{O}\imag$ is~$\smallo+\smallo \imag=\c[\imag]^{\prec}\cap H[\imag]$.
The asymptotic relations $\preceq$, $\prec$, $\sim$ on $\c[\imag]$
restricted to~$H[\imag]$ are exactly the asymptotic relations $\preceq$, $\prec$, $\sim$ on $H[\imag]$ that $H[\imag]$ has as a valued field, with $f \asymp \abs{f}$ in $\c[\imag]$ for all $f \in H[\imag]$.
In particular, the binary relation $\preceq$ on $\c[\imag]$  restricts to a dominance relation on each subfield of $H[\imag]$ (see [ADH, 3.1.1]).
Let~$K$ be a   subfield of $\c[\imag]$. We note that the following are equivalent:
\begin{enumerate}
\item The binary relation $\preceq$ on $\c[\imag]$ restricts to a dominance relation on $K$;
\item for all $f,g\in K$: $f\preceq g$ or $g\preceq f$;
\item for all $f\in K$: $f\preceq 1$ or $1\preceq f$.
\end{enumerate}
Moreover, the following are equivalent:
\begin{enumerate}
\item $K=H[\imag]$ for some Hausdorff field $H$; 
\item $\imag\in K$ and $\overline{f}\in K$ for each $f\in K$;
\item $\imag\in K$ and $\Re f, \Im f\in K$ for each $f\in K$.
\end{enumerate}

\subsection*{Composition} 
Let $g\in \c$, and suppose that $\lim\limits_{t\to +\infty} g(t)=+\infty$; equivalently,  $g\geq 0$ and $g\succ 1$. Then the 
composition operation
$$f\mapsto f\circ g\ :\ \c[\imag] \to \c[\imag], \qquad (f\circ g)(t)\ :=\  f\big(g(t)\big)\ \text{ eventually},$$
is an injective endomorphism of the ring $\c[\imag]$
that is the identity on the subring~$\C$. For $f_1, f_2\in \c[\imag]$ we have:  $f_1\preceq f_2 \Leftrightarrow f_1\circ g \preceq f_2\circ g$, and likewise with $\prec$, $\sim$. 
This endomorphism of $\c[\imag]$ commutes with the automorphism $f\mapsto\overline{f}$ of $\c[\imag]$, and maps each
subfield $K$ of $\c[\imag]$ isomorphically 
onto the subfield $K\circ g=\{f\circ g:f\in K\}$ of~$\c[\imag]$. Note that if the subfield~$K$ of $\c[\imag]$ contains $x$, then $K\circ g$ contains $g$. 
Moreover,   $f\mapsto f\circ g$ restricts to an endomorphism of the subring $\c$ of~$\c[\imag]$ such that if  $f_1, f_2\in \c$ and $f_1\leq f_2$, then $f_1\circ g \leq f_2\circ g$. This endomorphism of $\c$
maps each Hausdorff field $H$ isomorphically (as an ordered field)
onto the Hausdorff field~$H\circ g$.

\medskip\noindent
Occasionally it is convenient to extend the composition operation on~$\c$ to the ring~$\mathcal G$ of all (not necessarily continuous) germs.
Let~$g\in\mathcal G$ with~$\lim\limits_{t\to+\infty} g(t)=+\infty$. Then for $f\in\mathcal G$ we have the
germ $f\circ g\in\mathcal G$ with  
$$(f\circ g)(t)\ :=\  f\big(g(t)\big)\ \text{ eventually.}$$ 
The map $f\mapsto f\circ g$ is an endomorphism of the $\R$-algebra $\mathcal G$. 
Let $f_1, f_2\in \mathcal{G}$. Then~$f_1\leq f_2\Rightarrow f_1\circ g\leq f_2\circ g$, and likewise with $\preceq$ and $\prec$ instead of $\leq$, where
we extend the binary relations  $\preceq$, $\prec$
from $\c$ to $\mathcal G$ in the natural way: 
\begin{align*} f_1\preceq f_2\quad &:\Longleftrightarrow\quad \text{there exists $c\in \R^{>}$ such that $|f_1(t)|\le c|f_2(t)|$, eventually;}\\
f_1\prec f_2\quad &:\Longleftrightarrow\quad \text{$f_2\in {\mathcal G}^\times$ and $\lim_{t\to \infty} f_1(t)/f_2(t)=0$.}  
\end{align*}

\subsection*{Compositional inversion}
Suppose that $g\in \c$ is eventually strictly increasing such that~$\lim\limits_{t\to +\infty} g(t)=+\infty$.\label{p:ginv}
Then its compositional inverse $g^{\inv}\in \c$ is given by~$g^{\inv}\big(g(t)\big)=t$, eventually,
and  $g^{\inv}$ is also  eventually strictly increasing with~$\lim\limits_{t\to +\infty} g^{\inv}(t)=+\infty$.
Then $f\mapsto f\circ g$  is an automorphism of the ring  $\c[\imag]$, with inverse~$f\mapsto f\circ g^{\inv}$. 
In particular, $g\circ g^{\inv}=g^{\inv}\circ g=x$. Moreover, 
$f\mapsto f\circ g$ restricts to an automorphism of $\c$, and
if  $h\in\c$ is eventually strictly increasing with  $g\leq h$, then~$h^{\operatorname{inv}}\leq g^{\operatorname{inv}}$.

Let now  $f,g\in\c$ with $f,g\geq 0$, $f,g\succ 1$.  
It is not true in general that if~$f$,~$g$ are eventually strictly increasing and $f\sim g$, then $ f^{\operatorname{inv}}\sim g^{\operatorname{inv}}$. (Counterexample:  $f=\log x$, $g=\log 2x$.) Corollary~\ref{cor:Entr} below gives a useful condition on $f$, $g$ under which this implication does hold. 
In addition, let $h\in\c^\times$ be eventually monotone and continuously differentiable with~$h'/h\preceq 1/x$.

\begin{lemma}[Entringer~\cite{Entringer}]\label{lem:Entr}
Suppose ${f\sim g}$.  Then~$h\circ f\sim h\circ g$.
\end{lemma}
\begin{proof}
Replacing $h$ by $-h$ if necessary we arrange that $h\geq 0$, so $h(t)>0$ eventually. 
Set $p:=\min(f,g)\in\c$ and $q:=\max(f,g)\in\c$. Then $0\le p\succ 1$ and $f-g\prec p$.
The Mean Value Theorem gives $\xi\in\mathcal G$ such that 
$p\leq \xi \leq q$ (so $0\le \xi\succ 1$) and
$$h\circ f - h \circ g\  =\  (h'\circ \xi)\cdot (f-g).$$
From  $h'/h\preceq 1/x$ we obtain $h'\circ \xi \preceq ({h\circ \xi})/\xi\preceq (h\circ \xi)/p$, hence
$h\circ f - h \circ g  \prec h\circ \xi$. Set $u:=\max(h\circ p,h\circ q)$. Then 
$0\leq h\circ\xi \leq u$, hence~$h\circ f - h \circ g  \prec u$.
Also~$(u-h\circ f)(u-h\circ g)=0$, so Lemma~\ref{lem:fsimg criterion} yields~$h\circ f  \sim h\circ g$.
\end{proof}

\begin{cor}\label{cor:Entr}
Suppose $f$, $g$ are eventually strictly
increasing such that  $f\sim g$  and~$f^{\operatorname{inv}}\sim h$. Then~$g^{\operatorname{inv}}\sim h$.
\end{cor}
\begin{proof}
By the lemma above we have $h\circ f\sim h\circ g$, and from $f^{\operatorname{inv}}\sim h$
we obtain~$x=f^{\operatorname{inv}}\circ f\sim h\circ f$. Therefore
$g^{\operatorname{inv}}\circ g= x\sim h\circ g$ and thus $g^{\operatorname{inv}}\sim h$.
\end{proof}


\subsection*{Extending ordered fields inside an ambient partially ordered ring}
Let $R$ be a commutative ring with $1\ne 0$, equipped with a translation-invariant
partial ordering $\le$ such that $r^2\ge 0$ for all $r\in R$, and
$rs\ge 0$
for all $r,s\in R$ with $r,s\ge 0$. It follows that for $a,b,r\in R$ we have: 
\begin{enumerate}
\item if $a \le b$ and $r\ge 0$, then $ar\le br$; 
\item if $a$ is a unit and $a>0$, then $a^{-1}=a\cdot (a^{-1})^2>0$; 
\item if $a$,~$b$ are units and $0 < a \le b$, then $0 < b^{-1}\le a^{-1}$. 
\end{enumerate}
Relevant cases: 
$R=\mathcal{G}$ and $R=\mathcal{C}$, with partial ordering given by \eqref{eq:germs partial ordering}. 

\medskip
\noindent
An {\em ordered subring of $R$} is a subring  
of $R$ that is totally ordered by the partial ordering of $R$. An {\em ordered subfield of $R$} is an ordered subring $H$ of $R$
which happens to be a field; then $H$ equipped with the induced ordering
is indeed an ordered field, in the usual sense of that term.
(Thus any Hausdorff field is an ordered subfield of the partially ordered ring $\mathcal{C}$.) 
We identify $\Z$ with its image in $R$ via the unique ring embedding $\Z \to R$, and this makes $\Z$ with its usual ordering into an ordered subring of $R$.

\begin{lemma}\label{pr0} Assume $D$ is an ordered
subring of $R$ and every nonzero element of~$D$ is a unit of $R$.
Then $D$ generates an ordered subfield $\Frac{D}$ of $R$.
\end{lemma}
\begin{proof} It is clear that $D$ generates a subfield $\Frac{D}$ of $R$. For $a\in D$, $a>0$, we have $a^{-1}>0$. It follows that $\Frac{D}$ is totally ordered.
\end{proof}

\noindent
Thus if every $n\ge 1$ is a unit of $R$, then we may identify $\Q$ with its image in $R$ via the unique
ring embedding $\Q\to R$, making $\Q$ into an ordered subfield of $R$.

\begin{lemma}\label{pr2} Suppose $H$ is an ordered subfield of $R$,
all $g\in R$ with $g>H$ are units of $R$, and $H<f\in R$. Then we have an ordered subfield $H(f)$ of $R$. 
\end{lemma}
\begin{proof} 
For $P\in H[Y]$ of degree $d\ge 1$
with leading coefficient $a>0$ we have~$P(f)=af^d(1+\epsilon)$
with $-1/n < \epsilon < 1/n$ for all $n\ge 1$, in particular,
$P(f)>H$ is a unit of~$R$. It remains to appeal to Lemma~\ref{pr0}.
\end{proof}

\begin{lemma}\label{pr1} Let $H$ be a real closed ordered subfield of $R$.
Let $A$ be a nonempty downward closed subset of $H$ such that
$A$ has no largest element and $B:=H\setminus A$ is nonempty and
has no least element. Let $f\in R$ be such that $A<f<B$. Then the subring $H[f]$ of $R$ has the following properties: \begin{enumerate}
\item[$\rm{(i)}$] $H[f]$ is a domain;
\item[$\rm{(ii)}$] $H[f]$ is an ordered subring of $R$;
\item[$\rm{(iii)}$] $H$ is cofinal in $H[f]$;
\item[$\rm{(iv)}$] for all $g\in H[f]\setminus H$ and $a\in H$, if $a<g$, then $a < b<g$ for some $b\in H$, and if $g<a$, then 
$g<b<a$ for some $b\in H$. 
\end{enumerate}
\end{lemma} 
{\sloppy
\begin{proof} Let
$P\in H[Y]^{\neq}$; to obtain (i) and (ii) it suffices to show that then
$P(f) < 0$ or $P(f)>0$.  We have $$P \ =\ c\,Q\,(Y-a_1)\cdots(Y-a_n)$$ where
$c\in H^{\ne}$,  $Q$  is  a product of monic quadratic irreducibles in $H[Y]$, and~$a_1,\dots, a_n\in H$. This gives $\delta\in H^{>}$ such that $Q(r)\ge \delta$ for all $r\in R$.  Assume~$c>0$. (The case $c<0$ is handled similarly.) We can arrange that $m\le n$ is such that~$a_i\in A$ for $1\le i \le m$ and $a_j\in B$ for $m < j\le n$.
Take $\epsilon>0$ in $H$ such that
$a_i + \epsilon \le f$ for~$1\le i\le m$ and 
$f\le a_j-\epsilon$ for~$m < j\le n$.
Then $$P(f)\ =\ c\,Q(f)\,(f-a_1)\cdots (f-a_m)(f-a_{m+1})\cdots(f-a_n),$$
and $(f-a_1)\cdots(f-a_m) \ge \epsilon^m$. If
$n-m$ is even, then $(f-a_{m+1})\cdots(f-a_n)\ge \epsilon^{n-m}$, so $P(f)\ge c\delta\epsilon^n >0$. If $n-m$ is odd, then
$(f-a_{m+1})\cdots(f-a_n)\le -\epsilon^{n-m}$,
so~$P(f) \le -c\delta\epsilon^n < 0$. These estimates also yield (iii) and (iv). 
\end{proof}}

\begin{lemma}\label{pr3} With $H$,~$A$,~$f$ as in Lemma~\ref{pr1}, suppose all $g\in R$ with $g\ge 1$ are units of $R$. Then we have an ordered subfield $H(f)$ of $R$ such that {\rm(iii)} and {\rm(iv)} of Lemma~\ref{pr1} go through for $H(f)$ in place of $H[f]$.
\end{lemma} 
\begin{proof} 
Note that if $g\in R$ and $g\ge \delta\in H^{>}$, then 
$g\delta^{-1}\ge 1$, so $g$ is a unit of $R$ and
$0<g^{-1}\le \delta^{-1}$. For   
$Q\in H[Y]^{\ne}$ with $Q(f)>0$ we can take $\delta\in H^{>}$ such that~$Q(f)\ge \delta$, so $Q(f)\in R^\times$  and $0 <Q(f)^{-1}\le \delta^{-1}$. Thus we have an
ordered subfield
$H(f)$ of $R$ by Lemma~\ref{pr0}, and 
the rest now follows easily.  
\end{proof}

\subsection*{Adjoining pseudolimits and increasing the value group} Let $H$ be a real closed Hausdorff field and  view $H$ as a valued ordered field as
before. Let $(a_{\rho})$ be a strictly increasing divergent pc-sequence in $H$.
Set 
$$A\ :=\ \{a\in H:\ \text{$a< a_{\rho}$ for some $\rho$}\}, \qquad
B\ :=\ \{b\in H:\ \text{$b>a_{\rho}$ for all $\rho$}\},$$ so
$A$ is nonempty and downward closed without a largest element.
Moreover,~$B=H\setminus A$ is nonempty and has no least element, since a
least element of $B$ would be a limit and thus a pseudolimit of
$(a_{\rho})$. Let $f\in \mathcal{C}$ satisfy $A<f<B$.
Then by Lem\-ma~\ref{pr3} for $R=\c$ we have an ordered subfield $H(f)$ of $\mathcal{C}$, and: 

\begin{lemma} \label{ps1} $H(f)$ is an immediate valued field extension of $H$ with $a_{\rho} \leadsto f$.
\end{lemma}
\begin{proof} We can assume that $v(a_{\tau}-a_{\sigma})>v(a_{\sigma}-a_{\rho})$ for all indices $\tau>\sigma>\rho$. Set~$d_{\rho}:= a_{s(\rho)}-a_{\rho}$ ($s(\rho):=$~successor of $\rho$).
Then $a_{\rho}+2d_{\rho}\in B$ for all indices~$\rho$; see the discussion
preceding~[ADH,~2.4.2]. It then follows from that fact that~$a_{\rho} \leadsto f$. Now $(a_{\rho})$ is a divergent pc-sequence in the henselian valued field $H$, so it is of transcendental type over $H$, and thus $H(f)$ is an immediate
extension of~$H$.  
\end{proof}

\begin{lemma}\label{ps2} Let $H$ be a Hausdorff field with divisible value group $\Gamma:=v(H^\times)$. Let $P$ be a nonempty upward closed subset of $\Gamma$, and let $f\in \mathcal{C}$ be such that 
$a < f$ for all~$a\in H^{>}$ with 
$va\in P$, and $f<b$ for all $b\in H^{>}$ with $vb < P$.
Then $f$ generates a Hausdorff field~$H(f)$ such that 
$P >vf > Q$ where $Q:= \Gamma\setminus P$, and $f$ is transcendental over~$H$.  
\end{lemma}
\begin{proof} For any positive $a\in H^{\text{rc}}$ there is
$b\in H^{>}$ with $a\asymp b$ and $a< b$, and also an element
$b\in H^{>}$ with $a\asymp b$ and $a>b$. Thus by Proposition~\ref{b1} we can replace~$H$ by~$H^{\text{rc}}$ and arrange in this way that $H$ is real closed.
 Set 
$$A\ :=\ \{a\in H:\ \text{$a\le 0$  or $va\in P$}\}, \qquad B:= H\setminus A.$$
Then we are in the situation of Lemma~\ref{pr1} for $R=\mathcal{C}$,   so by that lemma and Lemma~\ref{pr3} we have a Hausdorff field
$H(f)$. Clearly then $P> vf >Q$. In particular, $f\notin H$, so   $f$ is transcendental over~$H$.  
\end{proof}

\subsection*{Non-oscillation} A germ $f\in \c$ is said to {\bf oscillate\/}  if $f(t)=0$ for arbitrarily large~$t$ and $f(t)\ne 0$ for
arbitrarily large $t$. \index{germ!oscillating} Thus for $f,g\in \c$,
$$ f-g \text{  is non-oscillating }\quad \Longleftrightarrow\quad \begin{cases} &\parbox{17em}{ either $f(t)< g(t)$ eventually, or $f=g$, or $f(t)>g(t)$ eventually.}\end{cases}$$
In particular, $f\in\c$ does not oscillate iff $f=0$ or $f\in\c^\times$.
If $g\in\c$ and~$g(t)\to +\infty$ as $t\to+\infty$, then $f\in\c$ oscillates iff $f\circ g$ oscillates.  
The following two lemmas are included for use in a later paper:

\begin{lemma}\label{lem:nonosc, 1}
Let $f\in\c$ be such that for every $q\in\Q$ the germ $f-q$ is non-oscillating. 
Then~$\lim\limits_{t\to\infty} f(t)$ exists in $\R\cup\{-\infty,+\infty\}$.
\end{lemma}
\begin{proof} Set $S\, :=\, \{s\in\Q: f(t)> s \text{ eventually}\}$.
If $S=\emptyset$, then $\lim\limits_{t\to\infty} f(t)=-\infty$, whereas if $S=\Q$, then 
$\lim\limits_{t\to\infty} f(t)=+\infty$.
If $S\neq\emptyset,\Q$, then for $\ell:=\sup S\in\R$ we have
$\lim\limits_{t\to\infty} f(t)=\ell$.
\end{proof}

\begin{lemma}\label{lem:nonosc, 2}
Let $H$ be a real closed Hausdorff field and $f\in\c$. Then $f$ lies in a Hausdorff field extension of $H$ iff 
$f-h$ is non-oscillating for all $h\in H$. 
\end{lemma}
\begin{proof}
The forward direction is clear. For the converse, suppose $f-h$ is non-oscillating  for all $h\in H$.
We assume $f\notin H$, so $h<f$ or $h>f$ for all $h\in H$. 
Set~${A:=\{h\in H:h < f\}}$, a downward closed subset of $H$.
If $A=H$, then we are done by Lemma~\ref{pr2} applied to $R=\c$; if $A=\emptyset$ then 
we apply the same lemma to~$R=\c$ and $-f$ in place of $f$.
Suppose $A\neq\emptyset,H$.
If~$A$ has a largest element $a$, then we replace $f$ by~$f-a$ to arrange $0 < f (t)< h(t)$ eventually, for all $h\in H^>$, and
then  Lemma~\ref{pr2} applied to $R=\c$, $f^{-1}$ in place of $f$ yields that $f^{-1}$, and hence also $f$, lies in a Hausdorff field
extension of $H$. The case that $B:=H\setminus A$ has a least element is handled in the same way. 
If $A$ has no largest element and $B$ has no least element, then
we are done by Lemma~\ref{pr3}.
\end{proof}

\section{Germs of Differentiable Functions} \label{sec:differentiable germs}

\noindent
In this section we fix notations and conventions concerning differentiable functions and summarize well-known
results on second-order linear differential equations as needed later. 
(Basic facts about linear  differential equations can be found in  \cite[Ch.~X]{Dieudonne}, \cite[Ch.~XI]{Hartman}, and \cite[Ch.~IV]{Walter}.)

\subsection*{Differentiable functions} 
Let $r$ range over $\N\cup\{\infty\}$, and let $U$ be a nonempty open subset of $\R$. \label{p:Cr(U)}
Then $\c^r(U)$ denotes the $\R$-algebra of $r$-times
continuously differentiable functions $U\to \R$, with the usual pointwise defined algebra operations. (We use ``$\c$'' instead of ``$C$'' since $C$ will often denote the constant field of a
differential field.) For $r=0$ this is the $\R$-algebra~$\c(U)$ of continuous real-valued functions on~$U$, so 
$$\c(U)\ =\ \c^0(U)\ \supseteq\ \c^1(U)\ \supseteq\ \c^2(U)\ \supseteq\ 
\cdots\ \supseteq\ \c^{\infty}(U).$$ For $r\ge 1$ we have the derivation $f\mapsto f'\colon \c^r(U)\to \c^{r-1}(U)$ (with $\infty-1:=\infty$). This makes $\c^{\infty}(U)$ a differential ring, with
its subalgebra
$\c^{\omega}(U)$ 
of real-analytic functions $U\to \R$ as a differential subring. The algebra operations
on the algebras below are also defined pointwise.  
Note that
$$\c^r(U)^\times\  =\  \big\{ f\in\c^r(U):\text{$f(t)\neq 0$ for all $t\in U$}\big\},$$ 
also for $\omega$ in place of $r$ \cite[(9.2), ex.~4]{Dieudonne}. 

{\sloppy
\medskip
\noindent
Let $a$ range over $\R$.  
Then $\Car$ denotes the $\R$-algebra of functions $[a,+\infty) \to \R$ that extend to a function in $\c^r(U)$ for some open $U\supseteq [a,+\infty)$. Thus $\Caz$ (also denoted by $\c_a$) is the $\R$-algebra of real-valued continuous functions on $[a,+\infty)$, and
$$\Caz\ \supseteq\ \Cao\ \supseteq\ \Cat\ \supseteq\ \cdots\ \supseteq \Cainf.$$  
We have the subalgebra $\Caom$ of $\Cainf$, consisting of the functions~${[a,+\infty)\to \R}$ that extend to a real-analytic function
$U \to \R$ for some open $U\supseteq [a,+\infty)$. 
For~${f\in\Cao}$ and~$g\in\c^1(U)$ extending $f$  
with open $U\subseteq \R$ containing $[a,+\infty)$, the restriction of~$g'$ to $[a,+\infty)\to\R$ depends only on $f$, not on $g$, and so
we may de\-fine~${f':=g'|_{[a,+\infty)}\in\c_a}$.
For~$r\ge 1$ this gives the derivation $f\mapsto f'\colon \Car\to \Carl$. This makes $\Cainf$ a differential ring with~$\c^{\omega}_a$  as a differential subring.
} 

\medskip
\noindent
For each of the algebras $A$ above we also consider its complexification $A[\imag]$ which consists by definition of the
$\C$-valued functions $f=g+h\imag$ with $g,h\in A$, so
$g=\Re f$ and $h= \Im f$ for such $f$. We consider $A[\imag]$ as a $\C$-algebra with respect to the natural pointwise defined algebra operations. We identify each complex number with
the corresponding constant function to make $\C$ a subfield of
$A[\imag]$ and~$\R$ a subfield of~$A$. (This justifies the notation $A[\imag]$.) 
We have $\Car[\imag]^\times=\c_a[\imag]^\times\cap\Car[\imag]$ and
$(\Car)^\times=\c_a^\times\cap\Car$, and likewise with $r$ replaced by $\omega$.

\medskip
\noindent
For $r\ge 1$ we extend $g\mapsto g'\colon \Car\to \Carl$ to the derivation 
$$g+h\imag\mapsto g'+h'\imag\ :\  \Car[\imag] \to 
\Carl[\imag] \qquad (g,h\in \Car[\imag]),$$
which for $r=\infty$ makes $\Cainf$ a differential subring of $\Cainf[\imag]$.  
We shall use the map
$$f \mapsto f^\dagger:=f'/f\ \colon \ \Cao[\imag]^\times=\big(\Cao[\imag]\big)\!^\times \to \Caz[\imag],$$
with 
$$(fg)^\dagger=f^\dagger+g^\dagger\qquad\text{ for $f,g\in \Cao[\imag]^\times$,}$$
in particular the fact that $f\in \Cao[\imag]^\times$ 
and $f^\dagger\in \Caz[\imag]$ are related by
$$ f(t)\ =\ f(a)\exp\!\left[\int_a^t f^\dagger(s)\,ds\right] \qquad (t\ge a).$$ 
For $g\in \Caz[\imag]$, let $\exp\int g$  denote the function   $t\mapsto \exp\!\left[\int_a^tg(s)\,ds\right]$ in $\Cao[\imag]^\times$.
Then
$$(\exp\textstyle\int g)^\dagger=g\quad\text{ and }\quad \exp\int (g+h)=(\exp\int g)\cdot (\exp\int h)\qquad\text{ for $g,h\in \Caz[\imag]$.}$$ 
Therefore~$f\mapsto f^\dagger\colon \Cao[\imag]^\times\to \Caz[\imag]$ is surjective. 

\begin{notation}
For $b\geq a$ and $f\in \c_a[\imag]$ we set  $f|_b:=f|_{[b,+\infty)}\in \c_a[\imag]$. 
\end{notation}

\subsection*{Differentiable germs} 
Let $r\in \N\cup\{\infty\}$ and let $a$ range over $\R$.\index{germ!differentiable}\label{p:Gr}
Let $\Gr$ be the partially ordered subring of 
$\c$ consisting of the germs at $+\infty$ of the functions in~$\bigcup_a \Car$; thus $\Gz=\c$ consists of the germs at~$+\infty$ of the continuous real-valued functions on intervals $[a,+\infty)$, $a\in \R$. Note that $\Gr$ with its partial ordering satisfies the conditions on $R$ from  Section~\ref{sec:germs}. Also, every~$g\ge 1$ in~$\Gr$ is a unit of $\Gr$, so
Lemmas~\ref{pr2} and~\ref{pr3} apply to ordered subfields of $\Gr$. We have 
$$\Gz\ \supseteq\ \Go\ \supseteq\ \Gt\ \supseteq\ \cdots\ \supseteq\ \Ginf.$$
Each subring $\Gr$ of $\c$ yields the subring $\Gr[\imag]=\Gr+\Gr\imag$ of~$\Gz[\imag]=\c[\imag]$, with
$$\Gz[\imag]\ \supseteq\ \Go[\imag]\ \supseteq\ \Gt[\imag]\ \supseteq\ \cdots\ \supseteq\ \Ginf[\imag].$$
Suppose $r\geq 1$; then for $f\in\Car[\imag]$ the germ of $f'\in\Carl[\imag]$ only depends on the germ of $f$, and we thus obtain a derivation $g\mapsto g'\colon \Gr[\imag]\to\c^{r-1}[\imag]$
with $\text{(germ of $f$)}'=\text{(germ of $f'$)}$ for $f\in\bigcup_a \Car[\imag]$.
This derivation restricts to a derivation  $\Gr\to\c^{r-1}$. 
Note that
$\c[\imag]^\times\cap \c^r[\imag] = \c^r[\imag]^\times$, and hence 
$\c^\times\cap \c^r = (\c^r)^\times$. 
Given $g\in\c^r$ with~${g(t)\to+\infty}$ as $t\to+\infty$ and $f\in\c^r[\imag]$,   the germ $f\circ g$ (as defined  in Section~\ref{sec:germs}) also lies
in $\c^r[\imag]$, with $f\circ g\in\c^r$ if $f\in\c^r$.   

\medskip
\noindent
We set 
$$\Gi[\imag]\ :=\ \bigcap_{n}\, \Gn[\imag].$$ Thus $\Gi[\imag]$ is naturally a
differential ring with $\C$ as its ring of constants. \label{p:Gi} 
We also have the differential subring
$$\Gi \ :=\ \bigcap_{n}\, \Gn$$ of $\Gi[\imag]$, 
with $\R$ as its ring of constants and $\Gi[\imag]=\Gi+\Gi\imag$.  
Note that~$\Gi[\imag]$ has~$\Ginf[\imag]$ as a differential subring. Similarly, $\Gi$ has $\Ginf$ as a differential subring,
and the differential ring $\Ginf$ has in turn
the differential subring~$\Gom$, whose elements are the
germs at~$+\infty$ of the functions in $\bigcup_a \Caom$.\index{germ!analytic}\index{germ!smooth} 
We have~$\c[\imag]^\times\cap\Calinf[\imag]=(\Calinf[\imag])^\times$ and  
$\c^\times  \cap \Calinf = (\Calinf)^\times$, and likewise with~$\Gom$ in place of $\Calinf$.  If $R$ is a subring of~$\Go$ such that $f'\in R$ for all~$f\in R$, then~$R\subseteq\Gi$ is a differential subring of~$\Gi$.

\subsection*{Compositional  conjugation of differentiable germs}  
Let~$\ell\in \Go$, $\ell'(t)>0$ eventually (so $\ell$ is eventually strictly increasing) and $\ell(t)\to+\infty$ as $t\to+\infty$. Then~$\phi:=\ell'\in\c^\times$,  and the compositional inverse 
$\ell^{\inv}\in \Go$ of $\ell$ satisfies 
$$\ell^{\inv}>\R, \qquad (\ell^{\inv})'\ =\ (1/\phi)\circ \ell^{\inv}\in \c.$$
The $\C$-algebra automorphism\label{p:fcirc}
$f\mapsto f^\circ:= f\circ \ell^{\inv}$ of $\mathcal{C}[\imag]$ (with inverse $g\mapsto g\circ\ell$) maps $\Go[\imag]$ onto itself and  
satisfies for $f\in\Go[\imag]$ a useful identity:
$$(f^\circ)'\ =\ (f\circ \ell^{\inv})'\ =\ (f'\circ \ell^{\inv})\cdot (\ell^{\inv})'\ =\ (f'/\ell')\circ \ell^{\inv}
\ =\ (\phi^{-1}f')^\circ.$$
Hence if $n\ge 1$ and $\ell\in\c^n$, then $\ell^{\inv}\in \c^n$ and $f\mapsto f^\circ$ maps~$\Caln[\imag]$ and $\Caln$
onto themselves, for each $n$. Therefore, if $\ell\in\Calinf$, then $\ell^{\inv}\in \Calinf$ and $f\mapsto f^\circ$   maps~$\Calinf[\imag]$ and $\Calinf$ onto themselves;
likewise with $\Ginf$ or $\Gom$ in place of $\Calinf$. 
{\it In the rest of this subsection we assume~$\ell\in\Calinf$.}\/ Denote the differential ring~$\Calinf[\imag]$ by~$R$, and as usual let $R^\phi$ be  $R$  with its derivation~$f\mapsto \der(f)=f'$ replaced by the  derivation~$f\mapsto \derdelta(f)=\phi^{-1}f'$~[ADH, 5.7]. Then  
$f\mapsto f^\circ \colon R^\phi \to R$
is an isomorphism of differential rings by the identity above. 
We extend it to the isomorphism
$$Q\mapsto Q^\circ\ \colon\ R^\phi\{Y\} \to R\{Y\}$$
of differential rings given by $Y^\circ=Y$.  
 Let $y\in R$. Then 
$$Q(y)^\circ \  = \  Q^\circ(y^\circ)\qquad\text{for $Q\in R^\phi\{Y\}$}$$
and thus
$$P(y)^\circ\ =\ P^\phi(y)^\circ\ =\ (P^\phi)^\circ(y^\circ) \qquad\text{for $P\in R\{Y\}$.}$$

\subsection*{Second-order differential equations} Let $f\in \c_a$, that is, $f\colon [a,\infty) \to \R$ is continuous. We consider the differential equation
\begin{equation}\label{eq:2nd order}\tag{L}
 Y'' + fY\ =\ 0. 
 \end{equation}
The solutions~${y\in \Cat}$ of \eqref{eq:2nd order} form
an $\R$-linear subspace $\Sol(f)$ of $\Cat$.   The
solutions~${y\in \Cat[\imag]}$ of \eqref{eq:2nd order} are the $y_1+y_2\imag$ with
$y_1, y_2\in \Sol(f)$ and form
a $\C$-linear subspace $\Sol_{\C}(f)$ of~$\Cat[\imag]$. For any complex numbers
$c$,~$d$ there is a unique solution~$y\in \Cat[\imag]$ of   \eqref{eq:2nd order} with~${y(a)=c}$ and $y'(a)=d$, and the map that assigns to $(c,d)$ in $\C^2$
this unique solution is an isomorphism $\C^2\to \Sol_{\C}(f)$ of 
$\C$-linear spaces; it restricts to an $\R$-linear bijection~${\R^2\to \Sol(f)}$. 
We have $f\in \Can\Rightarrow\Sol(f)\subseteq \c^{n+2}_a$ (hence~$f\in \c^\infty_a\Rightarrow {\Sol(f)\subseteq \c^\infty_a}$)  and $f\in \Caom\Rightarrow \Sol(f)\subseteq \Caom$. (See \cite[(10.5.3)]{Dieudonne}.)    

\medskip
\noindent
Let $y_1,y_2\in\Sol(f)$, with 
Wronskian $w=y_1y_2'-y_1'y_2$. Then $w\in\R$, and
$$w\neq 0\ \Longleftrightarrow\ \text{$y_1$,~$y_2$ are $\R$-linearly independent}.$$ 
By~\cite[Chapter~6, Lemmas~2 and 3]{Bellman} we have:

\begin{lemma}\label{lem:2nd order inhom}
Let $y_1,y_2\in\Sol(f)$ be $\R$-linearly independent and $g\in\c_a$. Then
$$t\mapsto y(t)\ :=\ -y_1(t)\int_a^t \frac{y_2(s)}{w}g(s)\,ds + y_2(t)\int_a^t \frac{y_1(s)}{w} g(s)\,ds\ \colon\  [a,+\infty)\to\R$$
lies in $\Cat$ and satisfies $y''+fy=g$, $y(a)=y'(a)=0$. 
\end{lemma}

\begin{lemma}\label{lem:2nd lin indep sol}
Let $y_1\in\Sol(f)$ with $y_1(t)\neq 0$ for $t\geq a$. Then the function
$$t\mapsto y_2(t):=y_1(t)\int_a^t \frac{1}{y_1(s)^2}\,ds\colon [a,+\infty)\to\R$$
also lies in $\Sol(f)$, and $y_1$, $y_2$ are $\R$-linearly independent.
\end{lemma}

\noindent
From~\cite[Chapter~2, Lemma~1]{Bellman} we also recall:

\begin{lemma}[Gronwall's Lemma] \label{lem:gronwall}
Let   $C\in\R^{\ge}$,~$v,y\in\mathcal C_a$ satisfy~$v(t),y(t)\geq 0$ for all $t\geq a$ and
$$y(t)\ \leq\ C+ \int_a^t v(s)y(s)\,ds\quad\text{for all $t\geq a$.}$$
Then
$$y(t)\ \leq\ C\exp\!\left[ \int_a^t v(s)\,ds\right]\quad\text{for all $t\geq a$.}$$
\end{lemma}

\noindent
{\it In the rest of this subsection we assume that
$a\geq 1$ and that $c\in\R^>$  is such that~$|f(t)| \leq c/t^2$ for all $t\geq a$.}
Under this hypothesis, Lemma~\ref{lem:gronwall} yields the following bound on the growth of the solutions $y\in\Sol(f)$; the proof we give is similar to that of~\cite[Chap\-ter~6, Theorem~5]{Bellman}.

\begin{prop}\label{prop:bound} Let $y\in\Sol(f)$. Then there is $C\in\R^\geq$ such that $|y(t)| \leq Ct^{c+1}$ 
and $|y'(t)|\leq Ct^c$
for all $t\geq a$.
\end{prop}
\begin{proof}
Let $t$ range over $[a,+\infty)$.
Integrating $y''=-fy$ twice between $a$ and $t$, we obtain constants $c_1$, $c_2$ such that
for all $t$,
$$y(t)\ =\ c_1+c_2t -\int_a^t\int_a^{t_1} f(t_2)y(t_2)\,dt_2\,dt_1\ =\ c_1+c_2t-\int_a^t (t-s)f(s)y(s)\, ds$$
and hence, with $C:=|c_1|+|c_2|$,
$$ |y(t)|\ \leq\ Ct + t\int_a^t |f(s)|\cdot|y(s)|\, ds,\ 
\text{ so }\ 
\frac{|y(t)|}{t}\ \leq\ C + \int_a^t s|f(s)|\cdot\frac{|y(s)|}{s}\, ds.$$
Then by Lemma~\ref{lem:gronwall},  
$$\frac{|y(t)|}{t}\ \leq\ C\exp\!\left[\int_a^t s|f(s)|\,ds\right]\ \leq\
C\exp\!\left[\int_1^t c/s\,ds\right]\ =\ Ct^c$$
and thus
$|y(t)|\leq Ct^{c+1}$. Now  
\begin{align*} y'(t)\ &=\ c_2-\int_a^t f(s)y(s)\,ds, \text{ so}\\
|y'(t)|\ &\leq\ |c_2|+\int_a^t |f(s)y(s)|\,ds\ \leq\ C+Cc\int_1^t s^{c-1}\,ds\\ &=\
C+Cc\left[\frac{t^c}{c}-\frac{1}{c}\right]\ =\ Ct^c. \qedhere
\end{align*}
\end{proof} 

\noindent
Let $y_1, y_2\in \Sol(f)$ be $\R$-linearly independent. Recall that $w=y_1y_2'-y_1'y_2\in\R^\times$. It follows that $y_1$ and~$y_2$ cannot be simultaneously very small:

\begin{lemma}\label{lem:bound} There is a positive constant $d$ such that
$$\max\!\big(|y_1(t)|, |y_2(t)|\big)\ \ge\ dt^{-c} \quad \text{ for all $t\ge a$.}$$
\end{lemma}
\begin{proof} 
Proposition~\ref{prop:bound} yields $C\in \R^>$ such that  $|y_i'(t)|\le Ct^c$ for
$i=1,2$ and all~$t\ge a$. Hence $|w|\le 2\max\!\big(|y_1(t)|, |y_2(t)|\big)Ct^c$ for $t\ge a$, so
\equationqed{ \max\!\big(|y_1(t)|, |y_2(t)|\big)\ \ge\ \frac{|w|}{2C}t^{-c}\qquad(t\ge a).}
\end{proof}

\begin{cor}\label{cor:bound}
Set $y:=y_1+y_2\imag$ and $z:=y^\dagger$. Then for some $D\in\R^>$,  
$$|z(t)|\ \leq\ Dt^{2c}\quad\text{ for all $t\geq a$.}$$
\end{cor}
\begin{proof}
Take $C$ as in the proof of Lemma~\ref{lem:bound}, and $d$ as in that lemma. Then
$$|z(t)|\ =\ \frac{|y_1'(t)+y_2'(t)\imag|}{|y_1(t)+y_2(t)\imag|}\ \leq\  \frac{|y_1'(t)|+|y_2'(t)|}{\max\!\big(|y_1(t)|, |y_2(t)|\big)}\ \leq\ \left(\frac{2C}{d}\right)t^{2c}$$
for $t\geq a$.
\end{proof}

\subsection*{Oscillation} 
Let  $y\in\c_a$. 
We say that $y$ {\bf oscillates} if its germ in $\c$ oscillates. So~$y$ does not oscillate iff~$\sgn y(t)$ is  constant, eventually.
If $y$ oscillates, then so does~$cy$ for $c\in\R^\times$.
If $y\in\Cao$ oscillates, then so does $y'\in\c_a$, by Rolle's Theorem.

We now continue with the study of \eqref{eq:2nd order}.
Let 
$y\in \Sol(f)^{\neq}$, and  let $Z:=y^{-1}(0)$ be the set of zeros of~$y$, so $Z\subseteq [a,+\infty)$ is closed in $\R$.  
Moreover, $Z$ has no limit points: for all $t_0<t_1$ in $Z$ there is an $s\in(t_0,t_1)$ with $y'(s)=0$  (by Rolle), hence
if $t$ is a limit point of $Z$, then $t\ge a$ and $y(t)=y'(t)=0$, so~$y=0$, a contradiction.
In particular, $Z\cap[a,b]$  is finite for every $b\ge a$. Thus
$$\text{$y$ does not oscillate}
\quad \Longleftrightarrow\quad   \text{$Z$ is finite}\quad \Longleftrightarrow\quad \text{$Z$ is bounded.}$$
If $t_0\in Z$    is not the largest element of $Z$, then $Z\cap (t_0,t_1)=\emptyset$ for some $t_1> t_0$ in~$Z$.
We say that a pair of zeros~$t_0<t_1$ of $y$ is {\bf consecutive} if $Z\cap(t_0,t_1)=\emptyset$. Sturm's Separation Theorem says that if $y,z\in\Sol(f)$ are $\R$-linearly independent and   $t_0<t_1$ are consecutive zeros of $z$,  then $(t_0,t_1)$ contains a unique 
zero  of $y$ \cite[\S{}27, VI]{Walter}.  Thus:

\begin{lemma}\label{lem:some vs all oscillate}
Some $y\in\Sol(f)^{\neq}$ oscillates  $\ \Longleftrightarrow\ $ every $y\in\Sol(f)^{\neq}$ oscillates.
\end{lemma}

\noindent
We say that {\bf $f$ generates oscillation} if some element of $\Sol(f)^{\ne}$  oscillates.\index{germ!generates oscillations}\index{generates oscillations!germ}

\begin{lemma}
Let $b\in\R^{\geq a}$. Then 
$$\text{$f$ generates oscillation}\quad\Longleftrightarrow\quad\text{$f|_{b}\in\c_b$ generates oscillation.}$$
\end{lemma}
\begin{proof}
The forward direction is obvious. For the backward direction, use that every~$y\in\c^2_b$ with $y''+gy=0$ for $g:=f|_{b}$
extends uniquely to a solution of~\eqref{eq:2nd order}.
\end{proof}

\noindent
By this lemma, whether $f$ generates oscillation depends only on its germ in $\c$. So this induces the notion of an element of $\c$ generating oscillation. 
Here is another result of Sturm~\cite[loc.~cit.]{Walter} that we use below:

\begin{theorem}[Sturm's Comparison Theorem] \label{thm:Sturm Comp}
Let $g\in\c_a$ with $f(t)\geq g(t)$ for all~$t\geq a$. Let $y\in\Sol(f)^{\ne}$ and $z\in\Sol(g)^{\neq}$, and let $t_0<t_1$ be consecutive zeros of~$z$. Then either $(t_0,t_1)$ contains a zero of~$y$,  or on $[t_0, t_1]$ we have
 $f=g$ and~$y=cz$ for some constant $c\in\R^\times$. 
\end{theorem}

\noindent
Here is an immediate consequence:

\begin{cor} \label{cor:gen osc closed upward}
If $g\in\c_a$  generates oscillation and   $f(t)\geq g(t)$, eventually, then~$f$ also generates oscillation.
\end{cor}

\begin{example}
For $k\in\R^\times$ we have the differential equation of the harmonic oscillator,
$$y'' + k^2 y\ =\ 0.$$
A function $y\in \Cat$ is a solution iff for some real constants $c,t_0$ and all $t\ge a$,
$$y(t)\ =\  c\sin k(t-t_0).$$
For $c\ne 0$, any function $y\in \Cat$ as displayed oscillates.  
Thus
if $f(t)\geq \varepsilon$, eventually, for some constant~$\varepsilon>0$, then $f$ generates oscillation. 
\end{example}

\noindent
To \eqref{eq:2nd order} we associate the corresponding {\bf Riccati equation}
\begin{equation}\label{eq:Riccati}\tag{R}
z'+z^2+f\ =\ 0.
\end{equation}
Let $y\in\Sol(f)^{\neq}$ be a non-oscillating solution to \eqref{eq:2nd order}, 
and take $b\geq a$ with~${y(t)\neq 0}$ for $t\geq b$. Then the function
$$t\mapsto z(t)\ :=\ y'(t)/y(t)\ \colon\ [b,+\infty)\to\R$$ 
in $\mathcal C^1_b$ satisfies \eqref{eq:Riccati}.
Conversely, if $z\in \mathcal C^1_b$ ($b\geq a$) is a solution to \eqref{eq:Riccati}, then 
$$t\mapsto y(t)\ :=\ \exp \left(\int_b^t z(s) \,ds\right)\ \colon\ [b,+\infty)\to\R$$
is a non-oscillating solution to \eqref{eq:2nd order} with $y\in(\c^1_b)^\times$ and $z=y^\dagger$.

\medskip
\noindent
Let $g\in \Cao$, $h\in\Caz$ and consider the second-order linear differential equation 
\begin{equation}\label{eq:2nd order, gen}\tag{$\tilde{\operatorname{L}}$}
y''+gy'+hy\ =\ 0.
\end{equation}
For the next corollary, see also \cite[Chapter~6, Lemma~4]{Bellman}.

\begin{cor}\label{coroscgen}
Set $f:=-\frac{1}{2}g'-\frac{1}{4}g^2+h\in \c_a$. Then the following are equivalent:
\begin{enumerate}
\item[\textup{(i)}] some nonzero solution of \eqref{eq:2nd order, gen}  oscillates;
\item[\textup{(ii)}] all nonzero solutions of \eqref{eq:2nd order, gen} oscillate; 
\item[\textup{(iii)}] $f$ generates oscillation.
\end{enumerate}
\end{cor}
\begin{proof}
Let $G\in (\Cat)^\times$ be given by  
$G(t):=\exp\!\left(-\frac{1}{2}\int_a^t g(s)\,ds\right)$. 
Then $y\in\Cat$ is a solution to \eqref{eq:2nd order} iff $Gy$ is a solution to \eqref{eq:2nd order, gen}; cf. [ADH, 5.1.13].
\end{proof}

\subsection*{Non-oscillation}  We continue with \eqref{eq:2nd order}. 
Let~$y_1$,~$y_2$ range over elements of $\Sol(f)$, and recall that its  Wronskian~$w=y_1y_2'-y_1'y_2$ is a real constant.

\begin{lemma}\label{y1y2}
Suppose $b\geq a$ is such that $y_2(t) \neq 0$ for $t\geq b$. Then for~$q:=y_1/y_2\in\c^2_b$ we have
$q'(t)=-w/y_2(t)^2$  for $t\geq b$,
so $q$ is monotone and $\lim_{t\to\infty} q(t)$ exists in $\R\cup\{-\infty,+\infty\}$.
\end{lemma}

\noindent
This leads to:

\begin{cor}\label{cor:I1 I2}
Suppose  $b\geq a$ and $y_1(t),y_2(t)\neq 0$ for $t\geq b$. For $i=1,2$, set
$$h_i(t)\  :=\  \int_b^t \frac{1}{y_i(s)^2}  \, ds\quad \text{ for $t\ge b$, so $h_i\in \mathcal{C}_b^3$}.$$
Then:  $y_1\prec y_2\ \Longleftrightarrow\ h_1 \succ 1 \succeq h_2$.
\end{cor}
\begin{proof} 
Let $q:=y_1/y_2$; so $q'=-wh_2'$ by  Lemma~\ref{y1y2}, hence $q+wh_2$ is constant. Thus, if $w\ne 0$, then
$y_1\preceq y_2 \Leftrightarrow q\preceq 1\Leftrightarrow h_2\preceq 1$.

Suppose $y_1\prec y_2$. Then $y_1$, $y_2$ are $\R$-linearly independent, so $w\ne 0$,  and~$h_2\preceq 1$. 
Note that $h_1$ is strictly increasing. Also $h_1y_1\in \Sol(f)$ by a routine computation. If~$h_1(t)\to r\in \R$ as~$t\to +\infty$,  then
$z:=(r-h_1)y_1\in \Sol(f)$, and~$z\prec y_1$, so~$z=0$, hence $h_1=r$, a contradiction. Thus~$h_1\succ 1$. 

For the converse, suppose $h_1\succ 1\succeq h_2$. Then $y_1$, $y_2$ are $\R$-linearly independent, so $w\ne 0$, and $q\preceq 1$. If $q(t)\to r\in \R^{\ne}$ as $t\to +\infty$, then
$y_1=qy_2\asymp y_2$, and thus~$h_1\asymp h_2$, a contradiction. Hence $q\prec 1$, and thus $y_1\prec y_2$. 
\end{proof}

\noindent
The pair $(y_1,y_2)$ is said to be a {\bf principal system} of solutions of \eqref{eq:2nd order} if 
{\samepage 
\begin{enumerate}
\item $y_1(t),y_2(t)>0$ eventually, and
\item $y_1 \prec y_2$.
\end{enumerate}}
Then $y_1$, $y_2$ form a basis of the $\R$-linear space $\Sol(f)$, and $f$ does not generate oscillation, by Lemma~\ref{lem:some vs all oscillate}. Moreover, for  
$y=c_1 y_1+c_2 y_2$ with $c_1,c_2\in\R$, $c_2\neq 0$ we have $y\sim c_2 y_2$.
Here are some facts about this notion:

\begin{lemma}\label{lem:admissible pair, unique}
If $(y_1,y_2)$, $(z_1,z_2)$ are principal systems of solutions of \eqref{eq:2nd order}, then
there are $c_1,d_1,d_2\in\R$ such that
$z_1 = c_1 y_1$, $z_2 = d_1 y_1+d_2 y_2$, and $c_1,d_2>0$.
\end{lemma}

\begin{lemma}\label{lem:admissible pair}
Suppose $f$ does not generate oscillation. Then \eqref{eq:2nd order} has a principal system of solutions.
\end{lemma}
\begin{proof}
It suffices to find a basis $y_1$, $y_2$ of $\Sol(f)$ with $y_1\prec y_2$.  
Suppose $y_1$, $y_2$ is any basis of $\Sol(f)$, and set   $c:=\lim_{t\to\infty} y_1(t)/y_2(t)\in\R\cup\{-\infty,+\infty\}$.
If $c=\pm\infty$, then    interchange $y_1$, $y_2$, otherwise   replace $y_1$ by $y_1-cy_2$. Then $c=0$, so $y_1\prec y_2$.
\end{proof}

\noindent
One calls $y_1$  a {\bf principal} solution of \eqref{eq:2nd order} if $(y_1,y_2)$ is a principal system of solutions of \eqref{eq:2nd order}  for some $y_2$.\index{solution!principal}\index{principal!solution} (See \cite[Theorem~XI.6.4]{Hartman} and \cite{Leighton,LeightonMorse}.) By the previous two lemmas, \eqref{eq:2nd order} has a principal solution iff $f$ does not generate oscillation, and any two principal solutions differ by a multiplicative factor in $\R^>$.
If $y_1\in (\c_a)^\times$ and $y_2$ is given as in Lemma~\ref{lem:2nd lin indep sol}, then~$y_2$  is a non-principal solution of  \eqref{eq:2nd order} and~$y_1\notin \R y_2$.  

\begin{remarkNumbered}[Hardy-type inequality associated to \eqref{eq:2nd order}]\label{rem:Hardy inequ, 1}
Suppose $f(t)>0$ for all~${t>a}$ and \eqref{eq:2nd order} has a solution~$y$ such that  $y(t),y'(t)>0$ for all $t>a$. Then for some~$C=C_f\in\R^{\geq}$, every  $u\in\c^1_a$ with~$u(a)=0$ satisfies
\begin{equation}\label{eq:Hardy inequ}
\int_a^\infty \abs{u(t)}^2\,f(t)\, dt \ \leq\ C\int_a^\infty \abs{u'(t)}^2\,dt.
\end{equation}
For $a>0$ and $f(t):=\frac{1}{4t^2}$ for $t\geq a$ this was shown by Hardy~\cite{Hardy20,Hardy25}; here one can take $C=1$, and this is optimal~\cite{Landau}. For the general case, see \cite[Theorem~4.1]{OpicKufner}.
\end{remarkNumbered}

\section{Hardy Fields}\label{sec:Hardy fields}

\noindent
In this brief section we introduce Hardy fields and review some classical extension theorems for Hardy fields.

\subsection*{Hardy fields}
A {\em Hardy field\/} is a subfield of 
$\Gi$ that is closed under the derivation of $\Gi$; see also~[ADH, 9.1].\index{Hardy field} Let~$H$ be a Hardy field.
Then $H$ is considered as an ordered valued differential field in the obvious way; see Section~\ref{sec:germs} for the ordering and valuation on $H$.
The field of constants of $H$ is $\R\cap H$. Hardy fields are pre-$H$-fields, and $H$-fields if they contain~$\R$;
see [ADH, 9.1.9(i),~(iii)].
As in Section~\ref{sec:germs} we equip the differential subfield $H[\imag]$  of~$\Calinf[\imag]$ with the unique valuation ring 
whose intersection with $H$ is the valuation ring of~$H$. Then $H[\imag]$ is a pre-$\d$-valued field of $H$-type with small derivation and constant field $\C\cap H[\imag]$;
if $H\supseteq\R$, then $H[\imag]$ is $\d$-valued   with constant field~$\C$. 
(Section~\ref{sec:upper lower bds} has an example of a differential subfield of $\Gi[\imag]$ that is not contained in $H[\imag]$ for any Hardy field $H$.)
 
We also consider variants:  a {\em $\Ginf$-Hardy field\/} is a Hardy field $H\subseteq \Ginf$, and a {\em $\Gom$-Hardy field\/} (also called an {\em analytic Hardy field\/}) is a Hardy field~${H\subseteq \Gom}$.\index{Cr-Hardy field@$\c^r$-Hardy field}\index{Hardy field!analytic}\index{Hardy field!smooth} 
Most Hardy fields arising in practice are actually $\Gom$-Hardy fields.

\subsection*{Hardian germs} Let $y\in\mathcal G$. 
Following~\cite{Sj}  we call $y$  {\bf hardian} if it lies in a Hardy field (and thus $y\in\Calinf$).
We also say that $y$ is {\bf $\Ginf$-hardian} if $y$ lies in a $\Ginf$-Hardy field,
equivalently, $y\in\Ginf$ and $y$ is hardian; 
likewise with $\Gom$ in place of~$\Ginf$.\index{germ!hardian}\index{hardian}
Let $H$ be a Hardy field. Call $y\in\mathcal G$ {\bf $H$-har\-dian} (or {\bf hardian over $H$}) if $y$ lies in a Hardy field extension
of~$H$. (Thus $y$ is hardian iff~$y$ is $\Q$-hardian.) 
If $H$ is a $\Ginf$-Hardy field and $y\in\Ginf$ is hardian over $H$, then $y$ generates a $\Ginf$-Hardy field extension~$H\langle y\rangle$ of $H$; likewise with $\Gom$ in place of $\Ginf$.

\subsection*{Maximal and perfect Hardy fields}
Let $H$ be a Hardy field. Call~$H$  {\it maximal}\/ if no Hardy field properly contains $H$.\index{Hardy field!maximal}\label{p:E(H)}
Following Boshernitzan~\cite{Boshernitzan82} we denote by $\Ex(H)$ the intersection of all maximal Hardy fields containing $H$;
thus~$\Ex(H)$ is a Hardy field extension of $H$, and a maximal Hardy field contains $H$ iff it contains~$\Ex(H)$, so $\Ex(\Ex(H))=\Ex(H)$. 
If $H^*$ is a Hardy field extension of~$H$, then $\Ex(H)\subseteq\Ex(H^*)$; hence if
$H^*$ is a Hardy field with $H\subseteq H^*\subseteq\Ex(H)$, then~$\Ex(H^*)=\Ex(H)$. 
Note 
that~$\Ex(H)$ consists of the $f\in \mathcal G$ that are hardian over each Hardy field  $E\supseteq H$. Hence $\Ex(\Q)$ consists of the germs in $\mathcal G$  that are  hardian over each Hardy field. 
As in \cite{Boshernitzan82}  we also say that~$H$ is {\bf perfect} if~$\Ex(H)=H$.\index{Hardy field!perfect} (This terminology is slightly unfortunate, since Hardy fields,
being of characteristic zero, are perfect as  fields.)
Thus $\Ex(H)$ is the smallest perfect Hardy field extension of $H$. Maximal Hardy fields are perfect.

\subsection*{Differentially maximal Hardy fields} 
Let $H$ be a Hardy field.
We now define differentially-algebraic variants of the above:
call $H$   {\bf differentially maximal}, or {\bf $\d$-maximal}\/ for short, if $H$ has no proper $\d$-algebraic Hardy field extension.\index{Hardy field!d-maximal@$\d$-maximal}\index{Hardy field!differentially maximal} 
  Every maximal Hardy field is $\d$-maximal, so each Hardy field is contained in a $\d$-maximal one; in fact,
by Zorn, each Hardy field~$H$ has a $\d$-maximal Hardy field extension which is $\d$-algebraic over $H$. Let~$\Dx(H)$ be the intersection of all 
$\d$-maximal Hardy fields containing $H$. Then~$\Dx(H)$ is a $\d$-algebraic Hardy field extension of~$H$ with~$\Dx(H)\subseteq\Ex(H)$. We have  $\Dx(H)=\Ex(H)$ iff $\Ex(H)$ is $\d$-algebraic over~$H$:\label{p:D(H)}

\begin{lemma}\label{lem:Dx Ex}
$\Dx(H)=\big\{f\in\Ex(H):\text{$f$ is $\d$-algebraic over $H$}\big\}$.
\end{lemma}
\begin{proof}
We only need to show the inclusion ``$\supseteq$''. For this let~$f\in\Ex(H)$ be $\d$-algebraic over $H$, and let $E$ be a $\d$-maximal Hardy field extension of $H$; we need to show $f\in E$. To see this extend $E$
 to a maximal Hardy field $M$; then $f\in M$, hence $f$ generates a Hardy field extension $E\langle f\rangle$ of $E$.
 Since $f$ is $\d$-algebraic over~$H$ and thus over~$E$, this yields $f\in E$ by $\d$-maximality of $E$, as required.
\end{proof}

\noindent
A $\d$-maximal Hardy field contains $H$ iff it contains~$\Dx(H)$, hence $\Dx(\Dx(H))=\Dx(H)$.
For any Hardy field $H^*\supseteq H$  we have $\Dx(H^*)\supseteq \Dx(H)$,   hence if also $H^*\subseteq\Dx(H)$, then $\Dx(H^*)=\Dx(H)$.
We say that $H$ is {\bf $\d$-per\-fect}\index{Hardy field!d-perfect@$\d$-perfect} if $\Dx(H)=H$.  
Thus~$\Dx(H)$ is the smallest $\d$-perfect Hardy field extension of~$H$. Every perfect Hardy field is $\d$-perfect, as is every $\d$-maximal Hardy field. The following diagram summarizes the various implications among these properties of Hardy fields:
$$\xymatrix{ \text{maximal} \ar@{=>}[r]\ar@{=>}[d]  &  \text{perfect} \ar@{=>}[d] \\ 
\text{$\d$-maximal} \ar@{=>}[r]  &  \text{$\d$-perfect}  }$$
We call $\Dx(H)$ the {\bf $\d$-perfect hull} of $H$, and $\Ex(H)$ the {\bf perfect hull} of $H$.\index{Hardy field!d-perfect hull@$\d$-perfect hull}\index{Hardy field!perfect hull}

\subsection*{Variants of the perfect hull}   
Let  $H$ be a $\c^r$-Hardy field where   $r\in\{\infty,\omega\}$. We say that  $H$ is {\bf $\c^r$-maximal} if no $\c^r$-Hardy field properly contains it. By Zorn, $H$ has a $\c^r$-maximal extension. In analogy with $\Ex(H)$,  define the {\bf $\c^r$-perfect hull~$\Ex^r(H)$} of $H$ to be the
intersection of all $\c^r$-maximal Hardy fields containing~$H$. \label{p:Er(H)} We say that  $H$ is
  {\bf $\c^r$-perfect} if $\Ex^r(H)=H$. The penultimate subsection goes through
  with {\em Hardy field}, {\em maximal}, {\em hardian}, {\em $\Ex(\,\cdot\,)$}, and {\em perfect} replaced by {\em $\c^r$-Hardy field},
  {\em $\c^r$-maximal},  {\em $\c^r$-hardian}, {\em $\Ex^r(\,\cdot\,)$}, and {\em $\c^r$-perfect}, respectively.    (In \cite{ADH5} we show that no analogue of $\Dx(H)$ is needed
  for the $\c^r$-category.)\index{Cr-Hardy field@$\c^r$-Hardy field!Cr-maximal@$\c^r$-maximal}\index{Cr-Hardy field@$\c^r$-Hardy field!Cr-perfect@$\c^r$-perfect}\index{Cr-Hardy field@$\c^r$-Hardy field!Cr-perfect@$\c^r$-perfect hull}

\subsection*{Some basic extension theorems} 
We summarize some well-known extension results for Hardy fields, cf.~\cite{Bou, Sj, Ros}:

\begin{prop}\label{prop:Hardy field exts} Any Hardy field $H$ has the following Hardy field extensions: \begin{enumerate}
\item[\textup{(i)}] $H(\R)$, the subfield of $\Gi$ generated by $H$ and $\R$;
\item[\textup{(ii)}] $H^{\operatorname{rc}}$, the real closure of $H$ as defined in Proposition~\ref{b1};
\item[\textup{(iii)}] $H(\ex^f)$ for any $f\in H$;
\item[\textup{(iv)}] $H(f)$ for any $f\in \Go$ with $f'\in H$;
\item[\textup{(v)}] $H(\log f)$ for any $f\in H^{>}$.
\end{enumerate}
If $H$ is contained in $\Ginf$, then
so are the Hardy fields in {\rm (i), (ii), (iii), (iv), (v)}; like\-wise with $\Gom$ instead of $\Ginf$. 
\end{prop}

\noindent
Note that (v) is a special case of (iv), since $(\log f)'=f^\dagger\in H$
for $f\in H^{>}$. Another special case of (iv) is that a Hardy field $H$ yields a Hardy field $H(x)$. It also yields the Hardy field $H_{\operatorname{LE}}$ from the introduction as the smallest (under inclusion) Hardy field extension of $\R(x)$ that is log-closed and exp-closed. 

 A consequence of the proposition is that any Hardy field $H$ has a smallest real closed Hardy field extension~$L$ with~$\R\subseteq L$ such that for all $f\in L$ we have $\ex^f\in L$
and $g'=f$ for some~$g\in L$. Note that then~$L$ is a Liouville closed $H$-field as defined in
Section~\ref{sec:prelim}. 
Let $H$ be a Hardy field with~$H\supseteq\R$.  As in~\cite{AvdD2} and [ADH, p.~460] we  then  denote the above
$L$ by~$\Li(H)$; so~$\Li(H)$ is the smallest   Liouville closed Hardy field containing $H$,
called the {\it Hardy-Liouville closure}\/ of $H$ in \cite{ADH2}. 
We have  $\Li(H)\subseteq \Dx(H)$, hence if $H$ is $\d$-perfect, then~$H$ is a Liouville closed $H$-field. Moreover, if $H\subseteq\Ginf$ then $\Li(H)\subseteq\Ginf$, and
similarly with $\Gom$ in place of~$\Ginf$.\index{Hardy field!Hardy-Liouville closure}\label{p:HL}

\medskip
\noindent
The next more general result in Rosenlicht~\cite{Ros} is attributed there to M. Singer: 

{\sloppy
\begin{prop}\label{singer} Let $H$ be a Hardy field and $p,q\in H[Y]$, $y\in \mathcal{C}^1$, such that~$y'q(y)=p(y)$ with $q(y)\in\mathcal{C}^\times$. Then~$y$ generates a Hardy field~$H(y)$ over~$H$.
\end{prop} }

\noindent
Note that for $H$, $p$, $q$, $y$ as in the proposition we have~$y\in\Dx(H)$.

\subsection*{Compositional  conjugation in Hardy fields} Let now $H$ be a Hardy field, and
let~$\ell\in \Go$ be such that $\ell>\R$ and $\ell'\in H$. Then $\ell\in \Calinf$,  $\phi:=\ell'$
is active in $H$, $\phi>0$, and we have a Hardy field $H(\ell)$. 
 The
$\C$-algebra automorphism~$f\mapsto f^\circ:= f\circ \ell^{\inv}$ of~$\mathcal{C}[\imag]$  restricts to an ordered field isomorphism
$$h\mapsto h^\circ\ :\ H \to H^\circ:=H\circ \ell^{\inv}.$$
The identity  
$(f^\circ)'  = (\phi^{-1}f')^\circ$, valid for each $f\in\Go[\imag]$, shows
that  $H^\circ$ is again a Hardy field. Conversely, if $E$ is a subfield of $\Calinf$ with $\phi\in E$ and $E^\circ:= E\circ \ell^{\inv}$ is a Hardy field, then $E$ is a Hardy field. If $H\subseteq \Ginf$ and $\ell\in \Ginf$, then
$H^\circ\subseteq \Ginf$; likewise with $\Gom$ instead of~$\Ginf$. 
If $E$ is a Hardy field extension of $H$, then
 $E^\circ$ is a Hardy field extension of $H^\circ$, and
$E$ is $\d$-algebraic over~$H$ iff $E^\circ$ is $\d$-algebraic over $H^\circ$.
Hence $H$ is maximal iff $H^\circ$ is maximal, and likewise with ``$\d$-maximal'' in place of ``maximal''.
So~$\Ex(H^\circ)=\Ex(H)^\circ$ and~$\Dx(H^\circ)=\Dx(H)^\circ$, and thus $H$ is perfect iff $H^\circ$ is perfect, and likewise with ``$\d$-perfect'' in place of ``perfect''. 
The next lemma is
 \cite[Corollary~6.5]{Boshernitzan81}; see also \cite[Theo\-rem~1.7]{AvdD4}.
 
\begin{lemma}\label{lem:Bosh6.5}
The germ $\ell^{\inv}$ is hardian. Moreover,
if $\ell$ is $\Ginf$-hardian, then $\ell^{\inv}$ is also $\Ginf$-hardian, and
likewise with  $\Gom$  in place of $\Ginf$. 
\end{lemma}
\begin{proof}
By Proposition~\ref{prop:Hardy field exts}(iv) we can arrange that our Hardy field $H$ contains both~$\ell$ and~$x$.
Then $\ell^{\inv}=x\circ \ell^{\inv}$ is an element of the Hardy field $H\circ\ell^{\operatorname{inv}}$.
\end{proof}

\noindent
Next we consider the pre-$\d$-valued field $K:=H[\imag]$ of $H$-type, which gives rise to
$$K^\circ:=K\circ \ell^{\inv}=H^\circ[\imag],$$ 
also a pre-$\d$-valued field of $H$-type, and we have the valued field isomorphism
\[h\mapsto h^\circ\ :\ K\to K^\circ.\]
Note: $h\mapsto h^\circ\colon H^{\phi} \to H^\circ$ is an isomorphism 
of pre-$H$-fields, and~${h\mapsto h^\circ\colon K^{\phi} \to K^\circ}$ is an isomorphism of valued differential fields.
Recall that $K$ and $K^\phi$ have the same underlying field.

\begin{lemma}\label{lemlioucomp} From the isomorphisms $H^\phi\cong H^\circ$ and $K^\phi\cong K^\circ$ we obtain:
If $H$ is Liouville closed, then so is $H^\circ$. If $\I(K)\subseteq K^\dagger$, then $\I(K^\circ)\subseteq (K^\circ)^\dagger$.
\end{lemma}

\section{Upper and Lower Bounds on the Growth of Hardian Germs}\label{sec:upper lower bds}

\noindent
In this section   we use logarithmic decompositions to simplify arguments in \cite{Boshernitzan82,Boshernitzan86,Rosenlicht83}.    
It is not used for proving our main theorem, but some of it is needed for its applications, in the proofs of Corollary~\ref{cor:17.11},    Pro\-po\-si\-tion~\ref{prop:translog}, and Theorem~\ref{thm:coinitial in E(H)}. 

{\sloppy
\subsection*{Generalizing logarithmic decomposition}
In this subsection $K$ is a differential ring and~${y\in K}$. 
In [ADH, p.~213] we defined the $n$th iterated logarithmic derivative of $y^{\langle n\rangle}$ when $K$ is a differential field. (See also Section~\ref{sec:prelim}.)
Generalizing this, set~$y^{\langle 0\rangle}:=y$, and recursively, if~$y^{\langle n\rangle}\in K$ is defined and a unit in $K$, then
 $y^{\langle n+1\rangle}:=(y^{\langle n\rangle})^\dagger$, while otherwise~$y^{\langle n+1\rangle}$ is not defined. (Thus if~$y^{\langle n\rangle}$ is defined, then so are~$y^{\langle 0\rangle},\dots,y^{\langle n-1\rangle}$.)
  With~$L_n$ in~$\Z[X_1,\dots,X_n]$ as in~[ADH, p.~213], if $y^{\langle n\rangle}$ is defined, then
$$y^{(n)}\  =\  y^{\langle 0\rangle} \cdot L_n(y^{\langle 1\rangle},\dots, y^{\langle n\rangle}).$$
 If $y^{\langle n\rangle}$ is defined and $\i=(i_0,\dots,i_n)\in\N^{1+n}$, we set
 $$y^{\langle \i \rangle}\ :=\  (y^{\langle 0\rangle})^{i_0} (y^{\langle 1\rangle})^{i_1} \cdots (y^{\langle n\rangle})^{i_n}\in K.$$
Hence if $H$ is a differential subfield of $K$, $P\in H\{Y\}$ has order at most~$n$ and
logarithmic decomposition $P=\sum_{\i} P_{\langle \i\rangle} Y^{\langle \i\rangle}$ 
($\i$ ranging over $\N^{1+n}$, all~$P_{\langle \i\rangle}\in H$, and~$P_{\langle \i\rangle}=0$ for all but finitely many $\i$),
and $y^{\langle n\rangle}$ is defined, then~$P(y)=\sum_{\i}  P_{\langle \i\rangle} y^{\langle \i\rangle}$.
Below we apply these remarks to~$K=\Calinf$, where for $y\in K^\times$ we have $y^\dagger=(\log |y|)'$, hence
$y^{\langle n+1\rangle}=(\log |y^{\langle n\rangle}|)'$ if $y^{\langle n+1\rangle}$ is defined. 
}

\subsection*{Transexponential germs}
For  $f\in\c$ we recursively define the germs~$\exp_n f$ in $\c$ by $\exp_0 f:=f$ and $\exp_{n+1} f:=\exp (\exp_n f)$.
Following \cite{Boshernitzan82} we say that a germ~$y\in\c$ is {\bf transexponential} if $y\geq\exp_n x$ for all~$n$.\index{germ!transexponential}\index{transexponential}
{\it In the rest of this subsection $H$ is a Hardy field.}\/
By Corollary~\ref{cor:val at infty} and Proposition~\ref{prop:Hardy field exts}:

\begin{lemma}\label{lem:bounded Hardy field ext, 1}
If the $H$-hardian germ $y$ is $\d$-algebraic over $H$, then $y\leq\exp_n h$ for some $n$ and some $h\in H(x)$.
\end{lemma}

\noindent
Thus each transexponential hardian germ is $\d$-transcendental (over $\R$).  
{\it In the rest of this subsection:  $y\in\Calinf$ is transexponential and hardian, and $z\in\Calinf[\imag]$.}\/
Then~$y^{\langle n\rangle}$ is defined, and~$y^{\langle n\rangle}$ is also transexponential and hardian,  for all  $n$.
Next some variants of results from Section~\ref{sec:prelim}. 
For this, let $n$ be given and let $f\in\Calinf$, not necessarily hardian, be such that~$f\succ 1$, $f\geq 0$, and~$y\succeq\exp_{n+1}f$.

\begin{lemma}\label{lem:expn variant}
We have   $y^\dagger\succeq\exp_n f$ and $y^{\langle n\rangle}\succeq\exp f$.
\end{lemma}
\begin{proof} Since $y\succeq \exp_2 x$, we have
 $\log y \succeq \exp x$ by Lemma~\ref{lem:log preceq},  and thus 
$y^\dagger=(\log y)'\succeq \log y$.  Since $y\succeq\exp_{n+1}f$, the same lemma gives $\log y \succeq \exp_n f$. Thus~$y^\dagger\succeq\exp_n f$.
Now the second statement follows by an easy induction. 
\end{proof}

\begin{cor}\label{cor:expn variant}
Let $\i\in\Z^{1+n}$  and suppose $\i>0$ lexicographically. Then $y^{\langle\i\rangle} \succ f$.
\end{cor}
\begin{proof}
Let~$m\in\{0,\dots,n\}$ be minimal such that $i_m\neq 0$; so $i_m\geq 1$.
The remarks after Corollary~\ref{cor:itpsi} then give $y^{\langle\i\rangle} \succ 1$ and
$[v(y^{\langle\i\rangle})]=[v(y^{\langle m\rangle})]$, so we have~$k\in\N$, $k\geq 1$, such that $y^{\langle\i\rangle} \succeq (y^{\langle m\rangle})^{1/k}$.
Then Lemma~\ref{lem:expn variant} gives
$y^{\langle\i\rangle} \succeq (y^{\langle m\rangle})^{1/k}\succeq (\exp f)^{1/k}\succ f$ as required.
\end{proof}

\noindent
In the next proposition and lemma~$P\in H\{Y\}^{\neq}$ has order at most~$n$, and~$\i$,~$\j$,~$\k$ range over~$\N^{1+n}$.
Let  $\j$ be lexicographically maximal such that~$P_{\<\j\>}\neq 0$, and 
choose~$\k$ so that
$P_{\<\k\>}$ has   minimal valuation.
If $P_{\<\k\>}/P_{\<\j\>} \succ x$,  set
$f:=|P_{\<\k\>}/P_{\<\j\>}|$; otherwise set $f:=x$.
 Then~$f\in H(x)$, $f>0$, $f\succ 1$, and $f\succeq P_{\<\i\>}/P_{\<\j\>}$ for all~$\i$.
 
\begin{prop}\label{prop:val at infty variant}
We have $P(y)\sim P_{\<\j\>}y^{\langle \j\rangle}$ and thus 
$$P(y)\in \big(\Calinf\big)^\times, \qquad \sgn P(y)\ =\ \sgn P_{\<\j\>}\neq 0.$$
\end{prop}
\begin{proof}
For $\i<\j$ we have  $y^{\langle \j-\i \rangle}\succ f\succeq P_{\<\i\>}/P_{\<\j\>}$ by Corollary~\ref{cor:expn variant},
hence~$P_{\<\j\>}y^{\langle \j\rangle} \succ P_{\<\i\>}y^{\langle \i\rangle}$.
Thus $P(y)\sim P_{\<\j\>}y^{\langle \j\rangle}$.
\end{proof}

{\sloppy  
\begin{lemma}\label{lem:transexp}
Suppose   that  $z^{\langle n\rangle}$  is defined and $y^{\langle i\rangle} \sim z^{\langle i\rangle}$ for $i=0,\dots,n$.
Then~${P(y) \sim P(z)}$.
\end{lemma}}
\begin{proof}
For all $\i$ with $P_{\langle \i\rangle}\neq 0$ we have $P_{\langle \i\rangle}y^{\langle\i\rangle} \sim P_{\langle \i\rangle}z^{\langle\i\rangle}$, by Lem\-ma~\ref{lem:sim props}.
Now use that for   $\i\neq\j$ we have
$P_{\langle \i\rangle} y^{\langle\i\rangle} \prec P_{\langle \j\rangle} y^{\langle\j\rangle}$ by the proof of Proposition~\ref{prop:val at infty variant}.
\end{proof}

\noindent
From here on $n$ is no longer fixed. 

\begin{cor}[{Boshernitzan~\cite[Theorem~12.23]{Boshernitzan82}}] \label{cor:Bosh 12.23}
If   $y\geq \exp_n h$ for all~$h\in H(x)$  and all $n$, then $y$ is $H$-hardian.
\end{cor}

\noindent
This is an immediate consequence of Proposition~\ref{prop:val at infty variant}. (In \cite{Boshernitzan82}, the proof of this fact  is only indicated.)
From Lemma~\ref{lem:transexp} we also obtain: 

\begin{cor}\label{cor:transexp}
Suppose that $y$ is as in Corollary~\ref{cor:Bosh 12.23} and $z\in\Calinf$, and 
$z^{\langle n\rangle}$ is defined and $y^{\langle n\rangle} \sim z^{\langle n\rangle}$, for all $n$. Then
$z$ is $H$-hardian, and there is a unique ordered differential field isomorphism $H\langle y\rangle \to H\langle z\rangle$ over $H$ which sends~$y$ to $z$.
\end{cor}

\noindent
Lemma~\ref{lem:y H-hardian crit} below contains another criterion for $z$ to be $H$-hardian. This involves a certain binary relation $\sim_\infty$ on germs defined in the next subsection.
Lemma~\ref{lem:transexp} also yields a complex version of Corollary~\ref{cor:transexp}:

\begin{cor}\label{cor:transexp, complex}
Suppose that $y$ is as in Corollary~\ref{cor:Bosh 12.23} and that
$z^{\langle n\rangle}$ is defined and~$y^{\langle n\rangle} \sim z^{\langle n\rangle}$, for all $n$. Then
$z$ generates a differential subfield $H\langle z\rangle$ of $\Calinf[\imag]$, 
and there is a unique  differential field isomorphism $H\langle y\rangle \to H\langle z\rangle$ over $H$ which sends~$y$ to $z$.
Moreover, the binary relation $\preceq$ on $\c[\imag]$ restricts to a dominance relation on $H\langle z\rangle$
which makes this an isomorphism of valued differential fields.
\end{cor}

\subsection*{A useful equivalence relation}
We set
$$\Calinf[\imag]^{\preceq}\ :=\ \big\{ f\in \Calinf[\imag]: f^{(n)}\preceq 1\text{ for all $n$} \big\}\ \subseteq\ \c[\imag]^{\preceq},$$
a differential $\C$-subalgebra of  $\Calinf[\imag]$, and
 $$\mathcal I\ :=\ \big\{ f\in \Calinf[\imag]: f^{(n)}\prec 1\text{ for all $n$} \big\}\ \subseteq\ \Calinf[\imag]^{\preceq},$$
a differential ideal of $\Calinf[\imag]^{\preceq}$
(thanks to the Product Rule).
Recall from the remarks preceding   Lemma~\ref{lem:sim props} that~$(\c[\imag]^{\preceq})^\times = \c[\imag]^{\asymp}$. 

\begin{lemma}\label{lem:units of cinf[i]preceq}
The group of units of $\Calinf[\imag]^{\preceq}$ is
$$\Calinf[\imag]^{\asymp}\ :=\ \Calinf[\imag]^{\preceq}\cap\c[\imag]^{\asymp}\ =\ \big\{ f\in \Calinf[\imag]: f \asymp 1,\ f^{(n)}\preceq 1\text{ for all $n$} \big\}.$$
Moreover, $1+\mathcal I$ is a subgroup of~$\Calinf[\imag]^\asymp$.
\end{lemma}
\begin{proof}
It is clear that
$$(\Calinf[\imag]^{\preceq})^\times\ \subseteq\ \Calinf[\imag]^{\preceq}\cap  (\c[\imag]^{\preceq})^\times\ =\  \Calinf[\imag]^{\preceq}\cap\c[\imag]^{\asymp}\ =\ \Calinf[\imag]^{\asymp}.$$
Conversely, suppose $f\in \Calinf[\imag]$ satisfies $f\asymp 1$ and $f^{(n)}\preceq 1$ for all $n$.
For each~$n$ we have $Q_n\in\Q\{X\}$ such that
$(1/f)^{(n)} =  Q_n(f)/f^{n+1}$, hence~$(1/f)^{(n)} \preceq 1$. Thus~$f\in (\Calinf[\imag]^{\preceq})^\times$.
This shows the first statement.  
Clearly~${1+\mathcal I}\subseteq \Calinf[\imag]^{\asymp}$, and  $1+\mathcal I$ is closed under multiplication.
If $\delta\in\mathcal I$, then~$1+\delta$ is a unit of $\Calinf[\imag]^{\preceq}$ and~$(1+\delta)^{-1}=1+\varepsilon$ where $\varepsilon=-\delta(1+\delta)^{-1}\in\mathcal I$. 
\end{proof}

\noindent
For $y,z\in\c[\imag]^\times$ we  define
$$y\sim_\infty z \quad:\Longleftrightarrow\quad y\in z\cdot (1+\mathcal I);$$
hence $y\sim_\infty z\Rightarrow y\sim z$.
Lemma~\ref{lem:units of cinf[i]preceq} yields that $\sim_\infty$ is  an equivalence relation on~$\c[\imag]^\times$, and
for $y_i,z_i\in\c[\imag]^\times$ ($i=1,2$) we   have
$$y_1\sim_\infty y_2\quad \&\quad z_1\sim_\infty z_2 \qquad\Longrightarrow\qquad y_1z_1 \sim_\infty y_2z_2,\quad y_1^{-1}\sim_\infty y_2^{-1}.$$

\begin{lemma}\label{lem:y' siminfty z'}
Let $y,z\in\c^1[\imag]^\times$ with $y\sim_\infty z$ and $z\in z'\,\Calinf[\imag]^{\preceq}$. Then 
$$y', z'\in \c[\imag]^\times, \qquad y'\sim_\infty z'.$$
\end{lemma}
\begin{proof}
Let $\delta\in\mathcal I$ and $f\in\Calinf[\imag]^{\preceq}$  with $y=z(1+\delta)$ and $z=z'f$.
Then $z'\in\c[\imag]^\times$ and
$y'=z'(1+\delta)+z\delta'=z'(1+\delta+f\delta')$ where $\delta+f\delta'\in\mathcal I$, so $y'\sim_\infty z'$.
\end{proof}

\noindent
If  $\ell\in\c^n[\imag]$ and $f\in\c^n$ with $f\geq 0$, $f\succ 1$, then $\ell\circ f\in\c^n[\imag]$.
In fact, for $n\geq 1$ and~$1\le k\le n$ we have
a differential polynomial~$Q^n_k\in\Q\{X'\}\subseteq\Q\{X\}$ of order~$\le n$,  isobaric of weight~$n$, and homogeneous of degree $k$, 
such that for all such $\ell$, $f$, 
$$(\ell\circ f)^{(n)}\  =\  (\ell^{(n)}\circ f)\,Q^n_n(f)+\cdots+(\ell'\circ f)\,Q^n_1(f).$$
For example,
$$Q^1_1=X',\quad Q^2_2=(X')^2,\ Q^2_1=X'',\quad Q^3_3=(X')^3,\ Q^3_2=3X'X'',\ Q^3_1=X'''.$$
The following lemma is only used in the proof of Theorem~\ref{thm:coinitial in E(H)} below.

\begin{lemma}\label{lem:difference in I} Let $f,g\in\Calinf$ be such that $f,g\geq 0$ and $f,g\succ 1$, and set~$r:={g-f}$.
Suppose  $P(f)\cdot Q(r)\prec 1$ for all $P,Q\in\Q\{Y\}$ with $Q(0)=0$,
and let~$\ell\in\Calinf[\imag]$ be such that $\ell'\in\mathcal I$. Then $\ell\circ g - \ell\circ f\in \mathcal I$.
\end{lemma}
\begin{proof}
Treating real and imaginary parts separately we arrange $\ell\in\Calinf$. Note that $r\prec 1$. 
Taylor expansion [ADH, 4.2] for $P\in \Q\{X\}$ of order~$\leq n$ gives
$$P(g)-P(f)\  =\  \sum_{|\i|\geq 1} \frac{1}{\i!} P^{(\i)}(f) \cdot r^{\i} \qquad (\i\in \N^{1+n}),$$
and thus $P(g)-P(f)\prec 1$ and $rP(g)\prec 1$. 
The Mean Value Theorem yields a germ~$r_n\in\mathcal G$ such that
$$\ell^{(n)}\circ g - \ell^{(n)}\circ f\ =\ \big(\ell^{(n+1)}\circ (f+r_n)\big)\cdot r\quad\text{and}\quad |r_n| \leq |r|.$$
Now $r_0\prec 1$, so  $\ell'\circ (f+r_0)\prec 1$, hence $\ell\circ g-\ell \circ f\prec 1$.
For $1\le k\le n$,
\begin{multline*}
(\ell^{(k)}\circ g)\,Q^n_k(g) - (\ell^{(k)}\circ f)\,Q^n_k(f) = \\
 (\ell^{(k)}\circ f)\, \big( Q^n_k(g)-Q^n_k(f) \big) +  \big(\ell^{(k+1)}\circ(f+r_k)\big)\cdot r    Q^n_k(g),
\end{multline*}
so $(\ell^{(k)}\circ g)\,Q^n_k(g) - (\ell^{(k)}\circ f)\,Q^n_k(f)\prec 1$, and thus $\big( \ell\circ g - \ell\circ f \big)^{(n)}\prec 1$.
\end{proof}

\noindent
We consider next the differential $\R$-subalgebra $$(\Calinf)^{\preceq}\ :=\ \Calinf[\imag]^{\preceq}\cap\Calinf \ \subseteq\ \c^{\preceq}$$ of $\Calinf$. In the rest of this subsection $H$ is a Hardy field and $y,z\in\Calinf$, $y,z\succ 1$. 
Note that $(\Calinf)^{\preceq}\cap H=\mathcal O$ is the valuation ring of $H$ and $\mathcal I\cap H=\smallo$ is the maximal ideal of $\mathcal O$. This yields:

\begin{lemma}\label{lem:y siminf z}
Suppose $y-z\in (\Calinf)^{\preceq}$ and $z$ is hardian. Then $y\sim_\infty z$.
\end{lemma}
\begin{proof} From $y=z+f$ with $f\in (\Calinf)^{\preceq}$ we obtain $y=z(1+fz^{-1})$. Now $z^{-1}\in \mathcal{I}$, so
$fz^{-1}\in \mathcal{I}$, and thus $y\sim_\infty z$. 
\end{proof}

\noindent
We now formulate a sufficient condition involving $\sim_\infty$ for $y$  to be  $H$-hardian.

\begin{lemma}\label{lem:y H-hardian crit}
Suppose  
$z$ is $H$-hardian  with $z\geq\exp_n h$ for all $h\in H(x)$ and all~$n$, and~$y\sim_\infty z$. Then  $y$ is  $H$-hardian, and  there is a unique ordered differential field isomorphism~$H\langle y\rangle\to H\langle z\rangle$ which is the identity on $H$ and sends  $y$ to $z$.
\end{lemma}
\begin{proof}
By Lemma~\ref{lem:bounded Hardy field ext, 1} we may replace $H$ by the Hardy subfield
 $\Li\!\big(H(\R)\big)$ of~$\Ex(H)$ to arrange that $H\supseteq\R$ is Liouville closed.
By Corollary~\ref{cor:transexp} (with the roles of~$y$,~$z$ reversed) it is enough to show that for each $n$,
$y^{\langle n\rangle}$ is defined, $y^{\langle n\rangle}\succ 1$, and~${y^{\langle n\rangle}\sim_\infty z^{\langle n\rangle}}$. This holds by hypothesis for $n=0$. By Lemma~\ref{lem:expn, 1}, $z>H$ gives~$z^\dagger>H$, so~$z=z'f$ with $f \prec 1$ in the Hardy field $H\langle z\rangle$, hence $f^{(n)}\prec 1$ for all $n$. So by Lemma~\ref{lem:y' siminfty z'}, $y^{\langle 1\rangle}=y^\dagger$ is defined, $y^{\langle 1\rangle}\in (\Calinf)^\times$, $y^{\langle 1\rangle} \sim_\infty z^{\langle 1\rangle}$, and thus $y^{\langle 1\rangle}\succ 1$. 
Assume for a certain $n\geq 1$ that $y^{\langle n\rangle}$ is defined,  $y^{\langle n\rangle}\succ 1$, and~$y^{\langle n\rangle}\sim_\infty z^{\langle n\rangle}$. Then~$z^{\langle n\rangle}$ is $H$-hardian and $H<z^{\langle n\rangle}$ by Lemma~\ref{lem:expn, 3}.
Hence  by the case $n=1$ applied to $y^{\langle n\rangle}$, $z^{\langle n\rangle}$ in place of~$y$,~$z$, respectively,  $y^{\langle n+1\rangle}=(y^{\langle n\rangle})^\dagger$ is defined, $y^{\langle n+1\rangle}\succ 1$, and~$y^{\langle n+1\rangle}\sim_\infty z^{\langle n+1\rangle}$.
\end{proof}

\noindent
The next two corollaries are Theorems~13.6 and 13.10, respectively, in  \cite{Boshernitzan82}.

\begin{cor}\label{cor:yhardian}
Suppose $z$ is trans\-ex\-po\-nen\-tial and  hardian, and~$y-z\in (\Calinf)^{\preceq}$.
Then  $y$  is hardian, and there is a unique  isomorphism~$\R\langle y\rangle\to\R\langle z\rangle$ of ordered differential fields that is the identity on $\R$ and sends  $y$ to $z$.
\end{cor}
\begin{proof}
Take $H:=\Li(\R)$. Then $z$ lies in a Hardy field extension of $H$, name\-ly~$\Li\!\big(\R\langle z\rangle\big)$, and $H<z$. So $y\sim_\infty z$ by Lemma~\ref{lem:y siminf z}. Now use Lemma~\ref{lem:y H-hardian crit}. \end{proof}

\begin{cor} \label{cor:Bosh13.10}
If~$z\in\Ex(H)^{>\R}$, then $z \leq \exp_n h$ for some $h\in H(x)$ and some~$n$. \textup{(}Thus if $x\in H$ and
$\exp H \subseteq H$, then $H^{>\R}$ is cofinal in $\Ex(H)^{>\R}$.\textup{)}
\end{cor}
\begin{proof}
Towards a contradiction, suppose $z\in\Ex(H)^{>\R}$ and $z>\exp_n h$ in $\Ex(H)$ for all $h\in H(x)$ and all $n$.   Set~$y:=z+\sin x$. Then 
 $y$ is $H$-hardian   by Lemmas~\ref{lem:y siminf z} and~\ref{lem:y H-hardian crit}, so $y$, $z$ lie in a common Hardy field
 extension of $H$, a contradiction. 
\end{proof}

\begin{remark}
The same proof shows that Corollary~\ref{cor:Bosh13.10} remains true if $H$ is  a $\Ginf$-Hardy field and
$\Ex(H)$ is replaced by~$\Ex^\infty(H)$; likewise for $\omega$ in place of~$\infty$. 
\end{remark}

\noindent
Next a lemma
similar to Lemma~\ref{lem:y H-hardian crit}, but obtained using Corollary~\ref{cor:transexp, complex} instead of Corollary~\ref{cor:transexp}:

\begin{lemma}\label{lem:y H-hardian crit, complex}
Let $H$ be a Hardy field, let  
$z\in\Calinf$ be $H$-hardian  with $z\geq\exp_n h$ for all $h\in H(x)$ and all~$n$, and~$y\in\Calinf[\imag]$ with $y\sim_\infty z$. Then  $y$ generates a differential subfield~$H\langle y\rangle$ of $\Calinf[\imag]$, and  there is a unique   differential field isomorphism~$H\langle y\rangle\to H\langle z\rangle$ which is the identity on $H$ and sends  $y$ to $z$. The binary relation $\preceq$ on $\c[\imag]$ restricts to a dominance relation on $H\langle y\rangle$
which makes this an isomorphism of valued differential fields.
\end{lemma}

\noindent
We use the above at the end of the next subsection to produce a differential subfield of $\Calinf[\imag]$
that is not contained in $H[\imag]$ for any Hardy field $H$.

\subsection*{Boundedness} 
Let $H\subseteq\c$.
We say that $b\in\c$ {\bf bounds $H$} if $h\leq b$ for each~$h\in H$. We call $H$ {\bf bounded}\index{Hardy field!bounded} if some $b\in\c$ bounds $H$,
and we    call $H$ {\bf unbounded}\index{Hardy field!unbounded} if $H$ is not bounded.
If $H_1, H_2\subseteq\c$ and for each $h_2\in H_2$ there is an $h_1\in H_1$ with~$h_2\leq h_1$, then 
any $b\in\c$ bounding $H_1$ also  bounds $H_2$.
Every bounded subset of $\c$ is bounded by a germ in~$\Gom$; this follows from \cite[Lemma~14.3]{Boshernitzan82}:

\begin{lemma} \label{lem:Bosh 14.3} 
For every $b\ge 0$ in $\c^\times$ there is a~$\phi\ge 0$ in $(\c^\omega)^\times$ such that $\phi^{(n)}\prec b$ for all $n$.
\end{lemma}

\noindent
Every countable subset of~$\c$ is bounded, by  
du Bois-Reymond~\cite{dBR75}; see also~\cite[Chap\-ter~II]{Ha} or~\cite[Cha\-pi\-tre~V, p.~53, ex.~8]{Bou}.
Thus $H\subseteq\c$   is bounded if it is totally ordered by the partial ordering $\le$ of $\c$ and has countable cofinality. 
If $H$ is a Hausdorff field and $b\in \c$ bounds $H$, then~$b$ also bounds the real closure $H^{\operatorname{rc}}\subseteq\c$ of $H$ [ADH, 5.3.2]. {\it In the rest of this subsection~$H$ is a Hardy field.}\/

\begin{lemma}\label{lem:bounded Hardy field ext, 2}
Let $H^*$ be a  $\d$-algebraic Hardy field extension of $H$ and
suppose~$H$ is bounded. Then $H^*$ is also bounded.
\end{lemma}
\begin{proof} By [ADH, 3.1.11] we have $f\in H(x)^{>}$ such that for all~$g\in H(x)^\times$ there are~$h\in H^\times$ and $q\in \Q$ with $g\asymp hf^q$. Hence $H(x)$ is bounded. 
Replacing $H$, $H^*$ by~$H(x)^{\operatorname{rc}}$,~$\Li\!\big(H^*(\R)\big)$, respectively, we arrange that $H$ is real closed 
with $x\in H$, and
$H^*\supseteq\R$ is Liouville closed.
Let $b\in\c$  bound~$H$. Then any  $b^*\in\c$ such that~$\exp_n b \leq b^*$ for all $n$ 
bounds $H^*$, by Lemma~\ref{lem:bounded Hardy field ext, 1}.
\end{proof}

\noindent
In particular, if $H$ is bounded, then so is $\Li\!\big(H(\R)\big)$. We use this   to show:

\begin{lemma}\label{lem:bounded Hardy field ext, 3} Suppose that $H$ is bounded and  $f\in \Calinf$ is hardian over  $H$. Then~$H\<f\>$ is bounded. 
\end{lemma}
\begin{proof}  
Using 
that $f$ remains hardian over the bounded Hardy field~$\Li\!\big(H(\R)\big)$,
  we arrange that $H$ is Liouville closed. The case that~$H\<f\>$ has no element~$>H$ is trivial, so assume we have $y\in H\langle f \rangle$ with~$y>H$. Then $y$ is $\d$-transcendental over~$H$ and the sequence
$y, y^2, y^3,\dots$  is cofinal in $H\<y\>$, by Corollary~\ref{cor:val gp at infty}, so $H\<y\>$ is bounded. Now use that  $f$ is
$\d$-algebraic over $H\<y\>$. 
\end{proof}

\begin{theorem}[{Boshernitzan~\cite[Theorem~14.4]{Boshernitzan82}}]\label{thm:Bosh 14.4}
Suppose $H$ is bounded. Then the perfect hull~$\Ex(H)$ of $H$ is $\d$-algebraic over $H$ and hence  bounded. If
$H\subseteq\Ginf$, then~$\Ex^\infty(H)$ is $\d$-algebraic over $H$; likewise with $\omega$ in place of~$\infty$.
\end{theorem}

\noindent
Using the results above the proof is not difficult. It is omitted in \cite{Boshernitzan82}, but we include it here for the sake of completeness.  First, a lemma also needed for the proof of Theorem~\ref{thm:coinitial in E(H)}:

\begin{lemma}\label{lem:Q bound}
Let  $b\in\c^\times$ bound $H$, let $\phi\ge 0$ in $\Calinf$ satisfy
$\phi^{(n)}\prec b^{-1}$ for all~$n$, and let $r\in \phi\cdot  (\Calinf)^{\preceq}$. Then $Q(r)\prec 1$ for all $Q\in H\{Y\}$ with $Q(0)=0$.
\end{lemma}
 \begin{proof} From  $\phi\in \mathcal I$ we obtain $r\in\mathcal I$, so it is enough that $hr^{(n)}\prec 1$ for all
  $h\in H$ and all~$n$. Now use
the Product Rule and $h\phi^{(n)}\prec hb^{-1}\preceq 1$ for $h\in H^\times$. 
\end{proof}

\begin{proof}[Proof of Theorem~\ref{thm:Bosh 14.4}]
Using Lemma~\ref{lem:bounded Hardy field ext, 2}, replace $H$ by $\Li\!\big(H(\R)\big)$ to arrange that~$H\supseteq\R$ is Liouville closed.
Let $b\in\c$ bound~$H$. Then $b$ also bounds $\Ex(H)$,  by Corollary~\ref{cor:Bosh13.10}. 
 Lemma~\ref{lem:Bosh 14.3} yields~$\phi\ge 0$ in $(\Gom)^\times$ such that~$\phi^{(n)} \prec b^{-1}$ for all~$n$;
 set $r:=\phi\cdot\sin x\in\Gom$. Then~$Q(r) \prec f$ for all~$f\in\Ex(H)^\times$ and $Q\in \Ex(H)\{Z\}$ with~$Q(0)=0$, by Lem\-ma~\ref{lem:Q bound}.

Suppose towards a contradiction that $f\in \Ex(H)$ is   $\d$-transcendental over~$H$, and
set $g:=f+r\in\Calinf$. Then $f$, $g$ are not in a common Hardy field, so $g$ is not hardian over $H$. On the other hand, let $P\in H\{Y\}^{\neq}$.
Then $P(f)\in\Ex(H)^\times$, and by Taylor expansion,
$$P(f+Z)\ =\ P(f)+Q(Z)\quad\text{ where  $Q\in \Ex(H)\{Z\}$ with $Q(0)=0$,}$$ 
so  $P(g)=P(f+r) \sim P(f)$. Hence $g$ is hardian over $H$, a contradiction.

The proof in the case where $H\subseteq\Ginf$ is similar, using the version of Corollary~\ref{cor:Bosh13.10}
for $\Ex^\infty(H)$; similarly for $\omega$ in place of $\infty$.
\end{proof}

\noindent
As to
the existence of   trans\-ex\-po\-nen\-tial hardian germs, we have:

\begin{theorem}\label{thm:Bosh 1.2}
For every $b\in\c$ there is a $\Gom$-hardian germ $y\geq b$. 
\end{theorem}

\noindent
This is Boshernitzan~\cite[Theo\-rem~1.2]{Boshernitzan86}, and leads to~\cite[The\-o\-rem~1.1]{Boshernitzan86}:

\begin{cor}\label{cor:Bosh 1.1} No maximal Hardy field is bounded. 
\end{cor}
\begin{proof}
Suppose $x\in H$, and $b\in\c$ bounds $H$. Take $b^*\in\c$ such that~$b^*\geq \exp_n b$ for all~$n$.
Now Theorem~\ref{thm:Bosh 1.2} yields a $\Gom$-hardian germ $y\geq b^*$.
By Corollary~\ref{cor:Bosh 12.23}, $y$ is $H$-hardian, so $H\langle y\rangle$ is a proper Hardy field extension of~$H$.
\end{proof}

\noindent
The same proof shows also that no $\Ginf$-maximal Hardy field and no $\Gom$-maximal Hardy field is bounded.
In particular (Boshernitzan \cite[Theo\-rem~1.3]{Boshernitzan86}):

\begin{cor}\label{cor:Bosh 1.3}
Every  maximal Hardy field contains a transexponential  germ. Likewise with ``$\Ginf$-maxi\-mal'' or
``$\Gom$-maxi\-mal'' in place of ``maximal''.
\end{cor}

\begin{remark}
For $\Ginf$-Hardy fields, some of the above is in Sj\"odin's~\cite{Sj}, predating~\cite{Boshernitzan82,Boshernitzan86}: if $H$ is a bounded $\Ginf$-Hardy field, then so is~$\Li\!\big(H(\R)\big)$
\cite[The\-o\-rem~2]{Sj};   no maximal $\Ginf$-Hardy field is bounded \cite[The\-o\-rem~6]{Sj}; and  
$E:=\Ex^\infty(\Q)$ is bounded~\cite[The\-o\-rem~10]{Sj}. 
\end{remark}

\noindent
We can now  produce  a differential subfield $K$ of $\Gom[\imag]$ 
containing $\imag$ such that  $\preceq$  restricts to a dominance relation on $K$
making $K$ a $\d$-valued field of $H$-type with constant field~$\C$, yet $K\not\subseteq H[\imag]$ for every  $H$: 
Take a transexponential $\Gom$-hardian germ $z$, 
and~${h\in\R(x)}$ with $0\neq h\prec 1$. Then 
$\varepsilon:=h\ex^{x\imag}\in\mathcal I$, so~${y:=z(1+\varepsilon)\in\Gom[\imag]}$ with~$y\sim_\infty z$.
Lemma~\ref{lem:y H-hardian crit, complex} applied with $H=\R$ shows that~$y$ generates a differential subfield~$K_0:=\R\langle y\rangle$
of $\Gom[\imag]$, and   $\preceq$ restricts to a dominance relation on~$K_0$
making~$K_0$ a $\d$-valued field of $H$-type with constant field~$\R$.
Then~$K:=K_0[\imag]$ is a differential subfield of~$\Gom[\imag]$ with constant field $\C$.
Moreover,~$\preceq$ also restricts to a dominance relation on~$K$,   and this dominance
relation makes $K$
a $\d$-valued field of $H$-type~[ADH, 10.5.15]. We cannot have $K\subseteq H[\imag]$ for any $H$, since~$\Im y=zh\sin x\notin H$.

\subsection*{Lower bounds on $\d$-algebraic hardian germs} 
In this subsection $H$ is a Hardy field. Let $f\in\c$ and $f\succ 1$, $f\ge 0$. Then the germ $\log f\in\c$ also satisfies
$\log f\succ 1$, $\log f\ge 0$. So we may recursively define the germs $\log_n f$ in $\c$ by $\log_0 f:=f$, $\log_{n+1} f:=\log \log_n f$.
(So $\ell_n=\log_n x$ for each $n$.) Lemma~\ref{lem:bounded Hardy field ext, 1} gives exponential upper bounds on $\d$-algebraic $H$-hardian germs.
The next result leads to logarithmic lower bounds on such germs when $H$ is grounded.

\begin{theorem}[{Rosenlicht~\cite[Theorem~3]{Rosenlicht83}}]\label{thm:Ros83}
Suppose $H$ is grounded,  and let~$E$ be a Hardy field extension of $H$ such that $|\Psi_E\setminus\Psi_H|\leq n$ \textup{(}so $E$ is also 
grounded\textup{)}. Then there are $r,s\in\N$ with $r+s\leq n$ such that
\begin{enumerate}
\item[\textup{(i)}] for any $h\in H^>$ with $h\succ 1$ and $\max\Psi_H=v(h^\dagger)$, 
there exists~$g\in E^>$ such that $g\asymp \log_r h$
and
$\max\Psi_E = v(g^\dagger)$;
\item[\textup{(ii)}] for any $g\in E$ there exists $h\in H$ such that $g<\exp_s h$.
\end{enumerate}
\end{theorem}

\noindent
This theorem is most useful in combination with
the following lemma, which is~\cite[Pro\-po\-sition~5]{Rosenlicht83} (and also \cite[Lemma~2.1]{AvdD3} in the context of pre-$H$-fields).

\begin{lemma}\label{lemro83}
Let  $E$ be a Hardy field extension of $H$ such that $\operatorname{trdeg}(E|H)\leq n$. Then~$|\Psi_E\setminus\Psi_H|\leq n$. 
\end{lemma}

\noindent
From [ADH, 9.1.11] we recall that for $f,g\succ 1$ in a Hardy field we have $f^\dagger\preceq g^\dagger$ iff~$\abs{f}\leq\abs{g}^n$
for some~$n\geq 1$. 
Thus by Lemma~\ref{lemro83} and Theorem~\ref{thm:Ros83}:

\begin{cor}\label{cor:Ros83}
Let
$h\in H^>$, $h\succ 1$, and $\max\Psi_H=v(h^\dagger)$. 
Then for any Hardy field extension $E$ of
$H$ with $\operatorname{trdeg}(E|H)\leq n$: $E$ is grounded, and
for all $g\in E$ with~$g\succ 1$ there is an $m\geq 1$ such that $\log_n h  \preceq g^m$
\textup{(}and so~$\log_{n+1} h \prec g$\textup{)}.
\end{cor}

\noindent
Hence for $h$ as in Corollary~\ref{cor:Ros83} and $H$-hardian $y\in\c$, if $y$ is $\d$-algebraic over $H$, then the Hardy field $E=H\langle y\rangle$ is grounded, and there is an $n$ such that~$\log_n h  \prec g$ for all $g\in E$ with $g \succ 1$.
Applying this to  $H=\R(x)$, $h=x$   yields:

\begin{cor}[{Boshernitzan~\cite[Proposition~14.11]{Boshernitzan82}}]\label{cor:Bosh 14.11}
If $y\in\c$ is hardian and $\d$-algebraic over $\R$, then the Hardy field~$E=\R(x)\langle y\rangle$ is grounded, and there is an~$n$ such that $\ell_n  \prec g$
for all $g\in E$ with $g\succ 1$.
\end{cor}

\noindent
Following~\cite{Boshernitzan86} we say that~$y\in\c$ is {\bf translogarithmic} if~$r \leq y \leq \ell_n$ for all~$n$ and all $r\in\R$.
Thus for eventually strictly increasing~$y\succ 1$ in $\c$,  $y$ is
translogarithmic iff its compositional inverse~$y^{\operatorname{inv}}$ is transexponential.
By  Lemma~\ref{lem:Bosh6.5} and Corollary~\ref{cor:Bosh 1.3} there exist $\Gom$-hardian    translogarithmic germs; see also~[ADH, 13.9].\index{germ!translogarithmic}\index{translogarithmic}
 Translogarithmic hardian germs are $\d$-transcendental, by Corollary~\ref{cor:Bosh 14.11}.

\section{Second-Order Linear Differential Equations over Hardy Fields}\label{sec:order 2 Hardy fields}

\noindent
In this section we review Boshernitzan's work~\cite[\S{}16]{Boshernitzan82} on adjoining non-oscillating solutions of second-order linear differential equations to Hardy fields in the light of results from [ADH], and
deduce some consequences about complex exponentials over Hardy fields for use in \cite{ADH5}. {\it Throughout this section~$H$ is a Hardy field.}\/

\subsection*{Oscillation over Hardy fields}
In this subsection we assume $f\in H$ and consider the 
linear differential equation 
\begin{equation}\label{eq:2nd order, 2} \tag{4L}
4Y''+fY\ =\ 0
\end{equation}
over $H$. The factor $4$ is to simplify certain expressions, in conformity with [ADH, 5.2]. 
There  we also defined for any differential field $K$  the function $\omega\colon K \to K$ given by $\omega(z)=-2z'-z^2$, and
the function $\sigma\colon K^\times \to K$ given by
$\sigma(y)=\omega(z)+y^2$ for $z:= -y^\dagger$.
We define likewise 
$$\omega\ :\ \Go[\imag]\to \Gz[\imag], \qquad \sigma\ :\ \c^2[\imag]^\times \to \Gz[\imag]$$ 
by
$$\omega(z)\ =\ -2z'-z^2\quad\text{ and }\quad\sigma(y)\ =\ \omega(z)+y^2\text{ for $z:= -y^\dagger$.}$$
Note that $\omega(\Go)\subseteq\Gz$ and $\sigma\big( (\c^2)^\times \big) \subseteq\Gz$, and
$\sigma(y) = \omega(z+y\imag)$ for  $z:= -y^\dagger$. 
To clarify the
role of $\omega$ and $\sigma$ in connection with second-order linear differential equations, suppose
$y\in \Gt$  is a non-oscillating solution to \eqref{eq:2nd order, 2} with $y\neq 0$.
Then~${z:=2y^\dagger\in \Go}$ satisfies
$-2z'-z^2=f$, so $z$ generates a Hardy field~$H(z)$
with $\omega(z)=f$,  by Proposition~\ref{singer}, which in turn
yields a Hardy field $H(z,y)$ with~$2y^\dagger=z$.
Thus $y_1:=y$ lies in a Hardy field extension of $H$. From Lemma~\ref{lem:2nd lin indep sol} and 
Proposition~\ref{prop:Hardy field exts}(iv) we also obtain a solution $y_2$ of \eqref{eq:2nd order, 2}
in a Hardy field
extension of $H\<y_1\>=H(y,z)$ such that~$y_1$,~$y_2$ are $\R$-linearly independent; see also~\cite[Theo\-rem~2, Corollary~2]{Ros}. 
This shows: 
 
\begin{prop}\label{prop:2nd order Hardy field}
If $f/4$ does not generate oscillation, then $\Dx(H)$ contains $\R$-linearly independent solutions $y_1$, $y_2$ to \eqref{eq:2nd order, 2}.  
\end{prop}

\noindent
Indeed, if $f/4$ does not generate oscillation, then $\Dx(H)$ contains   solutions $y_1$, $y_2$ of~\eqref{eq:2nd order, 2}
with $y_1,y_2>0$ and $y_1\prec y_2$.  Here $y_1$ is determined up to multiplication by a factor in $\R^>$; we call such $y_1$ a {\bf principal solution} of~\eqref{eq:2nd order, 2}.\index{solution!principal} (Lemmas~\ref{lem:admissible pair, unique},~\ref{lem:admissible pair}.) 
 
\begin{notation}
Let $K$ be a differential field. Then $K[\der]$ denotes the ring of linear differential operators over $K$; see [ADH, 5.1]. Let $A\in K[\der]$.
The {\it twist}\/ of $A$ by $b\in K^\times$ is~$A_{\ltimes b}:=b^{-1}Ab\in K[\der]$. We say that
{\it $A$ splits over $K$}\/ if $A=a(\der-b_1)\cdots(\der-b_n)$ for some $a\in K^\times$, $b_1,\dots,b_n\in K$. If $A$ splits over $K$,
then so does $A_{\ltimes b}$ for each~$b\in K^\times$.
\end{notation}

\noindent 
By [ADH, p.~259], with~$A:=4\der^2+f\in H[\der]$ we have
$$\text{$4y''+fy=0$ for some $y\in H^\times$}\ \Rightarrow\ \text{$A$ splits over $H$} \ \Longleftrightarrow\ f\in\omega(H).$$
To simplify the discussion we now also introduce the subset \label{p:baromega}
$$\bar{\omega}(H):=  \big\{ f\in H: \text{$f/4$ does not generate oscillation} \big\}$$
of $H$. If $E$ is a Hardy field extension of $H$, then $\bar\omega(E)\cap H=\bar\omega(H)$.
By Corollary~\ref{cor:gen osc closed upward}, $\bar{\omega}(H)$ is downward closed, and  $\omega(H)\subseteq \bar{\omega}(H)$
by the discussion following~\eqref{eq:Riccati} in Section~\ref{sec:differentiable germs}.

\begin{cor}  \label{cor:omega(H) downward closed}
If $H$ is $\d$-perfect, then 
\[\omega(H)\	=\ \bar\omega(H) 
		\	=\ \big\{ f\in H :\ \text{$4y''+fy=0$ for some $y\in H^\times$} \big\},\]
and $\omega(H)$ is  downward closed in $H$.
\end{cor}

\noindent
Lemma~\ref{lem:2nd order inhom} and Proposition~\ref{prop:2nd order Hardy field} also yield:

\begin{cor}\label{cor:2nd order Hardy field, inhom}
If $f\in\bar\omega(H)$, then each~${y\in\Gt}$   such that $4y''+fy\in H$ is in~$\Dx(H)$.
\end{cor}

\noindent
Next some consequences of Proposition~\ref{prop:2nd order Hardy field} for more general linear differential equations of
order $2$: Let $g,h\in H$, and consider the linear differential equation
\begin{equation}\label{eq:2nd order, 3}\tag{$\tilde{\operatorname{L}}$}
Y''+gY'+hY\ =\ 0
\end{equation}
over $H$. An easy induction on $n$ shows that for a solution $y\in \Gt$ of \eqref{eq:2nd order, 3} we have~$y\in \mathcal{C}^n$ with $y^{(n)}\in Hy+Hy'$ for all $n$, so
$y\in\Calinf$. To reduce \eqref{eq:2nd order, 3}
to an equation~\eqref{eq:2nd order, 2} we take 
$$f\ :=\ \omega(g)+4h\ =\ -2g'-g^2+4h\in H,$$ take $a\in \R$, and take a representative of~$g$ in~$\Cao$, also denoted by $g$, and let $G\in (\Gt)^\times$ be 
the germ of 
$$t\mapsto \exp\!\left(-\frac{1}{2}\int_a^t g(s)\,ds\right)\qquad(t\ge a).$$ 
This gives an isomorphism $y\mapsto Gy$ from the $\R$-linear space of 
solutions of \eqref{eq:2nd order, 2} in~$\Gt$ onto the
$\R$-linear space of solutions of \eqref{eq:2nd order, 3} in~$\Gt$,
and $y\in \Gt$ oscillates iff~$Gy$ oscillates.
By Proposition~\ref{prop:Hardy field exts}, $G\in\Dx(H)$. 
Using $\frac{f}{4}=-\frac{1}{2}g'-\frac{1}{4}g^2+h$ we now obtain the following
germ version of Corollary~\ref{coroscgen}:
 
\begin{cor}\label{cor:char osc}
The following are equivalent:
\begin{enumerate}
\item[\textup{(i)}] some solution in $\Gt$ of \eqref{eq:2nd order, 3}  oscillates;
\item[\textup{(ii)}] all nonzero solutions in $\Gt$ of \eqref{eq:2nd order, 3} oscillate;
\item[\textup{(iii)}]  $-\frac{1}{2}g'-\frac{1}{4}g^2+h$ generates oscillation.
\end{enumerate}
Moreover, if $-\frac{1}{2}g'-\frac{1}{4}g^2+h$ does not generate oscillation, then all solutions of~\eqref{eq:2nd order, 3} in $\Gt$ belong to $\Dx(H)$. 
\end{cor}

\noindent
Set $A:=\der^2+ g\der + h$, and let $f=\omega(g)+4h$, $G$ be as above. Then $A_{\ltimes G}=\der^2 + \frac{f}{4}$.
Thus by combining Corollary~\ref{cor:2nd order Hardy field, inhom} and Corollary~\ref{cor:char osc} we obtain: 

\begin{cor}\label{cor:char osc, 1}
If \eqref{eq:2nd order, 3} has no oscillating solution in $\Gt$, and ${y\in\Gt}$ is such that ${y''+gy'+hy\in H}$, then $y\in\Dx(H)$.  
\end{cor}

\noindent
The next corollary follows from Proposition~\ref{prop:2nd order Hardy field} and [ADH, 5.1.21]: 
 
\begin{cor}\label{cor:char osc, 2}
The following are equivalent, for  $A\in H[\der]$ and $f$ as above:
\begin{enumerate}
\item[\textup{(i)}] $f/4$ does not generate oscillation;
\item[\textup{(ii)}] $A$ splits over some Hardy field extension of $H$;
\item[\textup{(iii)}] $A$ splits over $\Dx(H)$.
\end{enumerate}
\end{cor}

\noindent
For  $A\in H[\der]$ and $f$ as before we have  $A_{\ltimes G}=\der^2 + \frac{f}{4}$ and $G^\dagger=-\frac{1}{2}g\in H$, so: 

\begin{cor}\label{cor:char osc, 3}
$A$ splits over $H[\imag]$ $\Longleftrightarrow$   $\der^2+\frac{f}{4}$ splits over $H[\imag]$.
\end{cor}

\noindent
Proposition~\ref{prop:2nd order Hardy field} and its corollaries \ref{cor:2nd order Hardy field, inhom}--\ref{cor:char osc, 1}
are from \cite[Theorems~16.7, 16.8]{Boshernitzan82}, and Corollary~\ref{cor:omega(H) downward closed} is essentially 
\cite[Lemma~17.1]{Boshernitzan82}.

\medskip
\noindent
Proposition~\ref{prop:2nd order Hardy field} applies only when \eqref{eq:2nd order, 2} has a solution in $(\Gt)^\times$.
Such a solution might not exist, but \eqref{eq:2nd order, 2} does have $\R$-linearly independent solutions~${y_1, y_2\in \Gt}$, so~$w:= y_1y'_2-y'_1y_2\in \R^\times$. 
Set $y:= y_1+ y_2\imag$. Then~$4y'' + fy=0$ and~$y\in\Gt[\imag]^\times$, and for $z=2y^\dagger\in \Go[\imag]$   we have $-2z'-z^2=f$. Now
\begin{align*} z\ =\ \frac{2y_1'+2\imag y_2'}{y_1+\imag y_2}\ &=\
\frac{2y_1'y_1+2y_2'y_2- 2\imag(y_1'y_2-y_1y_2')}{y_1^2+y_2^2}\ =\ \frac{2(y_1'y_1+y_2'y_2)+2\imag w}{y_1^2+y_2^2},\\
\text{so }\ \Re z\ &=\ \frac{2(y_1'y_1+y_2'y_2)}{y_1^2+y_2^2}\in \Go, \qquad \Im z\ =\ \frac{2w}{y_1^2+y_2^2}\in \Gt.
\end{align*} 
Thus $\Im z\in (\Gt)^\times$ and $(\Im z)^\dagger=-\Re z$,  and so
$$\sigma(\Im z)\ =\ \omega\big({-(\Im z)^\dagger}+(\Im z)\imag\big)\ =\ 
\omega(z)\ =\ f\qquad\text{in $\Go$.}$$
Replacing~$y_1$ by $-y_1$ changes 
$w$ to $-w$; this way we can arrange ${w>0}$, so $\Im z> 0$.

\medskip
\noindent
Conversely, every $u\in  (\c^2)^\times$ such that $u>0$ and $\sigma(u)=f$ arises in this way. 
To see this, suppose we are given such $u$,  take
$\phi\in \c^3$ with $\phi'=\frac{1}{2}u$, and set
$$y_1\ :=\  \frac{1}{\sqrt{u}}\cos \phi, \qquad y_2\ :=\  \frac{1}{\sqrt{u}}\sin \phi \qquad \text{(elements of $\c^2$)}.$$
Then $\wr(y_1, y_2)=1/2$, and $y_1$, $y_2$ solve  \eqref{eq:2nd order, 2}.
To see the latter, consider  
$$y\ :=\ y_1+y_2\imag\ =\ \frac{1}{\sqrt{u}}\ex^{\phi\imag}\in \c^2[\imag]^\times$$ and note that $z:=2y^\dagger$ satisfies $$\omega(z)\ =\ \omega(-u^\dagger+u\imag)\ =\ \sigma(u)\ =\ f,$$ hence
$4y''+fy=0$. The computation above shows~$\Im z = 1/(y_1^2+y_2^2) = u$.
We have $\phi'> 0$, so either $\phi>\R$ or $\phi-c\prec 1$ for some $c\in\R$, with $\phi>\R$ iff $f/4$ generates oscillation. 
As to uniqueness of the above pair $(y_1,y_2)$, we have:

\begin{lemma}\label{lem:Im z}
Suppose $f\notin\bar\omega(H)$.
Let $\tilde y_1,\tilde y_2\in \c^2$ be  $\R$-linearly independent solutions of \eqref{eq:2nd order, 2} 
with $\wr(\tilde y_1, \tilde y_2)=1/2$.  Set $\tilde y:=\tilde y_1+\tilde y_2\imag$, $\tilde z:=2\tilde y^\dagger$. Then
$$ \Im \tilde z =u\quad\Longleftrightarrow\quad \tilde y=\ex^{\theta\imag}y \text{ for some $\theta\in\R$.}$$
\end{lemma}
\begin{proof}
If $\tilde y=\ex^{\theta\imag} y$ ($\theta\in\R$), then clearly $\tilde z = 2 \tilde y^\dagger = 2 y^\dagger = z$, hence $\Im z = \Im \tilde z$.
For the converse, let $A$ be the invertible $2\times 2$ matrix  with real entries and~$Ay=\tilde y$; here~$y=(y_1, y_2)^{t}$ and
$\tilde y=(\tilde y_1, \tilde y_2)^{t}$, column vectors with entries in $\c^2$.  As in the proof of~[ADH, 4.1.18], 
$\wr(y_1, y_2)=\wr(\tilde y_1, \tilde y_2)$ yields $\det A=1$.

Suppose $\Im \tilde z=u$, so~$y_1^2+y_2^2 = \tilde y_1^2+\tilde y_2^2$.  Choose~$a\in\R$ and representatives for $u$, $y_1$, $y_2$, $\tilde y_1$, $\tilde y_2$ in~$\c_a$, denoted by the same symbols, such that in~$\c_a$ we have $Ay=\tilde y$  and $y_1^2+y_2^2 = \tilde y_1^2+\tilde y_2^2$, and $u(t)\cdot \big(y_1(t)^2+y_2(t)^2\big)=1$ for all $t\ge a$. With $\dabs{\, \cdot\, }$  the usual euclidean norm on $\R^2$, 
we then have $\dabs{Ay(t)}=\dabs{y(t)}=1/\sqrt{u(t)}$ for~$t\geq a$.
Since~$f/4$ generates oscillation, we have
$\phi> \R$, and we conclude that~$\dabs{Av}=1$ for all $v\in\R^2$ with $\dabs{v}=1$. It is well-known that then $A=\left(\begin{smallmatrix} \cos \theta & -\sin\theta \\ \sin\theta & \cos\theta\end{smallmatrix}\right)$ with~$\theta\in\R$  (see, e.g., \cite[Chapter~XV, Exercise~2]{Lang}), so
$\tilde y=\ex^{\theta\imag}y$.
\end{proof}

\noindent
The observations above will be used in the proof  of Theorem~\ref{upo}   below. 
We finish this subsection with miscellaneous historical remarks  (not used later):

 \begin{remarks}
The connection between the second-order linear differential equation \eqref{eq:2nd order, 2} and the second-order
non-linear differential equation~$\sigma(u)=f$    
was first investigated by Kummer~\cite{Kummer} in 1834. Appell~\cite{Appell}
noted that the linear differential equation
$$Y'''+fY'+(f'/2)Y\ =\ 0$$ has $\R$-linearly independent solutions~$y_1^2, y_1y_2, y_2^2\in\Calinf$, though some cases were known earlier~\cite{Clausen,Liouville39}; in particular, $1/u=y_1^2+y_2^2$ is a solution.  
Hartman~\cite{Hartman61, Hartman73} investigates monotonicity properties of $y_1^2+y_2^2$.
Steen~\cite{Steen} in 1874, and independently Pinney~\cite{Pinney}, remarked that $r:=1/\sqrt{u}=\sqrt{y_1^2+y_2^2}\in\Calinf$ satisfies~$4r''+fr=4w^2/r^3$ where $w:=y_1y_2'-y_1'y_2\in \R^\times$. (See also~\cite{Redheffer}.)
\end{remarks}

\subsection*{Complex exponentials over Hardy fields}
We now use some  of the above  to prove an extension theorem for Hardy fields (cf. \cite[Lem\-ma~11.6(6)]{Boshernitzan82}).
Recall that the $H$-asymptotic field extension $K:=H[\imag]$ of~$H$ is a differential subring of $\Calinf[\imag]$.  
 
\begin{prop}\label{prop:cos sin infinitesimal, 1}
If $\phi\in H$ and $\phi\preceq 1$, then $\cos \phi,\sin\phi\in\Dx(H)$. 
\end{prop}
\begin{proof}
Replacing $H$ by $\Dx(H)$ we  arrange $\Dx(H)=H$.
Then by Pro\-po\-si\-tion~\ref{prop:Hardy field exts}, $H\supseteq\R$ is a Liouville closed $H$-field, and
by Corollary~\ref{cor:omega(H) downward closed}, $\omega(H)$ is downward closed.
Hence by \cite[Lemma~1.2.17]{ADH4},
there is for all $\phi\preceq 1$ in $H$ a (necessarily unique) $y\in K$ with $y\sim 1$ and $y^\dagger=\phi'\imag$.
Let now $\phi\in H$ and~${\phi\preceq 1}$.  Then $(\ex^{\phi\imag})^\dagger=\phi'\imag\in K^\dagger$, so~$\cos\phi+\imag\sin\phi=\ex^{\phi\imag}\in K$ using $K\supseteq\C$.  Thus~$\cos\phi,\sin\phi\in H$.
\end{proof}

\begin{cor}\label{cor:cos sin infinitesimal, 1}
Let  $\phi\in H$ and $\phi\preceq 1$. Then $\cos \phi$, $\sin\phi$ generate a $\d$-algebraic Hardy field extension 
$E:=H(\cos\phi,\sin\phi)$  of $H$. If
$H$ is a $\Ginf$-Hardy field, then so is $E$, and likewise with $\Gom$ in place of~$\Ginf$. 
\end{cor}

\noindent
In \cite[Chapter~4]{ADH4} we
sometimes assume $\I(K)\subseteq K^\dagger$, a condition that we consider more closely in 
the next proposition:

\begin{prop}\label{prop:cos sin infinitesimal, 2}
Suppose $H\supseteq\R$ is closed under integration, that is, $\der(H)=H$.  
Then the following conditions are equivalent:
\begin{enumerate}
\item[\textup{(i)}] $\I(K)\subseteq K^\dagger$;
\item[\textup{(ii)}] $\ex^{f}\in K$ for all $f\in K$ with $f\prec  1$;
\item[\textup{(iii)}] $\ex^\phi,\cos \phi, \sin\phi\in H$ for all $\phi\in H$ with $\phi\prec 1$.
\end{enumerate}
\end{prop}
\begin{proof}
Assume (i), and let $f\in K$, $f\prec  1$.
Then $f' \in \I(K)$, so we have~$g\in K^\times$ with~$f'=g^\dagger$
and thus $\ex^{f}=cg$ for some $c\in\mathbb C^\times$. Therefore~$\ex^{f}\in K$.
This shows~(i)~$\Rightarrow$~(ii), and (ii)~$\Rightarrow$~(iii) is clear.
Assume (iii), and let $f\in \I(K)$. Then~$f=g+h\imag$, $g, h\in \I(H)$. Taking 
$\phi, \theta\prec 1$ in~$H$
with~$\phi'=g$ and $\theta'=h$,   
$$\exp(\phi + \theta\imag)\ =\ \exp(\phi)\big(\!\cos (\theta) + \sin(\theta)\imag\big)\in H[\imag]\ =\ K$$
has the property that $f=\big(\!\exp(\phi+\theta\imag)\big){}^\dagger\in K^\dagger$. 
This shows (iii)~$\Rightarrow$~(i).
\end{proof}

\noindent
From Propositions~\ref{prop:cos sin infinitesimal, 1} and \ref{prop:cos sin infinitesimal, 2} we obtain:

\begin{cor}\label{cor:cos sin infinitesimal}
If $H$ is $\d$-perfect, then   $\I(K)\subseteq K^\dagger$. 
\end{cor}

\subsection*{Some special subsets of $H$}
Let $f\in H$.
Principal solutions of~\eqref{eq:2nd order, 2} arise from certain distinguished solutions of the non-linear (Riccati) equation $\omega(z)=f$.
To explain this and for later use we recall from [ADH, 11.8] some special subsets of $H$:   
\begin{align*}
\Upg(H)	&\ :=\  \big\{ h^{\dagger}:\ h\in H,\, h\succ 1\big\}\ \subseteq\  H^{>}, \\
\Upl(H)	&\ :=\  \big\{ {-h^{\dagger\dagger}}:\ h\in H,\, h\succ 1\big\}, \\
\Upd(H)	&\ :=\  \big\{ {-h^{\prime\dagger}}:\  h\in H,\, 0\neq h\prec 1\big\}.
\end{align*}
Here $\Upl(H)= - \Upg(H)^\dagger$ and $\Upd(H)$ are disjoint, and
$\omega\colon H\to H$ and $\sigma\colon H^\times\to H$ 
are strictly increasing on $\Upl(H)$ and $\Upg(H)$, respectively. If~$H\supseteq\R$ is Liouville closed, then~$\Upg(H)$ is upward closed, 
$\Upl(H)$ is downward closed, $H=\Upl(H)\cup\Upd(H)$, and~$\omega\big(\Upl(H)\big)=\omega\big(\Upd(H)\big)=\omega(H)$; see~[ADH, 11.8.13, 11.8.19, 11.8.20, 11.8.29].

\begin{lemma}\label{lem:omegabar(H)}
Suppose $H$ is $\d$-perfect and $f/4$ does not generate oscillation, and let~${y\in H}$ be a principal solution  of~\eqref{eq:2nd order, 2}.
Then $2y^\dagger$ is the unique $z\in \Upl(H)$ such that~$\omega(z)=f$.
\end{lemma}
\begin{proof} We already know $\omega(2y^\dagger)=f$, and as $\omega$ is strictly increasing on $\Upl(H)$, it remains to show that  $2y^\dagger\in \Upl(H)$.  For this
take~$h\in H$ with $h'=1/y^2$. Then $h\succ 1$ by Corollary~\ref{cor:I1 I2}, hence
$1/y^2\in\Upg(H)$, and thus
$2y^\dagger=-(1/y^2)^\dagger\in \Upl(H)$.
\end{proof}

\noindent
Combining Lemma~\ref{lem:omegabar(H)} with the remarks after Proposition~\ref{prop:2nd order Hardy field} yields:

\begin{cor}\label{cor:omegabar(H)}
If $H$ is $\d$-perfect and $f/4$ does not generate oscillation, then~\eqref{eq:2nd order, 2} has solutions $y_1, y_2\in H$ such that $y_1,y_2>0$, $y_1\prec y_2$,
 $2y_1^\dagger\in\Upl(H)$,  and $2y_2^\dagger\in\Upd(H)$ \textup{(}and thus $y_2'>0$ in view of 
 $-(\ex^x)^{\dagger\dagger}=0\in \Upl(H)$\textup{)}.
\end{cor}

\begin{remarkNumbered}\label{rem:Hardy inequ, 2}
Suppose $f/4>0$ does not generate oscillation. Remark~\ref{rem:Hardy inequ, 1} and Corollary~\ref{cor:omegabar(H)} yield $a\in\R$, a representative of $f$ in $\mathcal C_a$, also denoted by~$f$, and a constant  $C\in\R^{\geq}$
such that the inequality \eqref{eq:Hardy inequ} holds for all $u\in\mathcal C^1_a$ with~$u(a)=0$.
(This is not used later.)
\end{remarkNumbered}

\noindent
We have $\omega(H) < \sigma\big(\Upg(H)\big)$ by [ADH, remark before 11.8.29].
Recall that $\bar{\omega}(H)$ is downward closed and $\omega(H)\subseteq\bar{\omega}(H)$, with equality 
for $\d$-perfect~$H$. (Corollary~\ref{cor:omega(H) downward closed}.)
This yields a property of~$\bar{\omega}(H)$ used in the proof of Corollary~\ref{cor:FLW}:

\begin{lemma}\label{lem:Upgcapbaromega}
$\bar{\omega}(H) < \Upg(H)$.
\end{lemma}
\begin{proof}
We have $\bar{\omega}(H)\subseteq\bar{\omega}\big(\!\Dx(H)\big)$ and $\Upg(H)\subseteq \Upg\big(\!\Dx(H)\big)$.
Thus,
replacing $H$ by~$\Dx(H)$,  
we arrange that $H$ is $\d$-perfect. Hence $H\supseteq\R$ is Liouville closed and~$\bar{\omega}(H)=\omega(H)$. From $x^{-1}=x^\dagger\in \Upg(H)$ and $\sigma(x^{-1})=2x^{-2}\asymp (x^{-1})'\prec\ell^\dagger$ for all $\ell\succ 1$ in $H$ we obtain $\Upg(H)\subseteq \sigma\big(\Upg(H)\big){}^\uparrow$,
so $\omega(H)<\Upg(H)$.
\end{proof}

\noindent
Suppose $H$ has asymptotic integration. Then by [ADH, 11.8.16, 11.8.30]:
$H$ is $\upl$-free iff there is no~$\upl\in H$ such that~$\Upl(H)<\upl<\Upd(H)$, and
$H$ is $\upo$-free iff there is no~$\upo\in H$ such that~$\omega\big(\Upl(H)\big)<\upo<\sigma\big(\Upg(H)\big)$.   
By  [ADH, 11.7.3], if~$H$ is $\upo$-free, then $H$ is $\upl$-free.
 If $H\supseteq\R$ is
 Liouville closed, then $H$ is $\upl$-free, and 
 $$H \text{ is $\upo$-free}\ \Longleftrightarrow\  \omega(H)^\downarrow=H\setminus\sigma\big(\Upg(H)\big){}^\uparrow.$$

\subsection*{Determining $\bar{\omega}(H)$} 
In [ADH, 16.3] we   introduced the concept of a {\em $\HLO$-cut}\/ in a pre-$H$-field $F$:
these are the triples $(\I,\Upl,\Upo)$ of subsets of $F$ such that
$$(\I,\Upl,\Upo)={\big(\I(E)\cap F, \Upl(E)^\downarrow\cap F, \omega(E)^\downarrow\cap F\big)}$$
for some $\upo$-free $H$-field extension $E$ of $F$.
Every
pre-$H$-field has exactly one or exactly two $\HLO$-cuts~[ADH, remark before 16.3.19].
By~[ADH, 16.3.14, 16.3.16]: 

\begin{lemma}\label{lem:d-perfect HLO-cut}
Suppose $H$ is $\d$-perfect.  Then $\big(\I(H), \Upl(H), \bar{\omega}(H) \big)$ is a $\HLO$-cut in~$H$, and this is the unique $\HLO$-cut in $H$ iff $H$ is $\upo$-free.
\end{lemma}

\noindent
Thus in general, 
$$\big(\I\!\big(\!\Dx(H)\big)\cap H, \, \Upl\big(\!\Dx(H)\big)\cap H,\, \bar{\omega}(H) \big)$$   is a $\HLO$-cut  in $H$, and   hence
$\bar{\omega}(H) < \sigma\big(\Upg(H)\big){}^\uparrow$ (see~[ADH, p.~692]).
The classification of $\HLO$-cuts in $H$ from [ADH, 16.3]  can be used to narrow down the  possibilities for~$\bar{\omega}(H)$.
For this we recall the trichotomy from [ADH, 9.2.16]:  either $H$ has a gap, or $H$ is grounded, or $H$ has asymptotic integration.
With $\mathcal O=$~valuation ring of $H$ and $\smallo=$~maximal ideal of~$\mathcal O$, we have:

\begin{lemma}\label{lem:baromega(H) for H without as int} 
Let  $\phi\in H^{>}$ be such that $v\phi\notin (\Gamma_H^{\neq})'$. Then 
$$\bar{\omega}(H)\  =\  \omega(-\phi^\dagger) + \phi^2 \smallo^\downarrow \quad\text{or}\quad
\bar{\omega}(H)\  =\  \omega(-\phi^\dagger) + \phi^2 \mathcal O^\downarrow.
$$
The first alternative holds if $H$ is grounded, and the second alternative holds if 
$v\phi$ is a gap in  $H$ with $\phi\asymp b'$ for some $b\asymp 1$ in $H$.
\end{lemma}
\begin{proof} Either $v\phi=\max \Psi_H$ or $v\phi$ is a gap in $H$, by [ADH, 9.2]. 
The remark following Lemma~\ref{lem:d-perfect HLO-cut} yields  an $\HLO$-cut $(I,\Upl,\Upo)$ in $H$ where~$\Upo=\overline{\omega}(H)$. Now use the proofs of~[ADH, 16.3.11, 16.3.12, 16.3.13] together with the transformation formulas~[ADH, (16.3.1)] for $\HLO$-cuts.
\end{proof}

\noindent
By~[ADH, 16.3.15] we have:
 
\begin{lemma}\label{lem:baromega(H)}
If $H$ has   asymptotic integration and the set $2\Psi_H$ does  not have a  supremum in $\Gamma_H$, then
$$\bar{\omega}(H)\ =\ \omega\big(\Upl(H)\big){}^\downarrow\ =\ \omega(H)^\downarrow\quad\text{or}\quad \bar{\omega}(H)\ =\ H\setminus \sigma\big(\Upg(H)\big){}^\uparrow.$$
\end{lemma}

\begin{cor}\label{omuplosc}  
Suppose $H$ is $\upo$-free. Then 
$$\bar{\omega}(H)\ =\ \omega\big(\Upl(H)\big){}^\downarrow\ =\ \omega(H)^\downarrow \ = \  H\setminus \sigma\big(\Upg(H)\big){}^\uparrow.$$ 
\end{cor}
\begin{proof} By [ADH, 11.8.30] we have $\omega\big(\Upl(H)\big){}^\downarrow= \omega(H)^\downarrow=  H\setminus \sigma\big(\Upg(H)\big){}^\uparrow$. It follows from [ADH, 9.2.17] that $2\Psi_H$ has no  supremum in $\Gamma_H$. Now use Lem\-ma~\ref{lem:baromega(H)}.
\end{proof}

\subsection*{Non-oscillation and compositional conjugation} 
Which ``changes of variable'' preserve the general form of the linear
differential equation \eqref{eq:2nd order, 2}? The next lemma and Corollary~\ref{cor:chvar}   (used in the proof of our main Theorem~\ref{upo}) give an answer.

\begin{lemma}\label{chvar} 
Let   $K$ be a differential ring, $f\in K$, and   $P:= 4Y''+fY\in K\{Y\}$.
Then for $g\in K^\times$ and $\phi:= g^{-2}$ we have
$$ g^3P_{\times g}^\phi \ =\ 4Y'' + g^3P(g)Y.$$
\end{lemma} 
\begin{proof} Let $g,\phi\in K^\times$. Then 
\begin{align*} P_{\times g}\ &=\ 4gY'' + 8g'Y' + (4g'' + fg)Y\ =\ 4gY'' + 8g'Y' + P(g)Y, \quad \text{so}\\
 P_{\times g}^\phi\ &=\ 4g(\phi^2 Y'' + \phi'Y') + 8g'\phi Y' + P(g)Y\\ 
 &=\ 4g\phi^2 Y'' + (4g\phi' +8g'\phi)Y' + P(g)Y.\end{align*}
Now $4g\phi' + 8g'\phi=0$ is equivalent to $\phi^{\dagger}=-2g^\dagger$, which holds for $\phi=g^{-2}$.
For this~$\phi$ we get $P_{\times g}^\phi(Y)=
g^{-3}\big(4Y'' + g^3P(g)Y\big)$, that is, $g^3P_{\times g}^\phi(Y) = 4Y'' + g^3P(g)Y$. 
\end{proof}

\noindent
Now let  $\ell\in \Go$ be such that $\ell> \R$
and~$\phi:=\ell'\in H$. We use the superscript~$\circ$ as in the subsection
on compositional conjugation in Hardy fields of Section~\ref{sec:Hardy fields}. Let
 $P := 4Y''+fY$ where~$f\in H$. Recall that if $y\in \Gt[\imag]$ and   $4y'' + fy=0$,  then $y\in \Calinf[\imag]$.   Towards using Lemma~\ref{chvar}, 
suppose $\phi=g^{-2}$, $g\in H^\times$, and set $h:= \big(g^3P(g)\big){}^\circ\in H^\circ$.
We then obtain the following reduction of solving the differential equation  \eqref{eq:2nd order, 2} to solving a similar equation over $H^\circ$:

\begin{cor}\label{cor:chvar}
Let $y\in\Gt[\imag]$. Then $z:=(y/g)^\circ\in \Gt[\imag]$, and  
$$4y''+fy\ =\ 0\quad \Longleftrightarrow\quad 
4z''+h z\ =\ 0.$$
In particular, $f/4$ generates oscillation iff 
$h/4$
does. 
\end{cor}

\begin{remark}
Corollary~\ref{cor:chvar} is a special case of a result of Kummer~\cite{Kummer}; cf.~\cite[\S{}11]{Boruvka} and~\cite[\S{}2]{Willett}.
Lie~\cite{Lie} and St\"ackel~\cite{Staeckel} proved uniqueness results about the transformation $y\mapsto (y/h)^\circ$; see~\cite[\S{}22]{Boruvka}.
\end{remark}

\noindent
Consider   the increasing bijection
$$f\mapsto \Phi(f)\ :=\ \big(\big(f-\omega(-\phi^\dagger)\big)/\phi^2\big){}^\circ\ \colon\ H\to H^\circ,$$ 
and note that
$$g^3P(g)\  =\  g^3(4g''+fg)\ =\ \big(f-\omega(-\phi^\dagger)\big)/\phi^2,$$ so $h=\Phi(f)$.

\begin{lemma}\label{lem:baromega comp conj}
 $\Phi\big(\bar{\omega}(H)\big)=\bar{\omega}(H^\circ)$.
\end{lemma}
\begin{proof}
First replace $H$ by its real closure to arrange that $H$ is real closed,   then take~$g\in H^\times$ with $g^{-2}=\phi$, and use
the remarks above. 
\end{proof}

\noindent
The  bijection
$$y\mapsto (y/\phi)^\circ\ \colon\  H\to H^\circ$$ restricts to  bijections
$\I(H)\to\I(H^\circ)$ and $\Upg(H)\to\Upg(H^\circ)$, and the
  bijection
$$z\mapsto  \big( (z+\phi^\dagger) /\phi\big){}^\circ\ \colon\ H\to H^\circ$$
restricts to bijections $\Upl(H)\to\Upl(H^\circ)$ and $\Upd(H)\to\Upd(H^\circ)$.
(See the   transformation formulas in [ADH, p.~520].)
Then for $y\in H^\times$, $z\in H$ we have
$$\sigma\big((y/\phi)^\circ\big)\ =\  \Phi\big(\sigma(y)\big),\qquad \omega\big( \big( (z+\phi^\dagger) /\phi\big){}^\circ \big)\  =\  \Phi\big(\omega(z)\big).$$
(See the formulas in [ADH, pp.~518--519].) 
Hence $\Phi$ also restricts to  bijections  
$$\sigma(H^\times)\to\sigma\big((H^\circ)^\times\big),\quad  
\sigma\big(\I(H)^{\ne}\big)\to \sigma\big(\I(H^\circ)^{\ne}\big), \quad \sigma\big(\Upg(H)\big)\to \sigma\big(\Upg(H^\circ)\big),$$
and
$$\omega(H)\to\omega(H^\circ),\quad \omega\big(\Upl(H)\big)\to \omega\big(\Upl(H^\circ)\big), \quad 
\omega\big(\Upd(H)\big)\to \omega\big(\Upd(H^\circ)\big).$$

\noindent
To illustrate the above, we  use it to prove the Fite-Leighton-Wintner oscillation criterion for self-adjoint second-order linear differential
equations over $H$ \cite{Fite,Leighton50,Wintner49}. (See also \cite[\S{}2]{Hinton} and \cite[p.~45]{Swanson}.)
For this, let  
  $A=f\der^2+f'\der+g$ where~$f\in H^\times$ and~$g\in H$. For each $h\in\c$ we choose  a germ $\int h$ in $\c^1$ such that~$(\int h)'=h$.

\begin{cor}\label{cor:FLW}
Suppose~$\int f^{-1}>\R$ and $\int g>\R$. Then $A(y)=0$ for some oscillating $y\in\Calinf$.
\end{cor}
\begin{proof}
We arrange that $H\supseteq\R$ is  Liouville closed. Then $f^{-1},g\in\Upg(H)$ by [ADH, 11.8.19].   Note that  $\phi:=f^{-1}$ is active in $H$.
 Put~$B:=4\phi A_{\ltimes \phi^{1/2}}$, so $B=4\der^2+h$ with~$h:=\omega(-\phi^\dagger)+4g\phi$. 
Then
\begin{align*}
\text{$A(y)=0$ for some oscillating~$y\in\mathcal C^2$}	&\ \Longleftrightarrow\ \text{${B(z)=0}$ for some oscillating $z\in\mathcal C^2$}\\
													&\ \Longleftrightarrow\ \text{$h\notin\bar{\omega}(H)$,}
\end{align*}
by Corollary~\ref{cor:char osc}.   The latter is equivalent to~$(4g/\phi)^\circ\notin\bar{\omega}(H^\circ)$, by Lemma~\ref{lem:baromega comp conj} applied to $h$ in place of $f$.  Now~$\Upg(H^\circ)\cap\bar{\omega}(H^\circ)=\emptyset$ by Lemma~\ref{lem:Upgcapbaromega}, so it remains to note that~$4g\in\Upg(H)$ yields~$(4g/\phi)^\circ\in\Upg(H^\circ)$. 
\end{proof}

\section{Extending Hardy Fields to $\upo$-free Hardy Fields}\label{sec:upo-free Hardy fields}

\noindent
{\it In this section $H$ is a Hardy field.}\/
We first discuss in more detail the fundamental property of $\upo$-freeness and then prove a conjecture  from~\cite[\S{}17]{Boshernitzan82}. Next we establish the main result of the paper, Theorem~\ref{upo}, to the effect hat every maximal Hardy field is $\upo$-free.
We finish with some complements to the theorem.

\subsection*{Iterated logarithms and  $\upo$-freeness}
{\it In this subsection $H\supseteq\R(x)$ is log-closed.}\/
As in the introduction this yields a log-sequence $(\ell_{\rho})$ in $H^{>\R}$ from which we obtain sequences~$(\upg_\rho)$, $(\upl_\rho)$, $(\upo_{\rho})$ in $H$ as follows:
$$
\upg_\rho\ :=\ \ell_{\rho}^\dagger\in \Upg(H), \qquad
\upl_\rho\  :=\ -\upg_\rho^\dagger\in \Upl(H),\qquad \upo_{\rho}\ :=\ \omega(\upl_\rho)\in \omega\big(\Upl(H)\big).$$
Also $\upl_{\rho+1}=\upl_\rho+\upg_{\rho+1}$ and $\upo_{\rho+1}= \upo_\rho+ \upg_{\rho+1}^2$. The sequence~$\big(v(\upg_\rho)\big)$ is strictly increasing  and cofinal in
$\Psi_H$ by [ADH, beginning of 11.5], so the sequence $(\upg_{\rho})$ is strictly decreasing and coinitial in $\Upg(H)$. Using in addition that for all $f,g\in H^\times$ we have $f\prec g\Rightarrow f^\dagger<g^\dagger$, it follows that the pc-sequence $(\upl_\rho)$ is strictly increasing and cofinal in~$\Upl(H)$
 and the pc-sequence $(\upl_\rho+\upg_\rho)=\big({-(1/\ell_{\rho})'^{\dagger}}\big)$ is strictly decreasing and coinitial in~$\Upd(H)$; cf.~[ADH, proof of 11.8.15].
Recall that $\omega\colon H \to H$ is strictly increasing on~$\Upl(H)$; so the pc-sequence  $(\upo_{\rho})$ is strictly increasing
and cofinal in~$\omega\big(\Upl(H)\big)$.  
Hence  using also Corollary~\ref{cor:omega(H) downward closed}, we obtain:

\begin{cor}\label{hstar} 
If $H$ is $\d$-perfect, then $(\upo_\rho)$ is cofinal in $\bar{\omega}(H)$.
\end{cor} 
 
\noindent
Note that Corollary~\ref{hstar} applies in particular to the $\d$-perfect hull of any Hardy field. 
We have~$\sigma(\upg_{\rho})=\upo_\rho+\upg_{\rho}^2$, and by [ADH, 11.8.29]  the sequence~$\big(\sigma(\upg_{\rho})\big)$ is strictly decreasing and coinitial in $\sigma\big(\Upg(H)\big)$.  
Thus by Corollary~\ref{omuplosc}: 

{\it  if $H$ is $\upo$-free, then~$(\upo_\rho)$ is cofinal in $\bar{\omega}(H)$ and~$\big(\sigma(\upg_\rho)\big)$ is
coinitial in~${H\setminus\bar{\omega}(H)}$.}\/

{\samepage
\begin{lemma} \label{lem:17.7 generalized}
Suppose $H$ is $\upo$-free and $f\in H$.  Then the
following are equivalent:
\begin{enumerate}
\item[\textup{(i)}]  $f \in \bar{\omega}(H)$; 
\item[\textup{(ii)}]  $f < \upo_\rho$ for some $\rho$; 
\item[\textup{(iii)}] $f < \upo_\rho+c\upg_\rho^2$ for all $c\in\R^>$ and all $\rho$;
\item[\textup{(iv)}] there exists $c\in\R^>$ such that for all $\rho$ we have $f < \upo_\rho+c\upg_\rho^2$.
\end{enumerate}
\end{lemma}}
\begin{proof}
The equivalence (i)~$\Leftrightarrow$~(ii) holds by the sentence preceding the lemma.
The implication~(ii)~$\Rightarrow$~(iii) follows from
[ADH, 11.8.22], and~(iii)~$\Rightarrow$~(iv) is obvious.  
Since~$0<\upg_{\rho+1}\prec\upg_\rho$ we obtain for~$c\in\R^>$:  
$$\upo_{\rho+1}+c\upg_{\rho+1}^2\ =\  \upo_{\rho}+\upg_{\rho+1}^2+c\upg_{\rho+1}^2\ <\ \upo_\rho + \upg_\rho^2\ =\  \sigma(\upg_\rho).$$
This yields~(iv)~$\Rightarrow$~(i) by the sentence before the lemma.  
\end{proof}

{\samepage
\begin{cor}\label{bstar}  The following are equivalent:
\begin{enumerate}
\item[\textup{(i)}]
$H$ is $\upo$-free;
\item[\textup{(ii)}] $(\upo_\rho)$ is cofinal in $\bar{\omega}(H)$ and $\big(\sigma(\upg_\rho)\big)$ is
coinitial in~$H\setminus\bar{\omega}(H)$.
\end{enumerate}
\end{cor}}
\begin{proof}
The sentence preceding Lemma~\ref{lem:17.7 generalized} yields (i)~$\Rightarrow$~(ii). Now
assume (ii). As~$(\upo_{\rho})$ is cofinal in $\omega\big(\Upl(H)\big)$ and $\big(\sigma(\upg_\rho)\big)$ is
coinitial in $\sigma\big(\Upg(H)\big)$, we obtain
$$\bar{\omega}(H)\ =\ \omega\big(\Upl(H)\big){}^\downarrow, \qquad 
H\setminus \bar{\omega}(H)\ =\  \sigma\big(\Upg(H)\big){}^\uparrow,$$ 
hence
there is no $\upo\in H$ with~$\omega\big(\Upl(H)\big)< \upo < \sigma\big(\Upg(H)\big)$.
Thus if $H$ has asymptotic integration, then~$H$ is $\upo$-free by the remarks following Lemma~\ref{lem:Upgcapbaromega}. 
Suppose $H$ does not have asymptotic integration. As $H$ is log-closed, it is ungrounded, hence we have a gap~$v\upg$ in $H$, $\upg\in H^\times$. Then~$\upl_\rho\leadsto\upl:=-\upg^\dagger$ by [ADH, remark after~11.5.9], hence~$\upo_\rho\leadsto \upo:=\omega(\upl)$ by [ADH, 11.7.3], so 
$\upo_\rho<\upo\in\bar{\omega}(H)$ for all $\rho$,  and this contradicts (ii).
\end{proof}

\noindent
Note that $(\upo_\rho)$ is cofinal in $\bar{\omega}(H)$ iff for all $f\in H$ the equivalence~\eqref{eq:Hartman H} in
the introduction holds, and that $\big(\sigma(\upg_\rho)\big)$ is
coinitial in~$H\setminus\bar{\omega}(H)$ iff for all $f\in H$ the equivalence~\eqref{eq:Boshernitzan H} holds.
This justifies a remark before~\eqref{eq:Boshernitzan H}.
Lemma~\ref{lem:17.7 generalized} also yields a generalization of Hartman's \cite[Lemma~1]{Hartman48}:

\begin{cor}
Suppose $H$ is $\upo$-free, and $f\in H$. Then there exists $\rho$ such that~${f-\upo_{\rho'}\sim f-\upo_{\rho}}$ for all $\rho'\geq\rho$; for such $\rho$ we have $f<\upo_{\rho}$ iff $f\in\bar{\omega}(H)$.
\end{cor}
\begin{proof}
If $f\notin\bar{\omega}(H)$, then (i)$\Rightarrow$(iv) in Lemma~\ref{lem:17.7 generalized} yields $\rho$ such that
 $f-\upo_{\rho}\geq \upg_{\rho}^2$ and so~${f-\upo_{\rho}} \succ \upg_{\rho+1}^2$,
 hence  $f-\upo_{\rho} \succ  \upo_{\rho'}-\upo_{\rho}$ for $\rho'>\rho$, and thus~$f-\upo_{\rho'}\sim f-\upo_{\rho}$ 
for $\rho'>\rho$, as required. Suppose $f\in\bar{\omega}(H)$, and take an index $\sigma$ such that $f\leq \upo_{\sigma}$, and put $\rho:=\sigma+1$.
Then~$\upo_{\rho}-f\geq\upo_{\rho}-\upo_{\sigma}>0$ and thus for all $\rho'>\rho$ we have~$\upo_{\rho}-f\succeq \upg_{\rho}^2\succ\upg_{\rho+1}^2\asymp \upo_{\rho'}-\upo_{\rho}$
and so $\upo_{\rho}-f\sim \upo_{\rho'}-f$.
\end{proof}

\noindent
In the proof of Theorem~\ref{upo} we shall use:

\begin{lemma}\label{lem:upo}
Let $\upg\in(\Go)^\times$, $\upg > 0$, and $\upl:=-\upg^\dagger$ with $\upl_\rho < \upl <\upl_\rho+\upg_\rho$ in $\c$, for all $\rho$.
Then $\upg_\rho > \upg > \upg_{\rho}/\ell_{\rho}=(-1/\ell_\rho)'$ in $\c$, for all $\rho$.
\end{lemma}
\begin{proof}
Pick $a\in \R$ (independent of $\rho$) and functions in $\c_a$ whose germs at $+\infty$ are the elements
$\ell_{\rho}$,~$\upg_{\rho}$,~$\upl_{\rho}$ of $H$; denote these functions also by $\ell_{\rho}$,~$\upl_{\rho}$,~$\upg_{\rho}$. From $\ell_{\rho}^\dagger=\upg_{\rho}$ and $\upg_{\rho}^\dagger=-\upl_{\rho}$
 in~$H$ we obtain $c_{\rho}, d_{\rho}\in \R^{>}$ such that for all
 sufficiently large $t\ge a$, 
$$ \ell_{\rho}(t)\ =\ c_{\rho}\exp\!\left[\int_a^t \upg_{\rho}(s)\,ds\right], \quad \upg_{\rho}(t)\ =\ d_{\rho}\exp\!\left[-\int_a^t\upl_{\rho}(s)\,ds\right]. $$
(How large is ``sufficiently large'' depends on $\rho$.) Likewise we pick functions in $\c_a$ whose germ at $+\infty$ are $\upg$, $\upl$, and also denote these functions by $\upg$, $\upl$. From $\upg^\dagger=-\upl$ in $H$ we obtain a real constant $d>0$ such that for all sufficiently large $t\ge a$, 
$$ \upg(t)\ =\ d\exp\!\left[-\int_a^t\upl(s)\,ds\right].$$ 
Also, $\upl_\rho < \upl <\upl_\rho+\upg_\rho$ yields  constants $a_{\rho}, b_{\rho}\in \R$ such that for all $t\ge a$
$$\int_a^t \upl_{\rho}(s)\,ds\ <\ a_{\rho}+ \int_a^t \upl(s)\,ds\ <\ b_{\rho} +\int_a^t \upl_{\rho}(s)\,ds + \int_a^t \upg_{\rho}(s)\,ds,$$
 which by applying $\exp(-*)$ yields that for all sufficiently large $t\ge a$, 
$$\frac{1}{d_{\rho}}\upg_{\rho}(t)\ >\ \frac{1}{\ex^{a_{\rho}}d}\upg(t)\ >\ \frac{c_{\rho}}{\ex^{b_{\rho}}d_{\rho}} \upg_{\rho}(t)/\ell_{\rho}(t).$$
Here the positive constant factors don't matter,
since the valuation of $\upg_{\rho}$ is strictly increasing and that of $\upg_{\rho}/\ell_{\rho}=(-1/\ell_{\rho})'$
is strictly decreasing with $\rho$. Thus 
for all~$\rho$ we have 
$\upg_{\rho} > \upg > \upg_{\rho}/\ell_{\rho}= (-1/\ell_{\rho})'$, in $\c$. 
\end{proof}

\subsection*{Constructing $\upo$-free Hardy field extensions: special cases}
If $H$ is $\upo$-free, then so is every $\d$-algebraic Hardy field extension of~$H$, by [ADH, 13.6.1]. Thus if~$H$ is $\upo$-free, then
the Hardy-Liouville closure~$\Li\!\big(H(\R)\big)$ of $H(\R)$ is $\upo$-free.
Moreover, by \cite[Lemma~1.3.20]{ADH4}:

\begin{prop}\label{prop:1.3.20}
If $H$ is not $\upl$-free, then $\Li\!\big(H(\R)\big)$ is $\upo$-free.
\end{prop}

\noindent
Together with Corollary~\ref{cor:Ros83}, this yields:

\begin{cor}\label{cor:1.3.20}
Suppose $H$ is grounded. Then $L:=\Li\!\big(H(\R)\big)$ is $\upo$-free. Moreover, if for some $m$ we have $h\succ\ell_m$ for all $h\in H$
with $h\succ 1$, then for all $g\in L$ with~$g\succ 1$ there is an $n$ with $g\succ\ell_n$.
\end{cor}

\noindent
In  particular, $\Li(\R)=\Li\!\big(\R(x)\big)$ is $\upo$-free, and    $(\ell_n)$ is
coinitial in~$\Li(\R)^{>\R}$. We now use
these observations  to establish~\cite[Conjecture~17.11]{Boshernitzan82}: Corollary~\ref{cor:17.11}.
We first note that for $c\in\R$,   the germ
$(\upo_n+c\upg_{n}^2)/4$ 
generates oscillation iff~$c>0$. (See also the introduction.)
This follows  from the next corollary applied to~$f=\upo_n+c\upg_n^2$ and
 the grounded Hardy subfield~$H:=\R\langle\ell_n\rangle=\R(\ell_0,\dots,\ell_n)$ of~$\Li(\R)$:

 \begin{cor}\label{cor:17.7 generalized}
Suppose $H$ is grounded  and for some $m$ we have~${h\succ \ell_m}$ for all~${h\in H}$ with $h\succ 1$. Then for   $f\in H$, the
following are equivalent:
\begin{enumerate}
\item[\textup{(i)}] $f\in\bar{\omega}(H)$; 
\item[\textup{(ii)}]  $f < \upo_n$ for some $n$; 
\item[\textup{(iii)}] $f < \upo_n+c\upg_n^2$ for all $n$ and all  $c\in\R^>$;
\item[\textup{(iv)}] there exists $c\in\R^>$ such that for all $n$ we have $f < \upo_n+c\upg_n^2$.
\end{enumerate}
\end{cor}
\begin{proof}
By Corollary~\ref{cor:1.3.20}, $L:=\Li\!\big(H(\R)\big)$ is $\upo$-free, and  for all $g\in L$ with $g\succ 1$ there is an $n$ such that $\ell_n \prec g$.
Hence the corollary follows from Lemma~\ref{lem:17.7 generalized} applied to $L$ in place of $H$,
using $\bar{\omega}(H)=H\cap\bar{\omega}(L)$.
\end{proof}

\noindent
From the above equivalence (i)~$\Leftrightarrow$~(ii)  we recover~\cite[Theorem~17.7]{Boshernitzan82}:

\begin{cor}\label{cor:17.7}
Suppose $f\in\c$ is hardian and $\d$-algebraic over $\R$. Then
$$\text{$f$ generates oscillation}\ \Longleftrightarrow\  \text{$f>\upo_n/4$  for all $n$.}$$
\end{cor}
\begin{proof} By Corollary~\ref{cor:Bosh 14.11} the Hardy field~$H:=\R(x)\langle f\rangle$ satisfies the hypotheses of Corollary~\ref{cor:17.7 generalized}. Also, $f$ generates oscillation iff $4f\notin \bar{\omega}(H)$.  Now the equivalence follows from
(i)~$\Leftrightarrow$~(ii) in Corollary~\ref{cor:17.7 generalized}.
\end{proof}

\noindent
Using the above implication (iv)~$\Rightarrow$~(i) we obtain in the same way: 

\begin{cor}\label{cor:17.11}
Let $f\in\c$  be hardian and $\d$-algebraic over $\R$, and suppose for some~$c\in\R^>$ we have $f<\upo_n+c\upg_n^2$ for all $n$. Then $f/4$ does not generate oscillation. 
\end{cor}

\noindent
Next a variant of Proposition~\ref{prop:1.3.20}. Let $L\supseteq \R$ be a Liouville closed $\d$-algebraic Hardy field extension of $H$ such that $\omega(L)=\bar{\omega}(L)$.
(By Corollary~\ref{cor:omega(H) downward closed} this holds for~$L=\Dx(H)$.)
 Note that then $\bar{\omega}(L)=\omega\big(\Upl(L)\big)$ by [ADH, 11.8.20].
 
\begin{lemma}\label{lem:D(H) upo-free, 4} 
If $H$ is not $\upl$-free or $\bar{\omega}(H)=H\setminus\sigma\big(\Upg(H)\big){}^\uparrow$, then~$L$ is $\upo$-free.
\end{lemma}
\begin{proof}
If $H$ is $\upo$-free, then $L$ is $\upo$-free by [ADH, 13.6.1]. Hence, if $H$ is not $\upl$-free, then $L$ is 
$\upo$-free by Proposition~\ref{prop:1.3.20}. 
Suppose
$H$ is $\upl$-free but not $\upo$-free, and $\bar{\omega}(H)=H\setminus\sigma\big(\Upg(H)\big){}^\uparrow$. Then [ADH, 11.8.30] gives $\upo\in H$ with~$\omega\big(\Upl(H)\big) <  \upo  < \sigma\big(\Upg(H)\big)$. Hence 
$\upo\in\bar{\omega}(H)\subseteq \bar{\omega}(L)=\omega\big(\Upl(L)\big)$. Thus $L$ is $\upo$-free by 
\cite[Corollary~1.3.21]{ADH4}.
\end{proof}

\noindent
In \cite{ADH6} we show that for~$L=\Dx(H)$ the converse of Lemma~\ref{lem:D(H) upo-free, 4}   holds. 
Here are examples of (1) a non-$\upo$-free Liouville closed $\Gom$-Hardy field $L\supseteq \R(x)$, and (2) a non-$\upl$-free log-closed $\Gom$-Hardy field $M\supseteq \R(x)$ with asymptotic integration:

\begin{exampleNumbered}\label{ex:non-upo-free}
Remarks after Corollary~\ref{cor:Bosh 14.11} yield a translogarithmic and hardian germ $\ell\in\Gom$.
Then $E:=\R(\ell_0,\ell_1,\ell_2,\dots)$ is a Hardy subfield of $\Li\!\big(\R\langle\ell\rangle\big)$. Now $E$ is $\upo$-free by [ADH, 11.7.15], $\ell$ is $E$-hardian, and $\R<\ell<E^{>\R}$ by Corollary~\ref{cor:1.3.20}. 
Set
$$\upg\ :=\ \ell^\dagger,\qquad \upl\ :=\ -\upg^\dagger, \qquad \upo\ :=\ \omega(\upl),\qquad H\ :=\ E\langle\upo\rangle,
\qquad L\ :=\ \Li(H).$$
We have $\upl_n\leadsto\upl$ by [ADH, 11.5.7],  hence $\upo_n\leadsto\upo$ by [ADH,  11.7.3]. Then~$H$ is
an immediate $\upl$-free extension of $E$ by [ADH, 13.6.3, 13.6.4], and  $H^{>\R}$ is coinitial in~$L^{>\R}$
by   \cite[Proposition~1.3.15]{ADH4}. Now $(\ell_n)$ is coinitial in $E^{>\R}$, hence in $H^{>\R}$, and thus in $L^{>\R}$. It follows that $L$ is not $\upo$-free. Also $L\subseteq \Li\!\big(\R\langle\ell\rangle\big)\subseteq \Gom$. 
\end{exampleNumbered}

\begin{exampleNumbered}\label{ex:non-upl-free} Let $E$, etc.~be as in the previous example. Then
$E$ has asymptotic integration, and $E\<\upl\>$ is an immediate extension of $E$ by [ADH, 13.6.3], so $E\<\upl\>$ has asymptotic integration. Let $M$ be the smallest Hardy field extension of
$E\<\upl\>$ that is henselian as a valued field and closed under integration (hence log-closed). 
Then~$M\subseteq \Gom$, $M$ is an immediate extension of $E\<\upl\>$ by [ADH, 10.2.7], so has asymptotic integration, and is not $\upl$-free in view of $\upl_n\leadsto \upl$.  
\end{exampleNumbered} 

\subsection*{Proof of the main theorem} Here now is our main result:  

\begin{theorem}\label{upo} Every Hardy field has a $\d$-algebraic $\upo$-free Hardy field extension.
\end{theorem}

\begin{proof} 
It is enough to show that every $\d$-maximal Hardy field is $\upo$-free. That reduces to showing that
every non-$\upo$-free Liouville closed Hardy field containing $\R$
has a proper $\d$-algebraic Hardy field extension. So assume $H\supseteq \R$ is   Liouville closed  and  not $\upo$-free. We shall construct a proper $\d$-algebraic Hardy field extension
of~$H$.  As indicated in the remarks after Lemma~\ref{lem:Upgcapbaromega}, we have $\upo\in H$ such that 
$$\omega(H)\ <\  \upo\  <\ \sigma\big(\Upg(H)\big).$$ 
With $\upo$ in the role of $f$ in the discussion following Corollary~\ref{cor:char osc, 3}, we have  $\R$-linearly independent
solutions $y_1, y_2\in \Gt$ of the differential equation~${4Y'' + \upo Y=0}$; in fact, $y_1,y_2\in\Gi$.
Then the complex solution $y=y_1+y_2\imag$ is a unit
of~$\Gi[\imag]$, and so we have 
$z:=2y^\dagger\in \Gi[\imag]$.  We shall prove that the elements $\upl:=\Re z$ and~$\upg:= \Im z$ of $\Gi$ generate a proper $\d$-algebraic Hardy field extension $K=H(\upl, \upg)$ of~$H$
with~$\upo=\sigma(\upg)\in \sigma(K^\times)$. 
We can assume that~$w:= y_1y_2'-y_1'y_2\in \R^{>}$, so~$\upg=2w/|y|^2\in(\Gi)^\times$ and $\upg>0$.

Choose a log-sequence $(\ell_\rho)$ in $H$ and define $(\upg_\rho)$, $(\upl_\rho)$, $(\upo_\rho)$ as indicated
at the beginning of this section.
Then $\upo_{\rho} \leadsto \upo$, with~${\upo-\upo_{\rho}\sim \upg_{\rho+1}^2}$ by~[ADH,~11.7.1]. 
We aim to show:
\begin{equation}\label{eq:goal}
\upl-\upl_\rho\prec\upg_\rho\text{ and }\upg\prec\upg_\rho \text{ for all $\rho$.}
\end{equation}
For now  we fix $\rho$ and set
$g_{\rho}:=\upg_{\rho}^{-1/2}$, so $2g_{\rho}^\dagger=\upl_{\rho}=-\upg_{\rho}^\dagger$. 
For~$h\in H^\times$ we also have~$\omega(2h^\dagger)=-4h''/h$,
hence~$P:= 4Y'' + \upo Y\in H\{Y\}$ gives
$$P(g_{\rho})\ =\  
g_{\rho}(\upo-\upo_{\rho})\ \sim\ g_{\rho}\upg_{\rho+1}^2,$$ and so with an eye towards using Lemma~\ref{chvar}:
$$g_{\rho}^3P(g_{\rho})\ \sim\ g_{\rho}^4\upg_{\rho+1}^2\ \sim\ \upg_{\rho+1}^2/\upg_{\rho}^2\ \asymp\ 1/\ell_{\rho+1}^2.$$
Thus with
$g:= g_{\rho}=\upg_{\rho}^{-1/2}$, $\phi:=g^{-2}=\upg_{\rho}$ we  have $A_{\rho}\in \R^{>}$ such that
\begin{equation}\label{eq:bound}
g^3P_{\times g}^{\phi}(Y)\ =\ 4Y'' + g^3P(g)Y,\quad 
|g^3P(g)|\ \le\  A_{\rho}/\ell_{\rho+1}^2.
\end{equation}
From  $P(y)=0$ we get $P^{\phi}_{\times g}(y/g)=0$, that is,
$y/g\in \Gi[\imag]^\phi$ is a solution of $$4Y'' + g^3P(g)Y\ =\ 0,\ \text{ with }g^3P(g)\in H\subseteq \Gi.$$ 
Set $\ell:=\ell_{\rho+1}$, so~${\ell'=\ell_{\rho}^\dagger=\phi}$.
The subsection on compositional conjugation in Section~\ref{sec:Hardy fields}   yields the isomorphism
$$h\mapsto h^\circ=h\circ\ell^{\inv}\colon H^{\phi} \to H^\circ$$ of $H$-fields, where~$\ell^{\inv}$ is the compositional inverse of~$\ell$. Under this isomorphism the equation $4Y''+g^3P(g)Y=0$ corresponds to the equation
$$4Y'' + f_{\rho}Y\ =\ 0, \qquad f_{\rho}\ :=\ \big(g^3P(g)\big){}^\circ\in H^\circ\ \subseteq\ \Gi.$$ 
By Corollary~\ref{cor:chvar}, the equation $4Y'' + f_{\rho}Y=0$ has the ``real'' solutions
$$y_{j,\rho}\ :=\ (y_j/g)^\circ\in (\Gi)^\circ \ =\ \Gi \qquad(j=1,2),$$  and the ``complex''
solution $$y_\rho\ :=\ y_{1,\rho} + y_{2,\rho}\imag\ =\ (y/g)^\circ,$$ which is a unit of the ring
$\Gi[\imag]$. Set $z_{\rho}:=2 y_{\rho}^\dagger\in \Gi[\imag]$. 
The bound in \eqref{eq:bound} gives 
$|f_{\rho}|\ \le\  A_{\rho}/x^2$,
which by Corollary~\ref{cor:bound} yields  positive constants $B_{\rho}$, $c_\rho$ such that
$|z_{\rho}|\ \le\ B_{\rho}x^{c_\rho}$.  Using 
$(f^\circ)'=(\phi^{-1}f')^\circ$ for $f\in\Calinf[\imag]$ we obtain
$$z_{\rho}\ =\ 2\big((y/g)^\circ\big){}^\dagger \ =\ 2\big(\phi^{-1} (y/g)^\dagger\big)^\circ \ =\ \big((z-2g^\dagger)/\phi\big)^\circ$$
In combination with the  bound on $|z_\rho|$ this yields
\begin{align*} \left|\frac{z-2g^\dagger}{\phi}\right|\ &\le\ B_{\rho}\,\ell_{\rho+1}^{c_\rho}, \quad \text{hence}\\
|z- \upl_{\rho}|\ &\le\ B_{\rho}\,\ell_{\rho+1}^{c_\rho}\,\phi\ =\ 
  B_{\rho}\,\ell_{\rho+1}^{c_\rho}\,\upg_{\rho}, \quad \text{and so} \\ 
z\ &=\ \upl_{\rho}  + R_{\rho}\quad\text{where}\quad |R_{\rho}|\le B_{\rho}\,\ell_{\rho+1}^{c_\rho} \,\upg_{\rho}. 
\end{align*}
We now use this last estimate with $\rho+1$ instead of $\rho$, 
together with 
$$\upl_{\rho+1}\ =\ \upl_{\rho}+\upg_{\rho+1},\quad  \ell_{\rho+1}\upg_{\rho+1}\ =\ \upg_{\rho}.$$
This yields
\begin{align*} z &\, =\, \upl_{\rho} + \upg_{\rho+1} + R_{\rho+1}\\ & \qquad \text{with}\ 
 |R_{\rho+1}|\, \le\, B_{\rho+1}\,\ell_{\rho+2}^{c_{\rho+1}}\, \upg_{\rho+1} \, =\, B_{\rho+1}\big(\ell_{\rho+2}^{c_{\rho+1}}/\ell_{\rho+1}\big)\,\upg_{\rho},\\
\text{so}\quad z &\, =\, \upl_{\rho} + o(\upg_{\rho})\ \text{ that is, }z-\upl_{\rho}\prec \upg_{\rho}, \\
\text{and thus}\quad \upl\, &=\, \Re z\, =\, \upl_{\rho} + o(\upg_{\rho}), \quad \upg\, =\, \Im z\, \prec\, \upg_{\rho}, \quad\text{proving \eqref{eq:goal}.}
\end{align*}
Now varying $\rho$ again, $(\upl_{\rho})$ is a strictly increasing divergent
pc-sequence in $H$ which is cofinal in $\Upl(H)$, and 
$(\upl_\rho+\upg_\rho)$ is a strictly decreasing pc-sequence in $H$ which is coinitial in
$\Upd(H)=H\setminus\Upl(H)$. 
By the above, for each $\rho$ we have
$\upl=\upl_{\rho+1}+o(\upg_{\rho+1})=\upl_\rho+\upg_{\rho+1}+o(\upg_{\rho+1})$  
and hence $\upl_\rho<\upl<\upl_{\rho+1}+\upg_{\rho+1}$, thus
$\upl= \Re z$ satisfies $\Upl(H) < \upl < \Upd(H)$. This yields an ordered subfield 
$H(\upl)$ of $\Gi$, which by Lemma~\ref{ps1} is an immediate
valued field extension of $H$ with $\upl_{\rho} \leadsto \upl$.
Now~$\upl=-\upg^\dagger$ (see discussion before Lemma~\ref{lem:Im z}), so Lemma~\ref{lem:upo} gives~$\upg_{\rho} > \upg > (-1/\ell_{\rho})'$ in $\Gi$, for all $\rho$. In view of Lemma~\ref{ps2} applied to $H(\upl)$,~$\upg$ in the role of~$H$,~$f$ this yields an ordered subfield $H(\upl, \upg)$
of $\Gi$ where $\upg$ is transcendental over $H(\upl)$. Moreover,~$\upg$ 
satisfies the second-order differential equation
$2yy''-3(y')^2+y^4-\upo y^2=0$ over~$H$ (obtained from the relation 
$\sigma(\upg)=\upo$ by multiplication with $\upg^2$). It follows that~$H(\upl,\upg)$ is closed under the derivation of~$\Gi$, and hence~$H(\upl, \upg)=H\<\upg\>$ is a Hardy field that is $\d$-algebraic over $H$. 
\end{proof} 

\noindent
The proof also shows that every $\Ginf$-Hardy field
has an 
$\upo$-free $\d$-algebraic $\Ginf$-Hardy field extension, and the same with 
$\Gom$ instead of $\Gi$. In \cite{ADH6} we prove that the perfect hull of an $\upo$-free Hardy field  
is $\upo$-free, but that not every perfect Hardy field is $\upo$-free. 


\subsection*{Improving Theorem~\ref{upo}} 
{\it In this subsection, assume $H\supseteq\R$ is  Liouville closed  
and $\upo\in H$, $\upg\in (\c^2)^\times$ satisfy
$\omega(H) <  \upo < \sigma\big(\Upg(H)\big)$ and $\sigma(\upg)=\upo$.}\/
Proposition~\ref{prop:1.3.20} and results from [ADH, 13.7] lead to a more explicit version of Theorem~\ref{upo}: 

\begin{cor}\label{cor:upo}
The germ  $\upg$ generates a Hardy field extension~$H\langle\upg\rangle$ of $H$ with a gap $v\upg$,
and so $\operatorname{Li}\!\big(H\langle\upg\rangle\big)$ is an $\upo$-free Hardy field extension of $H$.
\end{cor}
\begin{proof}
Since $\sigma(-\upg)=\sigma(\upg)$, we may arrange $\upg>0$. The discussion before Lem\-ma~\ref{lem:Im z} with $\upo$, $\upg$ in the roles of $f$, $u$, respectively, yields
$\R$-linearly independent solutions $y_1,y_2\in\Calinf$ of the differential equation $4Y''+\upo Y=0$ with
Wronskian~$1/2$ such that~$\upg=1/(y_1^2+y_2^2)$.  
The  proof of Theorem~\ref{upo} shows that~$\upg$ generates a Hardy field extension~$H\langle\upg\rangle=H(\upl, \upg)$ of $H$. Recall that $v(\upg_{\rho})$ is strictly increasing as a function of $\rho$ and cofinal in $\Psi_H$; as $\upg \prec \upg_{\rho}$ for all $\rho$, this gives~$\Psi_H<v\upg$. Also~$\upg> (-1/\ell_{\rho})'>0$ for all $\rho$ and
$v(1/\ell_{\rho})'$ is strictly decreasing as a function of~$\rho$ and coinitial in $(\Gamma_H^>)'$, and so $v\upg < (\Gamma_H^{>})'$. Then
by [ADH, 13.7.1 and subsequent remark~(2) on p.~626], $v\upg$ is a gap in~$H\langle \upg\rangle$.
Thus $H\langle \upg\rangle$ does not have asymptotic integration and hence is not $\upl$-free,
so
 $\operatorname{Li}\!\big(H\langle\upg\rangle\big)$ is $\upo$-free by Proposition~\ref{prop:1.3.20}. 
\end{proof}

\begin{cor} 
Suppose $\upg>0$. Then with~$L:=\operatorname{Li}\!\big(H\langle\upg\rangle\big)$, 
$$\upo\notin\bar{\omega}(H)	    
\ \Longleftrightarrow\ \upg \in\Upg(L),  \qquad
\upo\in\bar{\omega}(H)	    
\ \Longleftrightarrow\ \upg \in\I(L).$$
\end{cor}
\begin{proof} If $\upg\notin \Upg(L)$, then $\upo\in \omega(L)^{\downarrow}$ by [ADH, 11.8.31], hence $\upo\in \bar{\omega}(H)$. If $\upg\in \Upg(L)$, then we can use Corollary~\ref{omuplosc} for $L$ to conclude $\upo\notin \bar{\omega}(H)$. The equivalence on the right now follows from that on the left and [ADH, 11.8.19]. 
\end{proof}

\noindent
We also note that if $\upo/4$ generates oscillation, then we have many choices for $\upg$:

\begin{cor}\label{cor:sigma(upg)=upo}
Suppose~$\upo/4$ generates oscillation. Then there are continuum many $\tilde\upg\in (\Calinf)^\times$
with $\tilde\upg>0$ and $\sigma(\tilde\upg)=\upo$; no Hardy field extension of~$H$ contains more than one such germ $\tilde\upg$.
\textup{(}Thus $H$ has at least continuum many maximal Hardy field extensions.\textup{)}
\end{cor} 

\begin{proof} As before we arrange $\upg>0$ and set $L:=\operatorname{Li}\!\big(H\langle\upg\rangle\big)$.
 Take~$\phi\in L$ with $\phi'= \frac{1}{2}\upg$ and consider the germs
$$y_1 :=  \frac{1}{\sqrt{\upg}}\cos \phi, \quad y_2 := \frac{1}{\sqrt{\upg}}\sin \phi\qquad\text{in $\Calinf$.}$$ 
The remarks preceding Lemma~\ref{lem:Im z} show that $y_1$, $y_2$ solve the differential equation~${4Y''+\upo Y=0}$, their
Wronskian equals $1/2$, and $\phi\succ 1$ (since $\upo/4$ generates oscillation). We now dilate~$y_1$,~$y_2$: let~$r\in\R^>$ be arbitrary and set
$$y_{1r}\ :=\ ry_1,\qquad   y_{2r}\ :=\ r^{-1}y_2.$$
Then $y_{1r}$, $y_{2r}$ still solve the   equation~$4Y''+\upo Y=0$, and their Wronskian is $1/2$.
Put~$\upg_r:=1/(y_{1r}^2+  y_{2r}^2)\in\Calinf$. Then $\sigma(\upg_r)=\upo$.
Let   $r,s\in\R^>$. Then
$$\upg_r=\upg_s \quad\Longleftrightarrow\quad y_{1r}^2+y_{2r}^2 = y_{1s}^2+y_{2s}^2
\quad\Longleftrightarrow\quad (r^2-s^2)\cos^2\phi + \big(\textstyle\frac{1}{r^2}-\frac{1}{s^2}\big)\sin^2\phi=0,$$
so $\upg_r=\upg_s$ iff $r=s$.
Next, suppose   $M\supseteq H$ is a $\d$-perfect Hardy field  containing both $\upg$ and $\tilde\upg \in (\Calinf)^\times$
with $\tilde\upg>0$ and $\sigma(\tilde\upg)=\upo$.
Corollary~\ref{cor:omega(H) downward closed} gives~${\upo\notin\omega(M)}$,
so
$\upg, \tilde\upg\in\Upg(M)$ by [ADH, 11.8.31],     hence~$\upg=\tilde\upg$ by~[ADH, 11.8.29]. 
\end{proof}

\section{A Question of Boshernitzan}\label{sec:applications}

\noindent
In this final section we apply Theorem~\ref{upo} to answer a question from~\cite{Boshernitzan86} and to generalize a theorem from \cite{Boshernitzan82}.  

\subsection*{Translogarithmic germs in maximal Hardy fields}
The following  analogue of Corollary~\ref{cor:Bosh 1.3} for translogarithmic germs  gives a positive answer to Question~4 in~\cite[\S{}7]{Boshernitzan86}: 

\begin{prop}\label{prop:translog}
Every maximal Hardy field contains a translogarithmic germ.
\end{prop}

\noindent
Let $H\supseteq\R$ be a Liouville closed Hardy field; then $H$ has no translogarithmic element iff $(\ell_n)$ is a logarithmic sequence for $H$ in the sense of [ADH, 11.5]. 
If in this case~$H$ is also $\upo$-free, then for each  translogarithmic $H$-hardian germ~$y$ the isomorphism type of the ordered  differential field
$H\langle y\rangle$ over $H$ is uniquely determined. This is part of the next lemma, which follows from [ADH, 13.6.7, 13.6.8]. We need to assume familiarity with the  {\em Newton degree $\ndeg(P)$}\/ and {\em Newton weight
$\nwt(P)$}\/ of~$P\in H\{Y\}^{\ne}$ for suitable $H$; see [ADH, 11.1]. 

\begin{lemma}
Let $H$ be an   $\upo$-free $H$-field, with asymptotic couple $(\Gamma,\psi)$, and let~${L=H\langle y\rangle}$ be a pre-$H$-field extension  of $H$ with $\Gamma^< < vy< 0$. Then for all~$P\in H\{Y\}^{\neq}$ we have
$$v\big(P(y)\big)\ =\ \gamma + \ndeg(P)vy+\nwt(P)\psi_{L}(vy)\qquad\text{where $\gamma=v^{\ev}(P)\in\Gamma$,}$$
and thus
$$\Gamma_L\ =\ \Gamma\oplus \Z vy \oplus \Z \psi_L(vy)\qquad \text{\textup{(}internal direct sum\textup{)}.}$$
Moreover, if $L^*=H\langle y^*\rangle$ is a  pre-$H$-field extension  of $H$ with $\Gamma^< < vy^*< 0$ and~$\sgn y=\sgn y^*$, then
there is a unique pre-$H$-field isomorphism $L\to L^*$ which is the identity on $H$ and sends $y$ to $y^*$.
\end{lemma}

\noindent
This lemma   suggests how to obtain Proposition~\ref{prop:translog}: follow the arguments in the proof of [ADH, 13.6.7]. In the rest of this subsection we carry out this plan.  
In the next lemma and corollary $H\supseteq\R$ is a Liouville closed Hardy field and
   $y\in\Calinf$. The proof of this lemma assumes familiarity with the binary relation  $\preceq^\flat$
   on $K$ and~$K\{Y\}$ for suitable $H$-asymptotic fields  $K$, and their variants; see [ADH, 9.4]. 
   

\begin{lemma}\label{lem:translog} 
Suppose $H$ is $\upo$-free and for all $\ell\in H^{>\R}$ we have, in $\c$:
\begin{enumerate}
\item[\textup{(i)}] $1\prec y\prec \ell$;
\item[\textup{(ii)}] $\derdelta^n(y)\preceq 1$ for all $n\ge 1$, where $\derdelta:=\phi^{-1}\der$, $\phi:=\ell'$;
\item[\textup{(iii)}] $y'\in\c^\times$ and $(1/\ell)' \preceq y^\dagger$.
\end{enumerate}
Let~$P\in H\{Y\}^{\neq}$. Then in~$\c$ we have
$$P(y) \sim a\,  y^d \,  (y^\dagger)^w \qquad\text{where
  $a\in H^\times$, $d=\ndeg(P)$, $w=\nwt(P)$.}$$
\textup{(}Hence $y$ is hardian over $H$ and $\d$-transcendental over $H$.\textup{)}
\end{lemma}

\begin{proof}
Since $H$ is real closed, it has a monomial group by [ADH, 3.3.32], so the material of [ADH, 13.3] applies. 
Then [ADH, 13.3.3] gives a monic~$D\in\R[Y]^{\neq}$, $b\in H^\times$, $w\in\N$, and an  active element $\phi$ of $H$ with $0<\phi\prec 1$ such that:
$$P^\phi\  =\ b \cdot D\cdot (Y')^w + R,\qquad   R\in H^{\phi}\{Y\},\ R  \prec_{\phi}^\flat  b.$$
Set $d:=\ndeg P$, and note that  by~[ADH, 13.1.9] we have $d=\deg D+w$ and~$w=\nwt P$.
Replace~$P$,~$b$,~$R$ by~$b^{-1}P$,~$1$,~$b^{-1}R$, respectively, to arrange $b=1$. Take~$\ell\in H$ with $\ell'=\phi$, so~$\ell>\R$; we use the superscript $\circ$ as in the subsection on compositional conjugation
of Section~\ref{sec:Hardy fields}; in particular, 
$y^\circ=y\circ\ell^{\operatorname{inv}}$ with~$(y^\circ)'=(\phi^{-1}y')^\circ$, so~$(y^\circ)^\dagger\succeq 1/x^2$
by   hypothesis (iii) of our lemma.
For~$f,g\in H$ we have
$$f\prec^\flat_\phi g \text{\ (in $H$)} \quad\Longleftrightarrow\quad 
f\prec^\flat g \text{\ (in $H^\phi$)} \quad\Longleftrightarrow\quad
f^\circ\prec^\flat g^\circ\text{\ (in $H^\circ$)}.$$
Hence
in $H^\circ\{Y\}$ we   have
$$(P^\phi)^\circ\  =\   D\cdot (Y')^w + R^\circ\qquad\text{where}\qquad    R^\circ  \prec^\flat 1.$$
Evaluating at $y^\circ$ we have $D(y^\circ)\big((y^\circ)'\big){}^w \sim (y^\circ)^d \big((y^\circ)^\dagger\big){}^w$ and so
$D(y^\circ)\big((y^\circ)'\big){}^w \succeq x^{-2w} \asymp^\flat 1$.
By (i) we have $(y^\circ)^m\prec x$ for $m\geq 1$, and
by~(ii)  we have
$(y^\circ)^{(n)}\preceq 1$ for $n\geq 1$.
Hence $R^\circ(y^\circ)\preceq  h^\circ$ for some  $h\in H$ with $h^\circ\prec^\flat 1$.
Thus  in $\c$ we have
$$(P^\phi)^\circ(y^\circ) \sim   (y^\circ)^d  \big((y^\circ)^\dagger\big){}^w.$$
Since $P(y)^\circ=(P^\phi)^\circ(y^\circ)$, this yields
$P(y)\sim a\cdot y^d\cdot (y^\dagger)^w$ for $a=\phi^{-w}$.
\end{proof}

\begin{cor}\label{cor:translog}
Suppose $H$ is $\upo$-free and $1\prec y\prec\ell$ for all $\ell\in H^{>\R}$.
Then the following are equivalent:
\begin{enumerate}
\item[\textup{(i)}] $y$ is hardian over $H$;
\item[\textup{(ii)}]  for all $h\in H^{>\R}$ there is an $\ell\in H^{>\R}$ such that $\ell\preceq h$ and
$y$, $\ell$ lie in a common Hardy field;
\item[\textup{(iii)}] for all $h\in H^{>\R}$ there is an $\ell\in H^{>\R}$ such that $\ell\preceq h$ and $y\circ\ell^{\operatorname{inv}}$ is hardian.
\end{enumerate}
\end{cor}
 
\begin{proof}
 (i)~$\Rightarrow$~(ii)~$\Rightarrow$~(iii) is clear. Let~$\ell\in H^{>\R}$ be such that~$y^\circ:=y\circ\ell^{\operatorname{inv}}$ lies in a Hardy field~$H_0$; we arrange $x\in H_0$. For $\phi:=\ell'$ we have 
 $$(\phi^{-1}y^\dagger)^\circ\  =\  
(y^\circ)^\dagger\  \succ\  (1/x)'\ =\ -1/x^2$$
and thus $y^\dagger \succ  -\phi/\ell^2=(1/\ell)'$.
Also~$y^\circ\prec x$, hence~$z:=(y^\circ)'\prec x'=1$ and so~$z^{(n)}\prec 1$ for all $n$.
With $\derdelta:=\phi^{-1}\der$ and $n\geq 1$
we   have~$\derdelta^{n}(y)^\circ=z^{(n-1)}$, and thus~$\derdelta^{n}(y)\prec 1$. 
 Moreover, for $h\in H^{>\R}$ with~$\ell\preceq h$ and $\theta:=h'$
 we have~$\theta^{-1}\der=f\derdelta$ where $f:=\phi/\theta\in H$, $f\preceq 1$.  Let $n\ge 1$. Then  
 \begin{multline*} (\theta^{-1}\der)^n=(f\derdelta)^n=G_n^n(f)\derdelta^n + \cdots + G^n_1(f)\derdelta\quad\text{ on $\Calinf$} \\ \text{ where
 $G^n_j\in \Q\{X\}\subseteq H^\phi\{X\}$ for $j=1,\dots,n$.}\end{multline*} 
As $\derdelta$ is small as a derivation on $H$, we have
 $G^n_j(f)\preceq 1$ for~$j=1,\dots,n$, and so~$(\theta^{-1}\der)^n(y) \prec  1$. 
 Thus (iii)~$\Rightarrow$~(i) by Lemma~\ref{lem:translog}.
 \end{proof}

\begin{proof}[Proof of Proposition~\ref{prop:translog}]
Let $H\supseteq\R$ be any $\upo$-free Liouville closed Hardy field not containing any translogarithmic element;
in view of Theorem~\ref{upo} it suffices to show that then some Hardy field extension of $H$ contains a translogarithmic element.  The remark 
after Corollary~\ref{cor:Bosh 14.11} yields a translogarithmic germ $y$ in a $\Gom$-Hardy field~${H_0\supseteq\R}$. 
Then for each~$n$, the germs~$y$,~$\ell_n$ are contained in a common Hardy field, namely $\Li(H_0)$.
Hence $y$ generates a proper Hardy field extension of~$H$ by~(ii)~$\Rightarrow$~(i) in Corollary~\ref{cor:translog}.
\end{proof}

\noindent
Proposition~\ref{prop:translog}  goes through
when ``maximal'' is replaced by ``$\Ginf$-maximal'' or ``$\Gom$-maximal''. This follows from its proof, using also
remarks after the proof of Theorem~\ref{upo}. Here is a conjecture that is much stronger than  Proposition~\ref{prop:translog}; it postulates an analogue of  Corollary~\ref{cor:Bosh 1.1} for infinite ``lower bounds'':

\begin{conjecture}
If $H$ is maximal, then there is no $y\in\c^1$ such that $1 \prec y \prec h$ for all~$h\in H^{>\R}$, and $y'\in\c^\times$.
\end{conjecture}

\noindent 
We observe that in this conjecture we may restrict attention to $\Gom$-hardian germs~$y$:

\begin{lemma}
Suppose  there exists
 $y\in\c^1$ such that $1\prec y\prec h$ for all $h\in H^{>\R}$ and~$y'\in\c^\times$. Then there exists such a germ~$y$ which is $\Gom$-hardian.
\end{lemma}
\begin{proof} Take $y$ as in the hypothesis.
Replace $y$ by $-y$ if necessary to arrange~$y>\R$.
Now Theorem~\ref{thm:Bosh 1.2} yields a $\Gom$-hardian germ  $z \geq y^{\operatorname{inv}}$. By 
Lemma~\ref{lem:Bosh6.5}, the germ~$z^{\operatorname{inv}}$ is also  
 $\Gom$-hardian, and~$\R< z^{\operatorname{inv}}\leq y\prec h$ for all $h\in H^{>\R}$.
\end{proof}

\subsection*{Lower bounds on the growth of germs in $\Ex(H)$}
{\it In this subsection $H$ is a Hardy field.}\/ Recall from Corollary~\ref{cor:Bosh13.10} that for all $f\in\Ex(H)$ there are~${h\in H(x)}$ and~$n$ such that $f\leq \exp_n h$.  In particular, the sequence $(\exp_n x)$ is cofinal in~$\Ex(\Q)$.
By Theorem~\ref{thm:Bosh 14.4} and Corollary~\ref{cor:Ros83}, the sequence $(\ell_n)=(\log_n x)$ is coinitial in~$\Ex(\Q)^{>\R}$; 
see also~\cite[Theorem~13.2]{Boshernitzan82}.
Thus for the Hardy field $H=\Li(\R)$,
the subset $H^{>\R}$ is coinitial  in  $\Ex(\Q)^{>\R}=\Ex(H)^{>\R}$, equivalently,
$\Gamma_H^<$ is  cofinal   in~$\Gamma_{\Ex(H)}^{<}$.  We now generalize this fact, also recalling from
the remarks after Corollary~\ref{cor:1.3.20} that~$\Li(\R)$ is $\upo$-free:  

\begin{theorem}\label{thm:coinitial in E(H)}
Suppose $H$ is $\upo$-free. Then $\Gamma_H^<$ is   cofinal   in $\Gamma_{\Ex(H)}^{<}$.
\end{theorem}

\begin{proof}
Replacing $H$ by $\Li\!\big(H(\R)\big)$ and using [ADH, 13.6.1]  we arrange that $H$ is
Liouville closed and $H\supseteq \R$.
Let $y\in\Ex(H)$ and  suppose towards a contradiction that~$\R < y < H^{>\R}$.
Then $f:=y^{\operatorname{inv}}$  is transexponential  and hardian (Lemma~\ref{lem:Bosh6.5}). 
Lemma~\ref{lem:bounded Hardy field ext, 3} gives a bound $b\in \c^\times$ for $\R\langle f \rangle$. 
Lem\-ma~\ref{lem:Bosh 14.3} gives~$\phi\in (\c^\omega)^\times$ such that~$\phi^{(n)}\prec 1/b$ for all~$n$; 
set $r:=\phi\cdot\sin x\in\c^\omega$. Then by Lem\-ma~\ref{lem:Q bound} (with~$\R\langle f\rangle$ in place of $H$)
we have $Q(r)\prec 1$ for all~$Q\in \R\langle f\rangle\{Y\}$ with~${Q(0)=0}$.
Hence~$g:=f+r$ is eventually strictly increasing with~$g\succ 1$, and~$y=f^{\operatorname{inv}}$ and~$z:=g^{\operatorname{inv}}\in\Calinf$ do not lie in a common Hardy field. 
Thus in order to achieve the desired contradiction   it suffices to show that $z$ 
is $H$-hardian.  For this we use Corollary~\ref{cor:translog}. It is clear that~$f\sim g$, so
$y\sim z$ by Corollary~\ref{cor:Entr}, and thus~$1\prec z\prec \ell$ for all~${\ell\in H^{>\R}}$.
Let~${\ell\in H^{>\R}}$ and $\ell \prec x$; we claim that~$z\circ\ell^{\operatorname{inv}}$ is hardian, equivalently, by Lemma~\ref{lem:Bosh6.5}, that~$\ell\circ g=(z\circ\ell^{\operatorname{inv}})^{\operatorname{inv}}$ is hardian.
Now~$\ell\circ f=(y\circ\ell^{\operatorname{inv}})^{\operatorname{inv}}$ is hardian and $\ell\circ f\succ 1$,
and Lemma~\ref{lem:difference in I} gives $\ell\circ f-\ell\circ g\in (\Calinf)^{\preceq}$. 
Hence~${\ell\circ f\sim_{\infty} \ell\circ g}$ by Lemma~\ref{lem:y siminf z}.  For all $n$ we have~${\ell_n\circ \ell=\log_n \ell\in H^{>\R}}$, so~$y\le \ell_n\circ \ell$, hence~${y\circ \ell^{\inv}\le \ell_n}$, which by
compositional inversion gives~${\ell\circ f\ge \exp_n x}$. So $\ell\circ g$ is hardian by Corollary~\ref{cor:yhardian}. 
Thus~$z$ is $H$-hardian by (iii)~$\Rightarrow$~(i) of Corollary~\ref{cor:translog}.
\end{proof}

\noindent
If $H\subseteq\Ginf$ is $\upo$-free, then  $\Gamma_H^<$ is also  cofinal   in $\Gamma_{\Ex^\infty(H)}^{<}$, and similarly
with $\omega$ in place of $\infty$. (Same proof as that of Theorem~\ref{thm:coinitial in E(H)}.) If $H$ is bounded, then $\Dx(H)=\Ex(H)$ by Theorem~\ref{thm:Bosh 14.4}, in which case
Theorem~\ref{thm:coinitial in E(H)} already follows from [ADH, 13.6.1]. 
Boshernitzan~\cite[p.~144]{Boshernitzan82} 
asked whether   $\Dx(H)=\Ex(H)$ in general, and
he gave Theorem~\ref{thm:Bosh 14.4} as support for a positive answer. Our
Theorem~\ref{thm:coinitial in E(H)}   can be seen as further evidence.

\newlength\templinewidth
\setlength{\templinewidth}{\textwidth}
\addtolength{\templinewidth}{-2.25em}

\patchcmd{\thebibliography}{\list}{\printremarkbeforebib\list}{}{}

\let\oldaddcontentsline\addcontentsline
\renewcommand{\addcontentsline}[3]{\oldaddcontentsline{toc}{section}{References}}

\def\printremarkbeforebib{\bigskip\hskip1em The citation [ADH] refers to our book \\

\hskip1em\parbox{\templinewidth}{
M. Aschenbrenner, L. van den Dries, J. van der Hoeven,
\textit{Asymptotic Differential Algebra and Model Theory of Transseries,} Annals of Mathematics Studies, vol.~195, Princeton University Press, Princeton, NJ, 2017.
}

\bigskip

} 

\bibliographystyle{amsplain}

\end{document}